\documentclass{amsart}
\usepackage[foot]{amsaddr}

\usepackage{graphicx}

\usepackage{multicol}
\usepackage[bottom]{footmisc}
\usepackage{bm}
\usepackage{cite}
\usepackage{url}

\usepackage[margin=3cm]{geometry}
\usepackage{lipsum}

\usepackage{float}
\usepackage{subfig}

\newcommand{\pa}{\partial}
\newcommand{\lt}{\left}
\newcommand{\rt}{\right}
\newcommand{\R}{\mathbb{R}}

\usepackage{mathtools}
\usepackage{graphicx}
\usepackage{multirow}
\usepackage[mathscr]{euscript}
\usepackage{enumitem}
\usepackage{changepage}
\usepackage{bbm}
\newtheorem{theorem}{Theorem}[section]

\theoremstyle{definition}

\theoremstyle{remark}
\newtheorem{remark}[theorem]{Remark}

\theoremstyle{corollary}

\usepackage{color}
\newcommand{\blue}[1]{\textcolor{black}{#1}}
\newcommand{\red}[1]{\textcolor{black}{#1}}

\newtheoremstyle{examplestyle}
   {}{}{}{}{\bfseries}{.}{.5em}{{\thmname{#1 }}{\thmnumber{#2}}{\thmnote{ (#3)}}}
\theoremstyle{examplestyle}
\newtheorem{examplecase}{Example}[section]

\makeindex             

\begin{document}

\title[  ]{Well-balanced finite volume schemes for hydrodynamic equations with general free energy}

\author{Jos\'e A. Carrillo}
\address[Jos\'e A. Carrillo]{Department of Mathematics, Imperial College London, SW7 2AZ, UK}
\curraddr{}
\email{carrillo@imperial.ac.uk}
\thanks{}

\author{Serafim Kalliadasis}
\address[Serafim Kalliadasis]{Department of Chemical Engineering, Imperial College London, SW7 2AZ, UK}
\curraddr{}
\email{s.kalliadasis@imperial.ac.uk}
\thanks{}

\author{Sergio P. Perez}
\address[Sergio P. Perez]{Departments of Chemical Engineering and Mathematics, Imperial College London, SW7 2AZ, UK}
\curraddr{}
\email{sergio.perez15@imperial.ac.uk}
\thanks{}

\author{Chi-Wang Shu}
\address[Chi-Wang Shu]{Division of Applied Mathematics, Brown University, Providence, RI 02912, USA}
\curraddr{}
\email{shu@dam.brown.edu}
\thanks{}

%

\begin{abstract}
Well-balanced and free energy dissipative first- and second-order accurate
finite volume schemes are proposed for a general class of hydrodynamic
systems with linear and nonlinear damping. The variation of the natural
Liapunov functional of the system, given by its free energy, allows, for a
characterization of the stationary states by its variation. An analog
property at the discrete level enables us to preserve stationary states at
machine precision while keeping the dissipation of the discrete free
energy. These schemes can accurately analyse the stability properties of
stationary states in challenging problems such as: phase transitions in
collective behavior, generalized Euler-Poisson systems in chemotaxis and
astrophysics, and models in dynamic density functional theories; having
done a careful validation in a battery of relevant test cases.
\end{abstract}

\maketitle


%
%
%
%

\section{Introduction}\label{sec:intro}

The construction of robust well-balanced numerical methods for conservation
laws has attracted a lot of attention since the initial works of LeRoux and
collaborators \cite{greenberg1996well, gosse1996well}. The well-balanced
property is equivalent to the exact C-property defined beforehand by
Berm\'udez and V\'azquez in \cite{bermudez1994upwind}, and both of them refer
to the ability of a numerical scheme to preserve the steady states at a
discrete level and to accurately compute evolutions of small deviations from
them. The historical evolution of well-balanced schemes is reviewed in \cite{gosse2013computing}. On the other hand, the derivation of numerical schemes preserving
structural properties of the evolutions under study such as dissipations or
conservations of relevant physical quantities is an important line of
research in hydrodynamic systems and their overdamped limits, see for
instance
\cite{de2012discontinuous,carrillo2015finite,sun2018discontinuous,pareschi2018structure}.
In the present work, we propose numerical schemes with well-balanced and free
energy dissipation properties for a general class of balance laws or
hydrodynamic models with attractive-repulsive interaction forces, and linear
or nonlinear damping effects, such as the Cucker-Smale alignment term in
swarming. The general hydrodynamic system has the form
\begin{equation}\label{eq:generalsys}
\begin{cases}
\partial_{t}\rho+\nabla\cdot\left(\rho \bm{u}\right)=0,\quad \bm{x}\in\R^{d},\quad t>0,\\[2mm]
{\displaystyle \pa_{t}(\rho \bm{u})\!+\!\nabla\!\cdot\!(\rho \bm{u}\otimes \bm{u})\!=-\nabla P(\rho)-\rho \nabla H(\bm{x},\rho) - \gamma\rho \bm{u}-\!\rho\!\!\!\int_{\R^{d}}\!\!\psi(\bm{x}-\bm{y})(\bm{u}(\bm{x})-\bm{u}(\bm{y}))\rho(\bm{y})\,d\bm{y}
,}
\end{cases}
\end{equation}
where $\rho = \rho(\bm{x},t)$ and $\bm{u} = \bm{u}(\bm{x},t)$ are the density and the velocity, $P(\rho)$ is the pressure, $H(\bm{x},\rho)$ contains the attractive-repulsive effects from external $V$ or interaction potentials $W$, assumed to be locally integrable, given by
\begin{equation*}\label{eq:pot}
H(\bm{x},\rho)=V(\bm{x})+W(\bm{x})\star \rho,
\end{equation*}
and $\psi(\bm{x})$ is a nonnegative symmetric smooth function called the communication function in the Cucker-Smale model \cite{cucker2007emergent,cucker2007mathematics} describing collective behavior of systems due to alignment \cite{carrillo2017review}.

The fractional-step methods \cite{leveque2002finite} have been the
widely-employed tool to simulate the temporal evolution of balance laws such
as \eqref{eq:generalsys}. They are based on a division of the problem in
\eqref{eq:generalsys} into two simpler subproblems: the homogeneous
hyperbolic system without source terms and the temporal evolution of density
and momentum without the flux terms but including the sources. These
subproblems are then resolved alternatively employing suitable numerical
methods for each. This procedure introduces a splitting error which is
acceptable for the temporal evolution, but becomes critical when the
objective is to preserve the steady states. This is due to the fact that the
steady state is reached when the fluxes are exactly balanced with the source
terms in each discrete node of the domain. However, when solving
alternatively the two subproblems, this discrete balance can never be
achieved, since the fluxes and source terms are not resolved simultaneously.

To correct this deficiency, well-balanced schemes are designed to discretely
satisfy the balance between fluxes and sources when the steady state is
reached \cite{bouchut2004nonlinear}. The strategy to construct well-balanced
schemes relies on the fact that, when the steady state is reached, there are
some constant relations of the variables that hold in the domain. These
relations allow the resolution of the fluxes and sources in the same level,
thus avoiding the division that the fractional-step methods introduce.
Moreover, if the system enjoys a dissipative property and it has a Liapunov
functional, obtaining analogous tools at the discrete level is key for the
derivation of well-balanced schemes. In this work the steady-state relations
and the dissipative property are obtained by means of the associated free
energy, which in the case of the system in \eqref{eq:generalsys} is
formulated as
\begin{equation}\label{eq:freeenergy}
  \mathcal{F}[\rho]=\int_{\R^{d}}\Pi(\rho)d\bm{x}+\int_{\R^{d}}V(\bm{x})\rho(\bm{x})d\bm{x}+\frac{1}{2}\int_{\R^{d}}\int_{\R^{d}}W(\bm{x}-\bm{y})\rho(\bm{x})\rho(\bm{y}) d\bm{x} d\bm{y},
\end{equation}
where
\begin{equation}\label{eq:PandPi}
\rho \Pi''(\rho)=P'(\rho).
\end{equation}
The pressure $P(\rho)$ and the potential term $H(\bm{x},\rho)$ appearing in the general system \eqref{eq:generalsys} can be gathered by considering the associated free energy. Taking into account that the variation of the free energy in  \eqref{eq:freeenergy} with respect to the density $\rho$ is equal to
\begin{equation}\label{eq:varfreeenergy}
\frac{\delta \mathcal{F}}{\delta \rho}=\Pi'(\rho)+H(\bm{x},\rho),
\end{equation}
it follows that the general system \eqref{eq:generalsys} can be written in a
compact form as
\begin{equation}\label{eq:generalsys2}
\begin{cases}
\partial_{t}\rho+\nabla\cdot\left(\rho \bm{u}\right)=0,\quad \bm{x}\in\R^{d},\quad t>0,\\[2mm]
{\displaystyle \pa_{t}(\rho \bm{u})\!+\!\nabla\!\cdot\!(\rho \bm{u}\otimes \bm{u})\!=-\rho \nabla \frac{\delta \mathcal{F}}{\delta \rho}-\gamma\rho \bm{u}-\!\rho\!\!\int_{\R^{d}}\!\!\psi(\bm{x}-\bm{y})(\bm{u}(\bm{x})-\bm{u}(\bm{y}))\rho(\bm{y})\,d\bm{y}
.}
\end{cases}
\end{equation}

The system in \eqref{eq:generalsys2} is rather general containing a wide
variety of physical problems all under the so-called density functional
theory (DFT) and its dynamic extension (DDFT) see e.g.
\cite{goddard2012general,goddard2012unification,goddard2012overdamped,yatsyshin2012spectral,yatsyshin2013geometry,duran2016dynamical}
and the references therein. A variety of well-balanced schemes have already
been constructed for specific choices of the terms $\Pi(\rho)$, $V(\bm{x})$
and $W(\bm{x})$ in the free energy \eqref{eq:freeenergy}, see
\cite{audusse2004fast, bouchut2004nonlinear, filbet2005approximation} for
instance. Here the focus is set on the free energy and the natural structure
of the system \eqref{eq:generalsys2}. It is naturally advantageous to
consider the concept of free energy in the construction procedure of
well-balanced schemes, since they rely on relations that hold in the steady
states, and moreover, the variation of the free energy with respect to the
density is constant when reaching these steady states, more precisely
\begin{equation}\label{eq:steadyvarener}
\frac{\delta \mathcal{F}}{\delta \rho}=\Pi'(\rho)+H(\bm{x},\rho)= \text{constant on each connected component of}\ \mathrm{supp}(\rho)\ \text{and}\ u=0,
\end{equation}
where the constant can vary on different connected components of
$\mathrm{supp}(\rho)$. As a result, the constant relations in the steady
states, which are needed for well-balanced schemes, are directly provided by
the variation of the free energy with respect to the density.

The steady state relations in \eqref{eq:steadyvarener} hold due to the
dissipation of the linear damping $-\rho \bm{u}$ or nonlinear damping in the
system \eqref{eq:generalsys}, which eventually eliminates the momentum of the
system. This can be justified by means of the total energy of the system,
defined as the sum of kinetic and free energy,
\begin{equation}\label{eq:totalenergy}
  E(\rho,\bm{u})=\int_{\R^{d}}\frac{1}{2}\rho \left|\bm{u}\right|^2 d\bm{x}+\mathcal{F}(\rho),
\end{equation}
since it is formally dissipated, see \cite{giesselmann2017relative,carrillo2017weak,carrillo2018longtime}, as
\begin{equation}\label{eq:equalenergy}
  \frac{dE(\rho,\bm{u})}{dt}=-\gamma\int_{\R^{d}}\rho \left|\bm{u}\right|^2 d\bm{x}-\!\int_{\R^{d}}\int_{\R^{d}}\!\!\psi(\bm{x}-\bm{y})\left|\bm{u}(\bm{y})-\bm{u}(\bm{x})\right|^2 \rho(\bm{x})\,\rho(\bm{y})\,d\bm{x}\,d\bm{y}.
\end{equation}
This last dissipation equation ensures that the total energy $E(\rho,\bm{u})$
keeps decreasing in time while there is kinetic energy in the system. At the
same time, since the definition of the total energy \eqref{eq:totalenergy}
also depends on the velocity $\bm{u}$, it results that the velocity
throughout the domain eventually vanishes. When $\bm{u}=\bm{0}$ throughout
the domain, the momentum equation in \eqref{eq:generalsys2} reduces to
\begin{equation*}
\rho \nabla \frac{\delta \mathcal{F}}{\delta \rho}=\bm{0},
\end{equation*}
meaning that in the support of the density the steady state relation \eqref{eq:steadyvarener} holds. However, for those points outside the support of the density and satisfying $\rho=0$, the variation of the free energy with respect to the density does not need to keep the constant value when the steady state is reached. A discussion of the resulting steady states depending on $\Pi(\rho)$ and $H(\bm{x},\rho)$ is provided in \cite{carrillo2015finite,carrillo2016nonlinear,hoffmann2017keller}.

The system \eqref{eq:generalsys} also satisfies an entropy identity
\begin{equation}\label{eq:entrineq}
 \pa_{t} \eta(\rho, \rho \textbf{u})+ \nabla \cdot \bm{G}(\rho, \rho \textbf{u})= -\rho \bm{u}\cdot \nabla H(\bm{x},\rho)-\gamma\rho \left|\bm{u}\right|^2-\rho \int_{\R^{d}}\!\!\psi(\bm{x}-\bm{y})\,\bm{u}(\bm{x}) \cdot (\bm{u}(\bm{x})-\bm{u}(\bm{y}))\rho(\bm{y})\,d\bm{y},
\end{equation}
where $\eta(\rho,\rho\bm{u})$ and $\bm{G}(\rho,\rho\bm{u})$ are the entropy and the entropy flux defined as
\begin{equation}\label{eq:entropy}
\eta(\rho, \rho \bm{u})=\rho \frac{\left|\bm{u}\right|^2}{2}+\Pi (\rho), \quad \bm{G} (\rho, \rho \bm{u})=\rho \bm{u}\left(\frac{\left|\bm{u}\right|^2}{2}+\Pi' (\rho) \right).
\end{equation}
From a physical point of view the entropy is always a convex function of the density
 \cite{lasota2013chaos}. As a result, from \eqref{eq:entropy} it is justified to assume that $\Pi(\rho)$ is convex, meaning that $\Pi'(\rho)$ has an inverse function for positive densities $\rho$. This last fact is a necessary requirement for the construction of the well-balanced schemes of this work, as it is explained in section \ref{sec:numsch}. Finally, notice that from the entropy identity \eqref{eq:entrineq}, one recovers the free energy dissipation \eqref{eq:equalenergy} by integration using the continuity equation to deal with the forces term $H(\bm{x},\rho)$ and using symmetrization of the nonlinear damping term due to $\psi$ being symmetric.

Let us also point out that the evolution of the center of mass of the density can be computed in some particular cases. In fact, it is not difficult to deduce from \eqref{eq:generalsys2} that
\begin{equation}\label{eq:moment}
\frac{d}{dt}\int_{\R^{d}}\bm{x} \rho d\bm{x} =\int_{\R^{d}} \rho \bm{u} d\bm{x} \quad \mbox{and} \quad \frac{d}{dt}\int_{\R^{d}} \rho \bm{u} d\bm{x}= - \int_{\R^{d}} \nabla V(\bm{x}) \rho d\bm{x} -\gamma \int_{\R^{d}} \rho \bm{u} d\bm{x}\,,
\end{equation}
due to the antisymmetry of $\nabla W(\bm{x})$ and the symmetry of $\psi(\bm{x})$. Therefore, in case $V(\bm{x})$ is not present or quadratic, \eqref{eq:moment} are explicitly solvable. Moreover, if the potential $V(\bm{x})$ is symmetric, the initial data for the density is symmetric, and the initial data for the velocity is antisymmetric, then the solution to \eqref{eq:generalsys2} keeps these symmetries in time, i.e., the density is symmetric and the velocity is antisymetric for all times, and the center of mass is conserved
$$
\frac{d}{dt}\int_{\R^{d}}\bm{x} \rho d\bm{x}= 0\,.
$$

The steady state relations \eqref{eq:steadyvarener} only hold when the linear damping term is included in system \eqref{eq:generalsys}. When only the nonlinear damping of Cucker-Smale type is present, the system has the so-called moving steady states, see \cite{carrillo2016pressureless,carrillo2017review,carrillo2018longtime}, which satisfy the more general relations
\begin{equation}\label{eq:steadyvarener2}
\frac{\delta \mathcal{F}}{\delta \rho}=\text{constant on each connected component of}\ \mathrm{supp}(\rho)\ \text{and}\ \bm{u}=\text{constant} .
\end{equation}
However, the construction of well-balanced schemes satisfying the moving steady state relations has proven to be more difficult than for the still steady states \eqref{eq:steadyvarener} without dissipation. For literature about well-balanced schemes for moving steady states without dissipation, we refer to \cite{noelle2007high,xing2011advantage}.

The most popular application in the literature for well-balanced schemes deals with the Saint-Venant system for shallow water flows with nonflat bottom \cite{audusse2004fast, bouchut2004nonlinear, canestrelli2009well, liang2009numerical, xing2005high, xing2014survey}, for which $\Pi(\rho)=\frac{g}{2}\rho^2$, with $g$ being the gravity constant, and $H(\bm{x},\rho)$ depends on the bottom. Here it is important to remark the work of Audusse et al. in \cite{audusse2004fast}, where they propose a hydrostatic reconstruction that has successfully inspired more sophisticated well-balanced schemes in the area of shallow water equations \cite{marche2007evaluation, noelle2006well}. Another area where well-balanced schemes have been fruitful is chemosensitive movement, with the works of Filbet, Shu and their collaborators \cite{filbet2005derivation,filbet2005approximation,gosse2012asymptotic,xing2006high}. In this case the pressure satisfies $\Pi(\rho)=\rho \left(\ln(\rho)-1\right)$  and $H$ depends on the chemotactic sensitivity and the chemical concentration. The list of applications of the system \eqref{eq:generalsys} continues growing with more choices of  $\Pi(\rho)$ and $H(\bm{x},\rho)$ \cite{xing2006high}: the elastic wave equation, nozzle flow problem, two phase flow model, etc.

The orders of accuracy from the finite volume well-balanced schemes presented
before range from first- and second-order \cite{audusse2004fast,
leveque1998balancing, xu2002well,kurganov2007second,liang2009numerical} to
higher-order versions \cite{xing2006high,
vukovic2004weno,noelle2006well,gallardo2007well}. Again, the most popular
application has been shallow water equations, and the survey from Xing and
Shu \cite{xing2014survey} provides a summary of all the shallow water methods
with different accuracies. Some of the previous schemes proposed were
equipped to satisfy natural properties of the systems under consideration,
such as nonnegativity of the density
\cite{audusse2005well,kurganov2007second} or the satisfaction of a discrete
entropy inequality \cite{audusse2004fast,filbet2005approximation}, enabling
also the computation of dry states \cite{gallardo2007well} . Theoretically
the Godunov scheme satisfies all these properties \cite{leroux1999riemann},
but its main drawback is its computationally expensive implementation. The
high-order schemes usually rely on the WENO reconstructions originally
proposed by Jiang and Shu \cite{jiang1996efficient}.

Other well-balanced numerical approaches employed to simulate the system
\eqref{eq:generalsys2} are finite differences
\cite{xing2006high2,xing2005high}, which are equivalent to the finite volume
methods for first-and second-order, and the discontinuous Galerkin methods
\cite{xing2006high}. The overdamped system of \eqref{eq:generalsys2} with
$\psi\equiv 0$, obtained in the free inertia limit where the momentum reaches
equilibrium on a much faster timescale than the density, has also been
numerically resolved for general free energies of the form \eqref{eq:freeenergy}, via finite volume schemes \cite{carrillo2015finite} or
discontinuous Galerkin approaches \cite{sun2018discontinuous}. This scheme for the overdamped system also conserves the dissipation of the free energy at the discrete approximation.

\red{ The novelty of this work is twofold. Foremost, all these previous
schemes were only applicable for specific choices of $\Pi(\rho)$ and
$H(\bm{x},\rho)$, meaning that a general scheme valid for a wide range of
applications is lacking. And while some previous schemes \cite{xing2006high}
could be employed in more general cases, the focus in the literature has been
on the shallow water and chemotaxis equations. In addition, the function
$H(\bm{x},\rho)$, which results from summing $V(\bm{x})$ and $W(\bm{x})\star
\rho$ as in \eqref{eq:pot}, has so far been taken as dependent on $\bm{x}$
only, unlike the present work where it depends on $\rho$ by means of the
convolution with an interaction potential $W(\bm{x})$.}

In this work we present a finite volume scheme for a general choice of $\Pi(\rho)$
and $H(\bm{x},\rho)$ which is first- and second-order accurate and satisfies
the nonnegativity of the density, the well-balanced property, the
semidiscrete entropy inequality and the semidiscrete free energy dissipation.
Furthermore, as it is shown in example \ref{ex:hardrods} of section
\ref{sec:numtest}, it can also be applied to more general free energies than
the one in \eqref{eq:freeenergy} and with the form
\begin{equation}\label{eq:freeenergygeneral}
\mathcal{F}[\rho]=\int_{\R^{d}}\Pi(\rho)d\bm{x}+\int_{\R^{d}}V(\bm{x})\rho(\bm{x})d\bm{x}+\frac{1}{2}\int_{\R^{d}}K\left(W(\bm{x})\star \rho(\bm{x})\right) \rho(\bm{x})d\bm{x},
\end{equation}
where $K$ is a function depending on the convolution of $\rho(\bm{x})$ with the kernel $W(\bm{x})$. Its variation with respect to the density satisfies
\begin{equation}\label{eq:varfreeenergygeneral}
\frac{\delta \mathcal{F}}{\delta \rho}=\Pi'(\rho)+V(\bm{x})+\frac{1}{2} K\left(W(\bm{x})\star \rho \right)+\frac{1}{2} K'\left(W(\bm{x})\star \rho\right) \left(W(\bm{x})\star \rho\right).
\end{equation}
These free energies arise in applications related to
(D)DFT\cite{duran2016dynamical,goddard2012unification}, see
\cite{carrillo2017blob} for other related free energies and properties.

\red{The other novel technical aspect of this work concerns the numerical
treatment of the different source terms in \eqref{eq:generalsys}. In fact, in
order to keep the well-balanced property and the decay of the free energy we
treat source terms differently. While the dissipative terms are harmless and
treated by direct approximations, the fundamental question is how to choose
the discretization of the potential term given by $H(\bm{x},\rho)
=V(\bm{x})+W(\bm{x})\star \rho$. For this purpose we appropriately extend the
ideas in \cite{bouchut2004nonlinear,filbet2005approximation} to our case to
keep the well-balanced property and the energy decay. The condition for
stationary states \eqref{eq:steadyvarener} is crucial in defining an
approximation of the term $-\rho \nabla H(\bm{x},\rho) $ by a dicretization
of $\nabla P(\rho)$ which is consistent when the new reconstructed values of
the density at the interfaces taking into account the potential
$H(\bm{x},\rho) $. This general treatment includes as specific cases both the
shallow-water equations \cite{audusse2004fast, bouchut2004nonlinear} and the
hyperbolic chemotaxis problem \cite{filbet2005approximation}.}

Section \ref{sec:numsch} describes the first- and second-order well-balanced
scheme reconstructions, and provides the proofs of their main properties.
Section \ref{sec:numtest} contains the numerical simulations, with a first
subsection \ref{subsec:val} where the validation of the well-balanced
property and the orders of accuracy is conducted, and a second subsection
\ref{subsec:numexp} with numerical experiments from different applications. A
wide range of free energies is employed to remark the extensive nature of our
well-balanced scheme. A short summary and conclusions are offered in section
4.
%
%
%
%
\section{Well-balanced finite volume scheme}\label{sec:numsch}
The terms appearing in the one-dimensional system \eqref{eq:generalsys2} are usually gathered in the form of
\begin{equation}\label{eq:compactsys}
\pa_t U + \pa_x F(U) = S(x,U),
\end{equation}
with
\begin{equation*}
U=\begin{pmatrix}
\rho \\ \rho u
\end{pmatrix}, \quad
F(U)=\begin{pmatrix}
\rho u \\ \rho u^2+P(\rho)
\end{pmatrix}
\end{equation*}
and
\begin{equation*}
S(x,U)=\begin{pmatrix}
0 \\ -\rho \pa_x H -\gamma\rho u-\rho \displaystyle\int_{\R}\psi(x-y)(u(x)-u(y))\rho(y)\,dy
\end{pmatrix},
\end{equation*}
where $U$ are the unknown variables, $F(U)$ the fluxes and $S(U)$ the sources.
The one-dimensional finite volume approximation of  \eqref{eq:compactsys} is obtained by breaking the domain into grid cells $\left(x_{i-1/2}\right)_{i\in \mathbb{Z}}$ and approximating in each of them the cell average of $U$. Then these cell averages are modified after each time step, depending on the flux through the edges of the grid cells and the cell average of the source term \cite{leveque2002finite}. Finite volume schemes for hyperbolic systems employ an upwinding of the fluxes and in the semidiscrete case they provide a discrete version of \eqref{eq:compactsys} under the form
\begin{equation}\label{eq:fvbasic}
\dfrac{d U_i}{dt}=-\dfrac{F_{i+\frac{1}{2}}-F_{i-\frac{1}{2}}}{\Delta x_i} + S_i,
\end{equation}
where the cell average of $U$ in the cell $\left(x_{i-\frac{1}{2}},x_{i+\frac{1}{2}}\right)$ is denoted as
\begin{equation*}
U_i=\begin{pmatrix}
\rho_i \\ \rho_i u_i
\end{pmatrix},
\end{equation*}
$F_{i+\frac{1}{2}}$ is an approximation of the flux $F(U)$ at the point $x_{i+\frac{1}{2}}$, $S_i$ is an approximation of the source term $S(x,U)$ in the cell $\left(x_{i-\frac{1}{2}},x_{i+\frac{1}{2}}\right)$ and $\Delta x_i$ is the possibly variable mesh size $\Delta x_i=x_{i+\frac{1}{2}}-x_{i-\frac{1}{2}}$.

The approximation of the flux $F(U)$ at the point $x_{i+\frac{1}{2}}$,
denoted as $F_{i+\frac{1}{2}}$, is achieved by means of a numerical flux
$\mathscr{F}$ which depends on two reconstructed values of $U$ at the left
and right of the boundary between the cells $i$ and $i+1$. These two values,
$U_{i+\frac{1}{2}}^-$ and $U_{i+\frac{1}{2}}^+$, are computed from the cell
averages following different construction procedures that seek to satisfy
certain properties, such as order of accuracy or nonnegativity. Two
widely-employed reconstruction procedures are the second-order finite volume
monotone upstream-centered scheme for conservation laws, referred to as MUSCL
\cite{osher1985convergence}, or the weighted-essentially non-oscllatory
schemes, widely known as WENO \cite{shu1998essentially}.

Once these two reconstructed values are computed, $F_{i+\frac{1}{2}}$ is obtained from
\begin{equation}\label{eq:numflux}
F_{i+\frac{1}{2}}=\mathscr{F}\left(U_{i+\frac{1}{2}}^-,U_{i+\frac{1}{2}}^+\right).
\end{equation}
The numerical flux $\mathscr{F}$ is usually denoted as Riemann solver,  since it
provides a stable resolution of the Riemann problem located at the cell interfaces, where the left value of the variables in $U_{i+\frac{1}{2}}^-$ and the right value $U_{i+\frac{1}{2}}^+$. The
literature concerning Riemann solvers is vast and there are different choices
for it \cite{toro2013riemann}: Godunov, Lax-Friedrich, kinetic, Roe, etc.
Some usual properties of the numerical flux that are assumed
\cite{audusse2004fast, bouchut2004nonlinear, filbet2005approximation} are:

\begin{enumerate}[label=\arabic*.]
\item It is consistent with the physical flux, so that $\mathscr{F}(U,U)=F(U)$.
\item It preserves the nonnegativity of the density $\rho_i (t)$ for the homogeneous problem, where the numerical flux is computed as in \eqref{eq:numflux}.
\item It satisfies a cell entropy inequality for the entropy pair \eqref{eq:entropy} for the homogeneous problem. Then, according to \cite{bouchut2004nonlinear}, it is possible to find a numerical entropy flux $\mathscr{G}$ such that
\begin{multline}\label{eq:cellentropyineq1}
\hspace{0.9cm} G(U_{i+1})+\nabla_U \,\eta(U_{i+1})\left(\mathscr{F}(U_{i},U_{i+1})-F(U_{i+1})\right)\\
\leq \mathscr{G}(U_{i},U_{i+1}) \leq G(U_{i})+\nabla_U \,\eta(U_{i})\left(\mathscr{F}(U_{i},U_{i+1})-F(U_{i})\right),
\end{multline}
where $\nabla_U\,\eta$ is the derivative of $\eta$ with respect to $U=\begin{pmatrix}
\rho \\ \rho u\end{pmatrix}$.
\end{enumerate}

The first- and second-order well-balanced schemes described in this section propose an alternative reconstruction procedure for $U_{i+\frac{1}{2}}^-$ and $U_{i+\frac{1}{2}}^+$ which ensures that the steady state in \eqref{eq:steadyvarener} is discretely preserved when starting from that steady state. Subsections \ref{subsec:firstorder} and \ref{subsec:secondorder} contain the first- and second-order schemes, respectively, together with their proved properties.

%
%
\subsection{First-order scheme}\label{subsec:firstorder}

The basic first-order schemes approximate the flux $F_{i+\frac{1}{2}}$ by a numerical flux $\mathscr{F}$ which depends on the cell averaged values of $U$ at the two adjacent cells, so that the inputs for the numerical flux in \eqref{eq:numflux} are
\begin{equation}\label{eq:numflux1}
F_{i+\frac{1}{2}}=\mathscr{F}\left(U_i,U_{i+1}\right).
\end{equation}

The resolution of the finite volume scheme in  \eqref{eq:fvbasic} with a numerical flux of the form in \eqref{eq:numflux1} and a cell-centred evaluation of $-\rho \partial_x H$ for the source term $S_i$ is not generally able to preserve the steady states, as it was shown in the initial works of well-balanced schemes \cite{greenberg1996well, gosse1996well}. These steady states are provided in \eqref{eq:steadyvarener}, and satisfy that the variation of the free energy with respect to the density has to be constant in each connected component of the support of the density. The discrete steady state is defined in a similar way,
\begin{equation}\label{eq:steadyvarenerdiscrete}
\left(\frac{\delta \mathcal{F}}{\delta \rho}\right)_i=\Pi'(\rho_i)+H_i= C_\Gamma \text{ in each}\  \Lambda_\Gamma, \Gamma\in\mathbb{N}\,,
\end{equation}
where $\Lambda_\Gamma$, $\Gamma\in\mathbb{N}$, denotes the possible infinite sequence indexed by $\Gamma$ of subsets $\Lambda_\Gamma$ of subsequent indices $i\in\mathbb{Z}$ where $\rho_i>0$ and $u_i=0$, and $C_\Gamma$ the corresponding constant in that connected component of the discrete support.

As it was emphasized above, the preservation of these steady states for particular choices of $\Pi'(\rho)$ and $H(x,\rho)$, such as shallow water \cite{audusse2004fast} or chemotaxis \cite{filbet2005approximation}, is paramount.
A solution to allow this preservation was proposed in the work of Audusse et al. \cite{audusse2004fast}, where instead of evaluating the numerical flux as in \eqref{eq:numflux}, they chose
\begin{equation}\label{eq:numfluxnew}
F_{i+\frac{1}{2}}=\mathscr{F}\left(U_{i+\frac{1}{2}}^-,U_{i+\frac{1}{2}}^+\right), \ \text{where}\ U_{i+\frac{1}{2}}^{\pm}=\begin{pmatrix}
\rho_{i+\frac{1}{2}}^{\pm} \\ \rho_{i+\frac{1}{2}}^{\pm} u_{i+\frac{1}{2}}^{\pm}
\end{pmatrix}.
\end{equation}
The interface values $U_{i+\frac{1}{2}}^{\pm}$ are reconstructed from $U_i$ and $U_{i+1}$ by taking into account the steady state relation in \eqref{eq:steadyvarenerdiscrete}.  Contrary to other works in which the interface values are reconstructed to increase the order of accuracy, now the objective is to satisfy the well-balanced property. Bearing this in mind, we make use of \eqref{eq:steadyvarenerdiscrete} to the cells with centred nodes at $x_i$ and $x_{i+1}$ to define the interface values such that
\begin{equation*}
\begin{gathered}
  \Pi'\left(\rho_{i+\frac{1}{2}}^{-}\right)+H_{i+\frac{1}{2}}=\Pi'\left(\rho_{i}\right)+H_{i},\\
  \Pi'\left(\rho_{i+\frac{1}{2}}^{+}\right)+H_{i+\frac{1}{2}}=\Pi'\left(\rho_{i+1}\right)+H_{i+1},
\end{gathered}
\end{equation*}
where the term $H_{i+\frac{1}{2}}$ is evaluated to preserve consistency and stability, with an upwind or average value obtained as
\begin{equation}\label{eq:hinterface}
H_{i+\frac{1}{2}}=\max\left(H_{i},H_{i+1}\right)\quad \textrm{or} \quad H_{i+\frac{1}{2}}=\frac{1}{2}\left(H_{i}+H_{i+1}\right).
\end{equation}
Then, by denoting as $\xi(s)$ the inverse function of $\Pi'(s)$ for $s>0$, we conclude that the interface values $U_{i+\frac{1}{2}}^{\pm}$ are computed as
\begin{equation}\label{eq:rhointerface}
\begin{gathered}
 \rho_{i+\frac{1}{2}}^{-}=\xi \left(\Pi'\left(\rho_{i}\right)+H_{i}-H_{i+\frac{1}{2}}\right)_+,\quad u_{i+\frac{1}{2}}^{-}=u_i,\\
  \rho_{i+\frac{1}{2}}^{+}=\xi \left(\Pi'\left(\rho_{i+1}\right)+H_{i+1}-H_{i+\frac{1}{2}}\right)_+,\quad u_{i+\frac{1}{2}}^{+}=u_{i+1}.
\end{gathered}
\end{equation}
The function $\xi(s)$ is well-defined for $s>0$ since $\Pi(s)$ is strictly convex, $\Pi''(s)>0$. This is always the case since, as mentioned in the introduction, the physical entropies are always strictly convex from \eqref{eq:entropy}. However, some physical entropies and applications allow for vacuum of the steady states, therefore we need to impose the value of $\rho_{i+\frac{1}{2}}^{\pm}$, given that they should be nonnegative. \red{Henceforth, $\xi(s)$ denotes} the extension by zero of the inverse of $\Pi'(s)$ whenever $s>0$.

Furthermore, the discretization of the source term is taken as
\begin{equation}\label{eq:sourcewellbalanced}
S_i=\frac{1}{\Delta x_i}\begin{pmatrix}
0 \\  P\left(\rho_{i+\frac{1}{2}}^{-}\right) - P\left(\rho_{i-\frac{1}{2}}^{+}\right)
\end{pmatrix}-\begin{pmatrix}
0 \\ \gamma \rho_i u_i + \rho_i \displaystyle\sum_{j} \Delta x_j (u_i-u_j)\rho_j \psi_{ij}
\end{pmatrix},
\end{equation}
which is motivated by the fact that in the steady state, with $u=0$ in \eqref{eq:compactsys}, the fluxes are balanced with the sources,
\begin{equation*}
\rho \partial_x \Pi'(\rho)=-\rho \partial_x H.
\end{equation*}
Here, $\psi_{ij}$ is an approximation of the average value of $\psi$ on the interval centred at $x_i-x_j$ of length $\Delta x_j$. From here, and integrating over the cell volume, it results that
\begin{equation}\label{eq:integratesource}
\int_{x_{i-\frac{1}{2}}}^{x_{i+\frac{1}{2}}}-\rho \partial_x H \,dx=\int_{x_{i-\frac{1}{2}}}^{x_{i+\frac{1}{2}}}\rho \partial_x \Pi'(\rho)\,dx=\int_{x_{i-\frac{1}{2}}}^{x_{i+\frac{1}{2}}}\partial_x P(\rho)\,dx=P(\rho_{i+\frac{1}{2}}^{-})-P(\rho_{i-\frac{1}{2}}^{+}),
\end{equation}
with the relation between $\Pi'(\rho)$ and $P(\rho)$ was given in  \eqref{eq:PandPi}. This idea of distributing the source terms along the interfaces has already been explored in previous works \cite{katsaounis2004upwinding}.

The discretization of the source term in \eqref{eq:sourcewellbalanced} entails that the discrete balance between fluxes and sources is accomplished when $F_{i+\frac{1}{2}}=P(\rho_{i+\frac{1}{2}}^{-})=P(\rho_{i+\frac{1}{2}}^{+})$. The computation of the numerical fluxes expressed in \eqref{eq:numfluxnew}, in which the interface values $U_{i+\frac{1}{2}}^{\pm}$ are considered, enables this balance if in the steady states $U_{i+\frac{1}{2}}^{-}=U_{i+\frac{1}{2}}^{+}=(\rho_{i+\frac{1}{2}}^{-},0)=(\rho_{i+\frac{1}{2}}^{+},0)$. Moreover, the discretization of the source term as in \eqref{eq:sourcewellbalanced} may seem counter-intuitive when the system is far away from the steady state, given that the balanced expressed in \eqref{eq:integratesource} only holds in those states. In spite of this, the consistency with the original system in \eqref{eq:compactsys} is not lost, as it will be proved in subsection \ref{subsec:firstorderprop}.

Let us finally discuss the discretization of the potential
$H(x,\rho)=V(x)+W\ast\rho (x)$. We will always approximate it as
$$
H_i=V_i+\sum_j \Delta x_j W_{ij} \rho_j \,, \mbox{ for all } i\in\mathbb{Z}\,,
$$
where $V_i=V(x_i)$ and $W_{ij}=W(x_i-x_j)$ in case the potential is smooth or choosing $W_{ij}$ as an average value of $W$ on the interval centred at $x_i-x_j$ of length $\Delta x_j$ in case of general locally integrable potentials $W$. Let us also point out that this discretization keeps the symmetry of the discretized interaction potential $W_{ij}=W_{ji}$ for all $i,j\in\mathbb{Z}$ whenever $W$ is smooth or solved with equal size meshes $\Delta x_i=\Delta x_j$ for all $i,j\in\mathbb{Z}$.

%
%
\subsection{Properties of the first-order scheme}\label{subsec:firstorderprop}

The first-order semidiscrete scheme defined in \eqref{eq:fvbasic}, constructed with \eqref{eq:numfluxnew}-\eqref{eq:sourcewellbalanced}, and for a numerical flux $\mathscr{F}\left(U_{i},U_{i+1}\right)=\left(\mathscr{F}^{\rho},\mathscr{F}^{\rho u}\right)\left(U_{i},U_{i+1}\right)$ satisfying the properties stated in the introduction of section \ref{sec:numsch}, satisfies:
\begin{enumerate}[label=(\roman*)]
\item preservation of the nonnegativity of $\rho_i(t)$;
\item well-balanced property, thus preserving the steady states given by \eqref{eq:steadyvarenerdiscrete};
\item consistency with the system \eqref{eq:generalsys2};
\item cell entropy inequality associated to the entropy pair \eqref{eq:entropy},
\begin{equation}\label{eq:cellentropyineq2}
\Delta x_i \frac{d\eta_i}{dt}+\Delta x_i H_i \frac{d\rho_i}{dt}+G_{i+\frac{1}{2}}-G_{i-\frac{1}{2}}=-u_i\left(\gamma\Delta x_i \rho_i u_i+\Delta x_i \rho_i \sum_j \Delta x_j \rho_j \left(u_i-u_j\right)\psi_{ij}\right)\,,
\end{equation}
where $\eta_i=\Pi\left(\rho_i\right)+\frac{1}{2}\rho_i u_i^2$ and
$$
G_{i+\frac{1}{2}}=\mathscr{G}\left(U^-_{i+\frac{1}{2}}, U^+_{i+\frac{1}{2}}\right)+\mathscr{F}^\rho\left(U^-_{i+\frac{1}{2}}, U^+_{i+\frac{1}{2}}\right)H_{i+\frac{1}{2}}.
$$
\item the discrete analog of the free energy dissipation property \eqref{eq:equalenergy} given by
\begin{equation}\label{eq:disfreeenergydiscrete}
\frac{d}{dt} E^\Delta(t)\leq -\gamma\sum_i \Delta x_i \rho_i u_i^2-\frac12 \sum_{i,j} \Delta x_i \Delta x_j \rho_i \rho_j \left(u_i-u_j\right)^2\psi_{ij}
\end{equation}
with
\begin{equation}\label{eq:freeenergydiscrete}
E^\Delta= \sum_i \frac{\Delta x_i}{2}\rho_i u_i^2 +
\mathcal{F}^\Delta \quad \mbox{ and } \quad
\mathcal{F}^\Delta = \sum_i \Delta x_i \left[\Pi\left(\rho_i\right)+ V_i\rho_i \right]+\frac12 \sum_{i,j} \Delta x_i \Delta x_j W_{ij} \rho_i \rho_j .
\end{equation}
\item the discrete analog of the evolution for centre of mass in \eqref{eq:moment},
\begin{equation}\label{eq:proofevolcentremass}
\frac{d}{dt}\left(\sum_i \Delta x_i \rho_i x_i  \right)=\sum_i \Delta x_i \mathscr{F}^\rho\left(U^-_{i+\frac{1}{2}}, U^+_{i+\frac{1}{2}}\right) ,
\end{equation}
which is reduced to
\begin{equation}\label{eq:proofevolcentremasssimp}
\sum_i \Delta x_i \rho_i x_i  =0\quad
\end{equation}
when the initial density is symmetric and the initial velocity antisymmetric. This implies that the discrete centre of mass is conserved in time and centred at $0$.
\end{enumerate}

\begin{proof} Some of the following proofs follow the lines considered in \cite{audusse2004fast,filbet2005approximation}.
\begin{enumerate}[label=(\roman*)]
\item

\red{If a first-order numerical flux
$\mathscr{F}\left(U_{i},U_{i+1}\right)=\left(\mathscr{F}^{\rho},\mathscr{F}^{\rho
u}\right)\left(U_{i},U_{i+1}\right)$ for the homogeneous problem, such as
the Lax-Friedrich scheme detailed in Appendix \ref{app:numerics},
satisfies the nonnegativity of the density $\rho_i(t)$,  then it
necessarily follows that
\begin{equation}\label{eq:proofnegat1}
\mathscr{F}^{\rho}((\rho_i=0,u_i),(\rho_{i+1},u_{i+1}))-\mathscr{F}^{\rho}((\rho_{i-1},u_{i-1}),(\rho_i=0,u_i))\leq 0\quad \forall (\rho_i,u_i)_i.
\end{equation}
In our case, the sources do not contribute to the continuity equation in
\eqref{eq:compactsys}, and for the numerical flux in
\eqref{eq:numfluxnew} we need to check that
\begin{equation}\label{eq:proofnegat2}
\mathscr{F}^{\rho}\left(U_{i+\frac{1}{2}}^-,U_{i+\frac{1}{2}}^+\right)-\mathscr{F}^{\rho}\left(U_{i-\frac{1}{2}}^-,U_{i-\frac{1}{2}}^+\right)\leq 0
\end{equation}
whenever $\rho_i=0$. When $\rho_i=0$, the reconstruction in
\eqref{eq:hinterface} and \eqref{eq:rhointerface} yields
$\rho_{i+\frac{1}{2}}^-=\rho_{i+\frac{1}{2}}^+=0$ since $\Pi(\rho)$ is
assumed to be convex, and \eqref{eq:proofnegat2} results in
\begin{equation}\label{eq:proofnegat3}
\mathscr{F}^{\rho}((0,u_i),(\rho_{i+\frac{1}{2} }^+,u_{i+1}))-\mathscr{F}^{\rho}((\rho_{i-\frac{1}{2} }^-,u_{i-1}),(\rho_i=0,u_i))\leq 0\quad \forall (\rho_{i+\frac{1}{2}}^+,\rho_{i+\frac{1}{2}}^-,u_i)_i.
\end{equation}
Then, given that the numerical scheme is chosen so that it preserves the
nonnegativity of the density for the homogeneous problem and
\eqref{eq:proofnegat1} holds, it follows that \eqref{eq:proofnegat3} is
satisfied too.}

\item To preserve the steady state the discrete fluxes and source need to be balanced,
\begin{equation}\label{eq:proof2step1}
F_{i+\frac{1}{2}}-F_{i-\frac{1}{2}}=\Delta x S_i.
\end{equation}
When the steady state holds it follows from \eqref{eq:rhointerface} that $\rho_{i+\frac{1}{2}}^-=\rho_{i+\frac{1}{2}}^+$ and $u_{i+\frac{1}{2}}^-=u_{i-\frac{1}{2}}^+=0$, and as a result $U_{i+\frac{1}{2}}^-=U_{i+\frac{1}{2}}^+$. Then, by consistency of the numerical flux $\mathscr{F}$,
\begin{equation}\label{eq:proof2step2}
F_{i+\frac{1}{2}}=\mathscr{F}\left((\rho_{i+\frac{1}{2}}^-,0),(\rho_{i+\frac{1}{2}}^+,0)\right)=F(U_{i+\frac{1}{2}}^-)=F(U_{i+\frac{1}{2}}^+)=\begin{pmatrix} 0 \\ P(\rho_{i+\frac{1}{2}}^-)\end{pmatrix} =\begin{pmatrix} 0 \\ P(\rho_{i+\frac{1}{2}}^+)\end{pmatrix} .
\end{equation}
Concerning the source term $S_i$ of \eqref{eq:sourcewellbalanced}, in the steady state it is equal to
\begin{equation}\label{eq:proof2step3}
\Delta x_i S_i=\begin{pmatrix}
0 \\  P\left(\rho_{i+\frac{1}{2}}^{-}\right) - P\left(\rho_{i-\frac{1}{2}}^{+}\right)
\end{pmatrix}.
\end{equation}
Then the balance in \eqref{eq:proof2step1} is obtained from \eqref{eq:proof2step2} and \eqref{eq:proof2step3}.

\item For the consistency with the original system of \eqref{eq:generalsys2} one has to apply the criterion in \cite{bouchut2004nonlinear}, by which two properties concerning the consistency with the exact flux $F$ and the consistency with the source term need to be checked. Before proceeding, the finite volume discretization in \eqref{eq:fvbasic} needs to be rewritten in a non-conservative form as
\begin{equation}
\begin{gathered}
\dfrac{d U_i}{dt}=-\dfrac{\mathscr{F}_{l}(U_i,U_{i+1},H_i,H_{i+1})-\mathscr{F}_{r}(U_{i-1},U_{i},H_{i-1},H_{i})}{\Delta x_i} \\-\begin{pmatrix}
0 \\ \gamma \rho_i u_i + \rho_i \sum_{j} (u_i-u_j)\rho_j \psi(x_i-x_j)
\end{pmatrix}
\end{gathered}
\end{equation}
where
\red{\begin{equation*}
\begin{gathered}
 \mathscr{F}_{l}(U_i,U_{i+1},H_i,H_{i+1})=F_{i+\frac{1}{2}}-\Delta x_i S_{i+\frac{1}{2}}^-,\\
 \mathscr{F}_{r}(U_{i-1},U_{i},H_{i-1},H_{i})=F_{i-\frac{1}{2}}+\Delta x_i S_{i-\frac{1}{2}}^+.
\end{gathered}
\end{equation*}}
Here the source term $S_i$ is considered as being distributed along the cells interfaces, satisfying
\red{\begin{equation*}
\begin{gathered}
S_i=S_{i+\frac{1}{2}}^-+S_{i-\frac{1}{2}}^+-\begin{pmatrix}
0 \\ \gamma \rho_i u_i + \rho_i \sum_{j} (u_i-u_j)\rho_j \psi(x_i-x_j)
\end{pmatrix},\\
 S_{i+\frac{1}{2}}^-=\frac{1}{\Delta x_i}\begin{pmatrix}
0\\ P(\rho_{i+\frac{1}{2}}^-) - P(\rho_i)
\end{pmatrix}\quad\text{and}\quad
 S_{i-\frac{1}{2}}^+=\frac{1}{\Delta x_i}\begin{pmatrix}
0\\ P(\rho_i)-P(\rho_{i-\frac{1}{2}}^+)
\end{pmatrix}.
\end{gathered}
\end{equation*}}
The consistency with the exact flux means that $ \mathscr{F}_{l}(U,U,H,H)=\mathscr{F}_{r}(U,U,H,H)=F(U)$. This is directly satisfied since $U_{i+\frac{1}{2}}^-=U_i$ and $U_{i+\frac{1}{2}}^+=U_{i+1}$ whenever $H_{i+1}=H_{i}$, due to \eqref{eq:rhointerface}.

For the consistency with the source term the criterion to check is
\begin{equation*}
\mathscr{F}_{r}(U_i,U_{i+1},H_i,H_{i+1})-\mathscr{F}_{l}(U_i,U_{i+1},H_i,H_{i+1})=\begin{pmatrix}0\\ -\rho (H_{i+1}-H_i)+o(H_{i+1}-H_i)\end{pmatrix}
\end{equation*}
as $U_i$, $U_{i+1}\rightarrow U$ and $H_i$, $H_{i+1}\rightarrow H$. \red{For this case,
\begin{equation}
\begin{gathered}
\mathscr{F}_{r}(U_i,U_{i+1},H_i,H_{i+1})-\mathscr{F}_{l}(U_i,U_{i+1},H_i,H_{i+1})=\begin{pmatrix}0\\S_{i+\frac{1}{2}}^++S_{i+\frac{1}{2}}^- \end{pmatrix}=\\\begin{pmatrix}0\\
 \\ -(P(\xi(\Pi'(\rho_{i+1})+H_{i+1}-H_{i+\frac{1}{2}})-P(\rho_{i+1}))+(P(\xi(\Pi'(\rho_i)+H_i-H_{i+\frac{1}{2}})-P(\rho_i))\end{pmatrix},
\end{gathered}
\end{equation}}
where $H_{i+\frac{1}{2}}=\max\left(H_{i},H_{i+1}\right)$. By assuming without loss of generality that $H_{i+\frac{1}{2}}=H_{i}$, the second term of the last matrix results in
\red{\begin{equation*}
-P(\xi(\Pi'(\rho_{i+1})+H_{i+1}-H_{i}))+P(\xi(\Pi'(\rho_{i})))=-P(\xi(\Pi'(\rho_{i+1})+H_{i+1}+H_{i}))-P(\rho_i)\,.
\end{equation*}}
This term can be further approximated as
\red{\begin{equation*}
-(P\circ\xi)'(\Pi'(\rho_{i+1}))\,(H_{i+1}-H_{i})+o(H_{i+1}-H_{i})=-\rho_{i+1}(H_{i+1}-H_{i})+o(H_{i+1}-H_{i})
\end{equation*}}
since
\begin{equation*}
(P\circ\xi)'(\Pi'(\rho_{i+1}))=P'(\rho_{i+1})\frac{1}{\Pi''(\rho_{i+1})}=\rho_{i+1}
\end{equation*}
by taking derivatives in $(\xi \circ \Pi')(\rho)=\rho$ and making use of \eqref{eq:PandPi}. Finally, since $\rho_{i+1}\rightarrow\rho$, the consistency with the source term is satisfied. An analogous procedure can be followed whenever $H_{i+\frac{1}{2}}=H_{i+1}$.

\item To prove \eqref{eq:cellentropyineq2} we follow the strategy from \cite{filbet2005approximation}. We first set $G_{i+\frac{1}{2}}$ to be
\begin{equation*}
G_{i+\frac{1}{2}}=\mathscr{G}\left(U^-_{i+\frac{1}{2}}, U^+_{i+\frac{1}{2}}\right)+\mathscr{F}^\rho\left(U^-_{i+\frac{1}{2}}, U^+_{i+\frac{1}{2}}\right)H_{i+\frac{1}{2}}.
\end{equation*}
Subsequently, and employing the inequalities for $\mathscr{G}\left(U^-_{i+\frac{1}{2}}, U^+_{i+\frac{1}{2}}\right)$ in \eqref{eq:cellentropyineq1}, it follows that
\begin{equation*}
\begin{split}
G_{i+\frac{1}{2}}-G_{i-\frac{1}{2}} \leq &\ G\left(U^-_{i+\frac{1}{2}}\right)+\nabla_U \eta \left(U^-_{i+\frac{1}{2}}\right) \left(\mathscr{F}\left(U^-_{i+\frac{1}{2}}, U^+_{i+\frac{1}{2}}\right)-F\left(U^-_{i+\frac{1}{2}}\right)\right)\\
 & -G\left(U^+_{i-\frac{1}{2}}\right)-\nabla_U \eta \left(U^+_{i-\frac{1}{2}}\right) \left(\mathscr{F}\left(U^-_{i+\frac{1}{2}}, U^+_{i+\frac{1}{2}}\right)-F\left(U^+_{i-\frac{1}{2}}\right)\right)\\ & +\mathscr{F}^\rho\left(U^-_{i+\frac{1}{2}}, U^+_{i+\frac{1}{2}}\right)H_{i+\frac{1}{2}}-\mathscr{F}^\rho\left(U^-_{i-\frac{1}{2}}, U^+_{i-\frac{1}{2}}\right)H_{i-\frac{1}{2}}.
 \end{split}
\end{equation*}
This last inequality can be rewritten after some long computations as
 \begin{equation*}
\begin{split}
G_{i+\frac{1}{2}}-G_{i-\frac{1}{2}} \leq& \left(\Pi'\left(\rho^-_{i+\frac{1}{2}}\right)-\frac{1}{2}u_i^2+H_{i+\frac{1}{2}}\right)\mathscr{F}^\rho\left(U^-_{i+\frac{1}{2}}, U^+_{i+\frac{1}{2}}\right)\\
 & -\left(\Pi'\left(\rho^+_{i-\frac{1}{2}}\right)-\frac{1}{2}u_i^2+H_{i-\frac{1}{2}}\right)\mathscr{F}^\rho\left(U^-_{i-\frac{1}{2}}, U^+_{i-\frac{1}{2}}\right)\\
  & +u_i\left(\mathscr{F}^{\rho u}\left(U^-_{i+\frac{1}{2}}, U^+_{i+\frac{1}{2}}\right)-\mathscr{F}^{\rho u}\left(U^-_{i-\frac{1}{2}}, U^+_{i-\frac{1}{2}}\right)+P\left(\rho^+_{i-\frac{1}{2}}\right)-P\left(\rho^-_{i+\frac{1}{2}}\right)\right).
\end{split}
\end{equation*}
From here, by bearing in mind the definition of $\rho^-_{i+\frac{1}{2}}$ and $\rho^+_{i-\frac{1}{2}}$ in \eqref{eq:rhointerface} and the definition of the scheme in \eqref{eq:fvbasic}-\eqref{eq:numfluxnew}-\eqref{eq:sourcewellbalanced}, we get
\begin{equation*}
\begin{split}
G_{i+\frac{1}{2}}-G_{i-\frac{1}{2}} \leq &\ \left(\Pi'\left(\rho_i\right)-\frac{1}{2}u_i^2+H_{i}\right)\left(\mathscr{F}^\rho\left(U^-_{i+\frac{1}{2}}, U^+_{i+\frac{1}{2}}\right)-\mathscr{F}^\rho\left(U^-_{i-\frac{1}{2}}, U^+_{i-\frac{1}{2}}\right)\right)\\
& +u_i\left(\mathscr{F}^{\rho u}\left(U^-_{i+\frac{1}{2}}, U^+_{i+\frac{1}{2}}\right)-\mathscr{F}^{\rho u}\left(U^-_{i-\frac{1}{2}}, U^+_{i-\frac{1}{2}}\right)+P\left(\rho^+_{i-\frac{1}{2}}\right)-P\left(\rho^-_{i+\frac{1}{2}}\right)\right)\\
= &\ -\left(\Pi'\left(\rho_i\right)-\frac{1}{2}u_i^2+H_{i}\right) \Delta x_i \frac{d\rho_i}{dt}-\Delta x_i u_i \frac{d}{dt}(\rho_i u_i)\\
&\, -u_i\left(\gamma \Delta x_i \rho_i u_i+\Delta x_i \rho_i \sum_j \rho_j \left(u_i-u_j\right)\psi_{ij}\right).\\
\end{split}
\end{equation*}
Finally, this last inequality results in the desired cell entropy inequality \eqref{eq:cellentropyineq2} by rearranging according to \eqref{eq:compactsys}, yielding
\begin{equation}\label{aux}
\Delta x_i \frac{d\eta_i}{dt}+\Delta x_i H_i \frac{d\rho_i}{dt}+G_{i+\frac{1}{2}}-G_{i-\frac{1}{2}}=-u_i\left(\gamma \Delta x_i \rho_i u_i+\Delta x_i \rho_i \sum_j \rho_j \left(u_i-u_j\right)\psi_{ij}\right).
\end{equation}

\item The last property of the scheme and formulas
    \eqref{eq:disfreeenergydiscrete}-\eqref{eq:freeenergydiscrete} follow
    by summing over the index $i$ over identity \eqref{aux}, collecting
    terms and symmetrizing the dissipation using the symmetry of $\psi$.

\item Starting from the finite volume equation for the density in \eqref{eq:compactsys},
\begin{equation*}
\Delta x_i \frac{d\rho_i}{dt} =-\mathscr{F}^\rho\left(U^-_{i+\frac{1}{2}}, U^+_{i+\frac{1}{2}}\right)+\mathscr{F}^\rho\left(U^-_{i-\frac{1}{2}}, U^+_{i-\frac{1}{2}}\right),
\end{equation*}
one can multiply it by $x_i$ and sum it over the index $i$, resulting in
\begin{equation*}
\frac{d}{dt}\left(\sum_i \Delta x_i \rho_i x_i  \right)=\sum_i x_i \left(-\mathscr{F}^\rho\left(U^-_{i+\frac{1}{2}}, U^+_{i+\frac{1}{2}}\right)+\mathscr{F}^\rho\left(U^-_{i-\frac{1}{2}}, U^+_{i-\frac{1}{2}}\right)\right).
\end{equation*}
By rearranging and considering, for instance, periodic or no flux boundary
conditions, we get~\eqref{eq:proofevolcentremass}.

On the other hand, the finite volume equation for the momentum in \eqref{eq:compactsys}, after summing over the index $i$, becomes
\begin{equation}\label{eq:proofmomentum}
\begin{gathered}
\frac{d}{dt}\left(\sum_i \Delta x_i \rho_i u_i \right)=\sum_i \left( P \left(\rho^-_{i+\frac{1}{2}}\right)-P \left(\rho^+_{i-\frac{1}{2}}\right)\right)-\gamma  \sum_i \Delta x_i \rho_i u_i\\ -\sum_{i,j} \Delta x_i \Delta x_j  \rho_i \rho_j (u_i-u_j) \psi_{ij} ,
\end{gathered}
\end{equation}
since the numerical fluxes cancel out due to the sum over the index $i$. In addition, the Cucker-Smale damping term also vanishes due to the symmetry in $\psi(x)$. Finally, if the initial density is symmetric and the initial velocity antisymmetric, the sum of pressures in the RHS of \eqref{eq:proofmomentum} is $0$, due to the symmetry in the density. This implies that the discrete solution for the density and momentum maintains those symmetries, since \eqref{eq:proofmomentum} is simplified as
\begin{equation*}
\sum_i \Delta x_i \rho_i u_i =0
\end{equation*}
and as a result \eqref{eq:proofevolcentremass} reduces to \eqref{eq:proofevolcentremasssimp}. This means that the discrete centre of mass is conserved in time and is centred at $0$, for initial symmetric densities and initial antisymmetric velocities.
\end{enumerate}
\end{proof}

\begin{remark}
As a consequence of the previous proofs, our scheme conserves all the
structural properties of the hydrodynamic system \eqref{eq:generalsys2} at
the semidiscrete level including the dissipation of the discrete free energy
\eqref{eq:equalenergy} and the characterization of the steady states. These
properties are analogous to those obtained for finite volume schemes in the
overdamped limit \cite{carrillo2015finite,sun2018discontinuous}.
\end{remark}

\begin{remark}
All the previous properties, which are applicable for free energies of the form \eqref{eq:freeenergy}, can be extended to the general free energies in \eqref{eq:freeenergygeneral}. It can be shown indeed that the discrete analog of the free energy dissipation in \eqref{eq:disfreeenergydiscrete} still holds for a discrete total energy defined as in \eqref{eq:freeenergydiscrete} and a discrete free energy of the form
\begin{equation}
\mathcal{F}^\Delta = \sum_i \Delta x_i \left[\Pi\left(\rho_i\right)+ V_i\rho_i \right]+\frac12 \sum_{i} \Delta x_i \rho_i K_i,
\end{equation}
where $K_i$ is a discrete approximation of $K(W(x)\star\rho)$ at the node $x_i$ and is evaluated as
\begin{equation}
K_i=K\left(\sum_{j}\Delta x_j W_{ij} \rho_j\right).
\end{equation}
\end{remark}

%
%
\subsection{Second-order extension}\label{subsec:secondorder}

The usual procedure to extend a first-order scheme to second order is by computing the numerical fluxes \eqref{eq:numflux} from limited reconstructed values of the density and momentum at each side of the boundary, contrary to the cell-centred values taken for the first order schemes \eqref{eq:numflux1}. These values are classically computed in three steps: prediction of the gradients in each cell, linear extrapolation and limiting procedure to preserve nonnegativity. For instance, MUSCL \cite{osher1985convergence} is a usual reconstruction procedure following these steps. From here the values $\rho_{i,l}$, $\rho_{i,r}$, $u_{i,l}$ and $u_{i,r}$ are obtained $\forall i$, where $l$ indicates at the left of the boundary and $r$ at the right. Then the inputs for the numerical flux in \eqref{eq:numflux}, for a usual second-order scheme, are
\begin{equation*}
F_{i+\frac{1}{2}}=\mathscr{F}\left(U_{i,r},U_{i+1,l}\right).
\end{equation*}
This procedure has already been adapted to satisfy the well-balanced property
and maintain the second order for specific applications, such as shallow
water \cite{audusse2004fast} or chemotaxis \cite{filbet2005approximation}. In
this subsection the objective is to extend the procedure to general free
energies of the form \eqref{eq:freeenergy}. As it happened for the
well-balanced first-order scheme, the boundary values introduced in the
numerical flux, which in this case are $U_{i,r}$ and $U_{i+1,l}$, need to be
adapted to satisfy the well-balanced property.

For the well-balanced scheme the first step is to reconstruct the boundary values $\rho_{i,l}$, $\rho_{i,r}$, $u_{i,l}$ and $u_{i,r}$ following the three mentioned steps. In addition, the reconstructed values of the potential $H(x,\rho)$ at the boundaries, $H_{i,l}$ and $H_{i,r}$ $\forall i$, have to be also computed. This is done as suggested in \cite{audusse2004fast}. Instead of reconstructing directly $H_{i,l}$ and $H_{i,r}$ following the three mentioned steps, for certain applications one has to reconstruct firstly $\left(\Pi'(\rho)+H(x,\rho)\right)_i$ to obtain $\left(\Pi'(\rho)+H(x,\rho)\right)_{i,l}$ and $\left(\Pi'(\rho)+H(x,\rho)\right)_{i,r}$, and subsequently compute $H_{i,l}$ and $H_{i,r}$ as
\red{\begin{equation*}
\begin{gathered}
 H_{i,l}=\left(\Pi'(\rho)+H(x,\rho)\right)_{i,l}-\Pi'\left(\rho_{i,l}\right),\\
  H_{i,r}=\left(\Pi'(\rho)+H(x,\rho)\right)_{i,r}-\Pi'\left(\rho_{i,r}\right).
\end{gathered}
\end{equation*}}
This is shown in \cite{audusse2004fast} to be necessary in order to  maintain nonnegativity and the steady state in applications where there is an interface between dry and wet cells. For instance, these interfaces appear when considering pressures of the form $P=\rho^m$ with $m>0$, as it is shown in examples \ref{ex:squarepot} and \ref{ex:squarepotdoublewell} of section \ref{sec:numtest}. For other applications where vacuum regions do not occur, the values $H_{i,l}$ and $H_{i,r}$ can be directly reconstructed following the three mentioned steps.

After this first step, the inputs for the numerical flux are updated from \eqref{eq:numflux} to satisfy the well-balanced property as
\begin{equation*}
F_{i+\frac{1}{2}}=\mathscr{F}\left(U_{i+\frac{1}{2}}^-,U_{i+\frac{1}{2}}^+\right), \ \text{where}\quad U_{i+\frac{1}{2}}^{-}=\begin{pmatrix}
\rho_{i+\frac{1}{2}}^{-} \\ \rho_{i+\frac{1}{2}}^{-} u_{i,r}
\end{pmatrix},\quad
U_{i+\frac{1}{2}}^{+}=\begin{pmatrix}
\rho_{i+\frac{1}{2}}^{+} \\ \rho_{i+\frac{1}{2}}^{+} u_{i+1,l}
\end{pmatrix}.
\end{equation*}
The interface values $\rho_{i+\frac{1}{2}}^{\pm}$ are reconstructed as in the first-order scheme, by taking into account the steady state relation in \eqref{eq:steadyvarenerdiscrete}.  The application of \eqref{eq:steadyvarenerdiscrete} to the cells with centred nodes $x_i$ and $x_{i+1}$ leads to
\begin{equation*}
\begin{gathered}
  \Pi'\left(\rho_{i+\frac{1}{2}}^{-}\right)+H_{i+\frac{1}{2}}=\Pi'\left(\rho_{i,r}\right)+H_{i,r},\\
  \Pi'\left(\rho_{i+\frac{1}{2}}^{+}\right)+H_{i+\frac{1}{2}}=\Pi'\left(\rho_{i+1,l}\right)+H_{i+1,l},
\end{gathered}
\end{equation*}
where the term $H_{i+\frac{1}{2}}$ is evaluated to preserve consistency and stability, with an upwind or average value obtained as
\begin{equation*}
H_{i+\frac{1}{2}}=\max\left(H_{i,r},H_{i+1,l}\right)\quad \textrm{or} \quad H_{i+\frac{1}{2}}=\frac{1}{2}\left(H_{i,r}+H_{i+1,l}\right).
\end{equation*}
Then, by denoting as $\xi(x)$ the inverse function of $\Pi'(x)$, the interface values $\rho_{i+\frac{1}{2}}^{\pm}$ are computed as
\begin{equation*}
\begin{gathered}
\rho_{i+\frac{1}{2}}^{-}=\xi \left(\Pi'\left(\rho_{i,r}\right)+H_{i,r}-H_{i+\frac{1}{2}}\right),\\
\rho_{i+\frac{1}{2}}^{+}=\xi \left(\Pi'\left(\rho_{i+1,l}\right)+H_{i+1,l}-H_{i+\frac{1}{2}}\right).
\end{gathered}
\end{equation*}
The source term is again distributed along the interfaces,
\begin{equation*}
S_i=S_{i+\frac{1}{2}}^{-}+S_{i-\frac{1}{2}}^{+}+S_i^c,
\end{equation*}
where
\begin{equation*}
S_{i+\frac{1}{2}}^{-}=\frac{1}{\Delta x_i}\begin{pmatrix}
0 \\  P\left(\rho_{i+\frac{1}{2}}^{-}\right) -P\left(\rho_{i,r}\right)
\end{pmatrix}, \quad
S_{i-\frac{1}{2}}^{+}=\frac{1}{\Delta x_i}\begin{pmatrix}
0 \\  P\left(\rho_{i,l}\right)- P\left(\rho_{i-\frac{1}{2}}^{+}\right)
\end{pmatrix}.
\end{equation*}
\red{The inclusion of the central source term $S_i^c$ is vital in order to
preserve the second-order accuracy and well-balanced property of the scheme.
This idea was firstly introduced in \cite{katsaounis2005first}, where second
order error estimates are derived under certain conditions for $S_i^c$.
Further works  customize this central source term $S_i^c$ for particular
applications such as shallow water equations
\cite{katsaounis2003second,audusse2004fast} or chemotaxis
\cite{filbet2005approximation}. There is some flexibility in the choice of
this term, as far as it satisfies two criteria for second-order accuracy and
well-balancing. In the following remark we summarize the two criteria, which
are described with more extend in Ref.~\cite{bouchut2004nonlinear}
(specifically, (4.187) for second-order accuracy, and (4.204) for
well-balancing).}

\begin{remark}
\red{The central source term $S_i^c$ preserves the second-order accuracy and well-balanced property of the scheme if the following two criteria are satisfied:
\begin{enumerate}[label=(\roman*)]
\item Second-order accuracy if
\begin{equation}\label{eq:criteriaSc1}
S_i^c\lt(\rho_{i,l},\rho_{i,r},H_{i,l},H_{i,r}\rt)=\begin{pmatrix}
0 \\ \lt(-\frac{\rho_{i,l}+\rho_{i,r}}{2}+\mathcal{O}\lt(\lt|\rho_{i,r}-\rho_{i,l}\rt|^2+\lt|H_{i,r}-H_{i,l}\rt|^2\rt)\rt)\lt(H_{i,r}-H_{i,l}\rt)
\end{pmatrix}
\end{equation}
as $\rho_{i,r}-\rho_{i,l}\rightarrow 0 $ and $H_{i,r}-H_{i,l}\rightarrow 0 $.
\item Well-balanced property if
\begin{equation}\label{eq:criteriaSc2}
S_i^c\lt(\rho_{i,l},\rho_{i,r},H_{i,l},H_{i,r}\rt)=F\lt(\rho_{i,r},H_{i,r}\rt)-F\lt(\rho_{i,l},H_{i,l}\rt),
\end{equation}
meaning that the steady states are let invariant.
\end{enumerate}}
\end{remark}

The objective here is to provide a general form of $S_i^c$ which applies to general free energies of the form \eqref{eq:freeenergy}. Following the strategy in \cite{bouchut2004nonlinear}, we propose to approximate the generalized centred sources as
\begin{equation*}
S_i^c=\frac{1}{\Delta x_i}\begin{pmatrix}
0 \\ P(\rho_{i,r})-P(\rho_{i,r}^*)-P(\rho_{i,l})+P(\rho_{i,l}^*)
\end{pmatrix}-\begin{pmatrix}
0 \\ \gamma \rho_i u_i + \rho_i \sum_{j} (u_i-u_j)\rho_j \psi(x_i-x_j)
\end{pmatrix},
\end{equation*}
where the values $\rho_{i,l}^*$ and $\rho_{i,r}^*$ are computed from the steady state relation \eqref{eq:steadyvarenerdiscrete} as
\begin{equation*}
\begin{gathered}
\rho_{i,l}^*=\xi \left(\Pi'\left(\rho_{i,l}\right)+H_{i,l}-H_i^*\right),\\
\rho_{i,r}^*=\xi \left(\Pi'\left(\rho_{i,r}\right)+H_{i,r}-H_i^*\right),
\end{gathered}
\end{equation*}
and $H_i^*$ is a centred approximation of the potentials satisfying
\begin{equation*}
H_i^*=\frac{1}{2}(H_{i,l}+H_{i,r}).
\end{equation*}
\red{The proposed structure of $S_i^c$ is suggested in \cite{bouchut2004nonlinear} and satisfies the two criteria for second-order accuracy \eqref{eq:criteriaSc1} and well-balanced property \eqref{eq:criteriaSc2}.}

Overall, the second-order semidiscrete scheme defined in \eqref{eq:fvbasic} and constructed as detailed in this subsection \ref{subsec:secondorder}, and for a numerical flux $\mathscr{F}$ satisfying the properties stated in the introduction of section \ref{sec:numsch}, satisfies:
\begin{enumerate}[label=(\roman*)]
\item preservation of the nonnegativity of $\rho_i(t)$;
\item well-balanced property, thus preserving the steady states given by \eqref{eq:steadyvarenerdiscrete};
\item consistency with the system \eqref{eq:generalsys2};
\item second-order accuracy.
\end{enumerate}
The proof of these properties is omitted here since it follows the same techniques from \cite{audusse2004fast,filbet2005approximation}, and the general procedure is very similar to the one from the first-order scheme in subsection \ref{subsec:firstorderprop}.

%
%
%
%
\section{Numerical tests}\label{sec:numtest}

This section details numerical simulations in which the first- and
second-order schemes from section \ref{sec:numsch} are employed. Firstly,
subsection \ref{subsec:val} contains the validation of the first- and
second-order schemes: the well-balanced property and the order of accuracy of
the schemes are tested in four different configurations. Secondly, subsection
\ref{subsec:numexp} illustrates the application of the numerical schemes to a
variety of choices of the free energy, leading to interesting numerical
experiments for which analytical results are limited in the literature.

Unless otherwise stated, all simulations contain linear damping with
$\gamma=1$ and have a total unitary mass. Only the indicated ones contain the
Cucker-Smale damping term, where the communication function satisfies
\begin{equation*}
\psi(x)=\frac{1}{\left(1+|x|^2\right)^\frac{1}{4}}.
\end{equation*}

The pressure function in the simulations has the form of $P(\rho)=\rho^m$,
with $m\geq1$. When $m=1$ the pressure satisfies the ideal-gas relation
$P(\rho)=\rho$, and the density does not develop vacuum regions during the
temporal evolution. For this case the employed numerical flux is the
versatile local Lax-Friedrich flux. For the simulations where
$P(\rho)=\rho^m$ and $m>1$ vacuum regions with $\rho=0$ are generated. This
implies that the hyperbolicity of the system \eqref{eq:generalsys2} is lost
in those regions, and the local Lax-Friedrich scheme fails. As a result, an
appropiate numerical flux has to be implemented to handle the vacuum regions.
In this case a kinetic solver based on \cite{perthame2001kinetic}, and
already implemented in previous works \cite{audusse2005well}, is employed.

The time discretization is acomplished by means of the third order \red{TVD}
Runge-Kutta method \cite{gottlieb1998total} and the CFL number is taken as
$0.7$ in all the simulations. The boundary conditions are chosen to be no
flux. For more details about the numerical fluxes, temporal discretization, \red{boundary conditions}
and CFL number, we remit the reader to Appendix \ref{app:numerics}.

Videos from all the simulations displayed in this work are available at  \cite{simulations}.
%
%
\subsection{Validation of the numerical scheme}\label{subsec:val}

The validation of the schemes from section \ref{sec:numsch} includes a test for the well-balanced property and a test for the order of accuracy \blue{in the transient regimes}. These tests are completed in four different examples with steady states satisfying \eqref{eq:steadyvarener}, which differ in the choice of the free energy, potentials and the inclusion of Cucker-Smale damping terms. An additional fifth example presenting moving steady states of the form \eqref{eq:steadyvarener2} is considered to show that our schemes satisfy the order of accuracy test even for this challenging steady states.

The well-balanced property test evaluates whether the steady state solution is preserved in time up to machine precision. As a result, the initial condition of the simulation has to be directly the steady state. The results of this test for the four examples of this section are presented in table \ref{table:preserve}. All the simulations are run from $t=0$ to $t=5$, \red{and the number of cells is 50}.

\begin{table}[!htbp]
\centering
\caption{Preservation of the steady state for the examples \ref{ex:idealpot}, \ref{ex:idealpotCS}, \ref{ex:idealker} and \ref{ex:squarepot} with the first- and second-order schemes and double precision, at $t=5$}\label{table:preserve}

\label{table:cproperty}
\begin{tabular}{cccc} \hline
                         & Order of the scheme & $L^1$ error & $L^{\infty}$ error \\
                         \hline
\multirow{2}{*}{Example \ref{ex:idealpot}} & 1\textsuperscript{st}            & 9.1012E-18  & 1.1102E-16     \\
                         & 2\textsuperscript{nd}         & 2.3191E-17  & 2.2843E-16     \\
                         \hline
\multirow{2}{*}{Example \ref{ex:idealpotCS}} & 1\textsuperscript{st}            & 7.8666E-18  & 1.1102E-16     \\
                         & 2\textsuperscript{nd}             & 1.4975E-17  & 1.5057E-16     \\
                         \hline
\multirow{2}{*}{Example \ref{ex:idealker}} & 1\textsuperscript{st}            & 5.5020E-17  & 6.6613E-16     \\
                         & 2\textsuperscript{nd}             & 6.4514E-17  & 7.2164E-16     \\
                         \hline
\multirow{2}{*}{Example \ref{ex:squarepot}} & 1\textsuperscript{st}            &    1.3728E-17         &       2.2204E-16         \\
                         & 2\textsuperscript{nd}             &   3.4478E-18          &      1.1102E-16          \\
                         \hline
\end{tabular}
\end{table}

The order of accuracy \blue{in the transient regimes} test is based on evaluating the $L^1$ error of a numerical solution for a particular choice of $\Delta x$ with respect to a reference solution, \blue{and for a time when the steady state is not reached yet}. Subsequent $L^1$ errors are obtained after halving the $\Delta x$ of the previous numerical solution, doubling in this way the total number of cells. The order of the scheme is then computed as
\begin{equation}
\text{Order of the scheme}=\ln_2\left(\frac{L^1\,\text{error} (\Delta x)}{L^1\,\text{error} (\Delta x/2)}\right),
\end{equation}
and the $\Delta x$ is halved four times.

The reference solution is frequently taken as an explicit solution of the
system that is being tested. In this case, the system in
\eqref{eq:generalsys2} does not have an explicit solution in time for the
free energies presented here, even though the steady solution can be
analytically computed. Since we are interested in evaluating the order of
accuracy away from equilibrium, the reference solution is computed from the
same numerical scheme but with a really small $\Delta x$, so that the numerical solution can be considered as the exact one. In all cases here the reference solution is obtained from a mesh with 25600 cells, while the numerical solutions employ a number of cells between $50$ and $400$.

The results from the accuracy tests are shown in the tables
\ref{table:idealpot}, \ref{table:idealpotCS}, \ref{table:idealker},
\ref{table:squarepot} and \ref{table:movingwater}. The simulations were run
with the configurations specified in each example and from $t=0$ to $t=0.3$,
unless otherwise stated. \blue{The final time of $t=0.3$ is taken so that all
examples are in the transient regime}.


\begin{examplecase}[Ideal-gas pressure and attractive potential]\label{ex:idealpot} In this example the pressure satisfies $P(\rho)=\rho$ and there is an external potential of the form $V(x)=\frac{x^2}{2}$. As a result, the relation holding in the steady state is
\begin{equation}\label{eq:constantidealpot}
\frac{\delta \mathcal{F}}{\delta \rho}=\Pi'(\rho)+H=\ln(\rho)+\frac{x^2}{2}=\text{constant}\ \text{on}\ \mathrm{supp}(\rho)\ \text{and}\ u=0.
\end{equation}
The steady state, for an initial mass $M_0$, explicitly satisfies
\begin{equation}\label{eq:steadyidealpot}
\rho_{\infty}=M_0 \frac{e^{-x^2/2}}{\int_\R e^{-x^2/2} dx}.
\end{equation}
For the order of accuracy test the initial conditions are
\begin{equation}\label{eq:icidealpot}
\rho(x,t=0)=M_0\frac{0.2+5\,\cos\left(\frac{\pi x}{10}\right)}{\int_\R\left(0.2+5\,\cos\left(\frac{\pi x}{10}\right)\right)dx},\quad \rho u(x,t=0)=-0.05 \sin\left(\frac{\pi x}{10}\right),\quad x\in [-5,5],
\end{equation}
with $M_0$ equals to $1$ so that the total mass is unitary.
The order of accuracy test from this example is shown in table \ref{table:idealpot}, and the evolution of the density, momentum, variation of the free energy with respect to the density, total energy and free energy are depicted in figure \ref{fig:idealpot}. From \ref{fig:idealpot} (D) one can notice how the discrete total energy always decreases in time, due to the discrete free energy dissipation property \eqref{eq:disfreeenergydiscrete}, and how there is an exchange between free energy and kinetic energy which makes the discrete free energy plot oscillate.

\begin{figure}[ht!]
\begin{center}
\subfloat[Evolution of the density]{\protect\protect\includegraphics[scale=0.4]{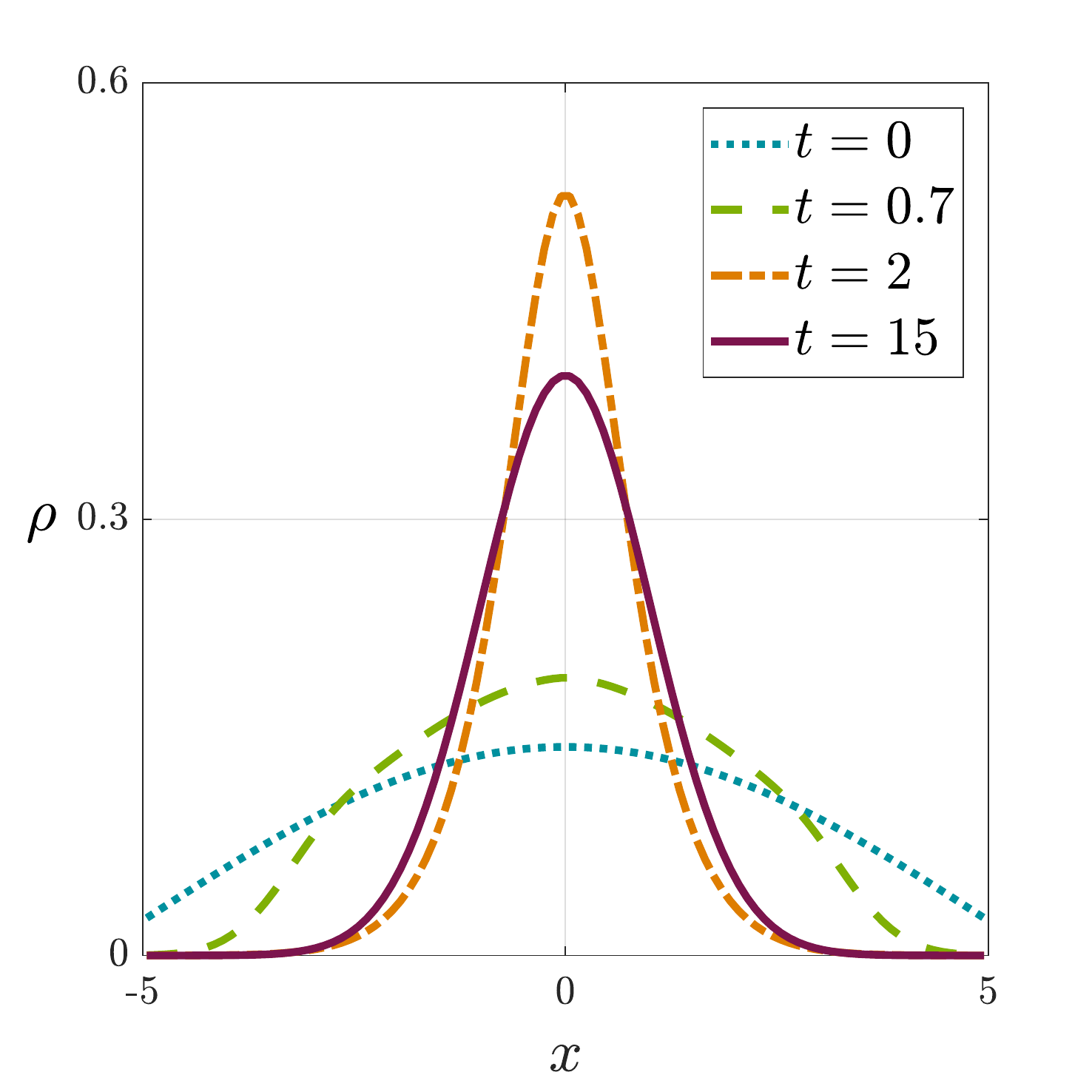}
}
\subfloat[Evolution of the momentum]{\protect\protect\includegraphics[scale=0.4]{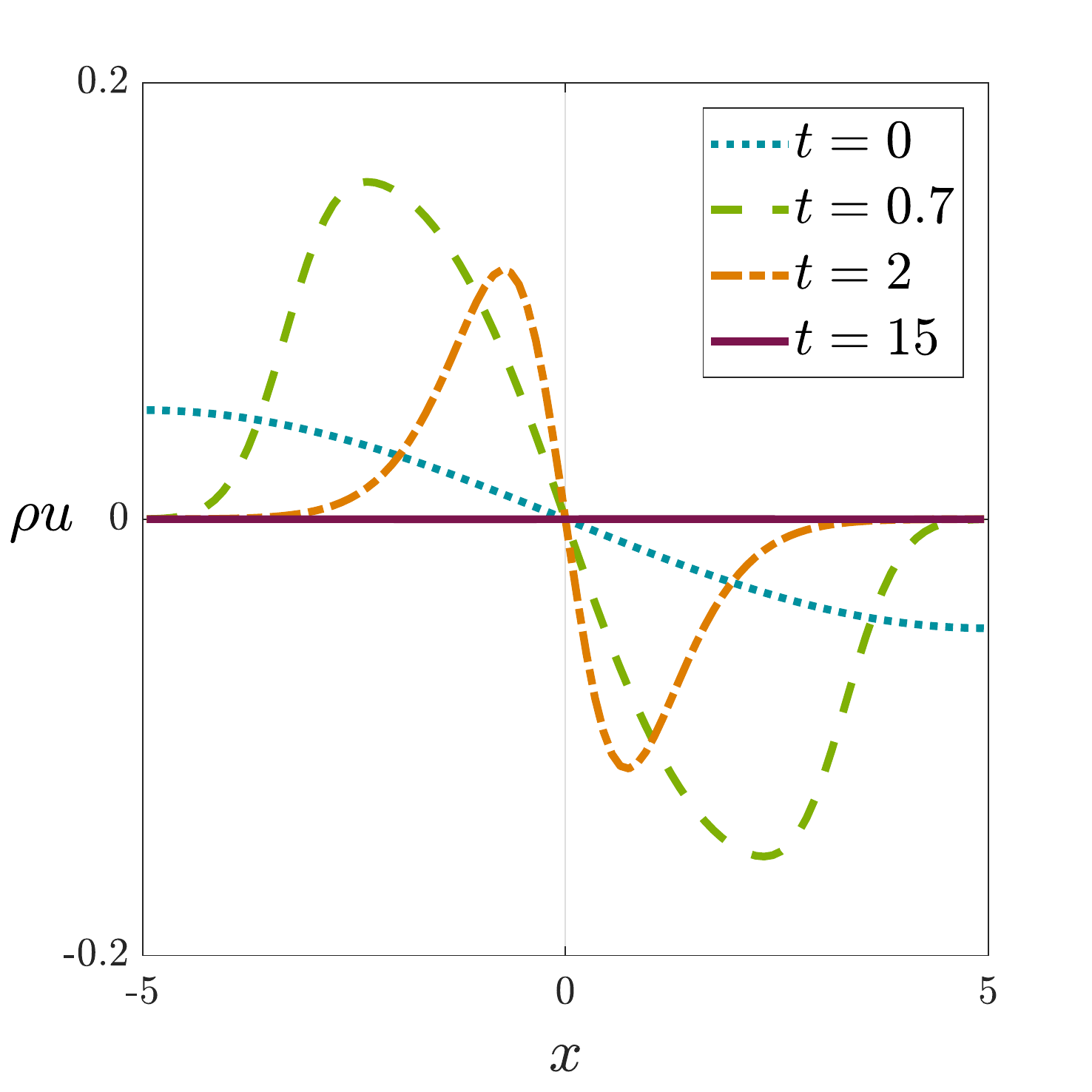}
}\\
\subfloat[Evolution of the variation of the free energy]{\protect\protect\includegraphics[scale=0.4]{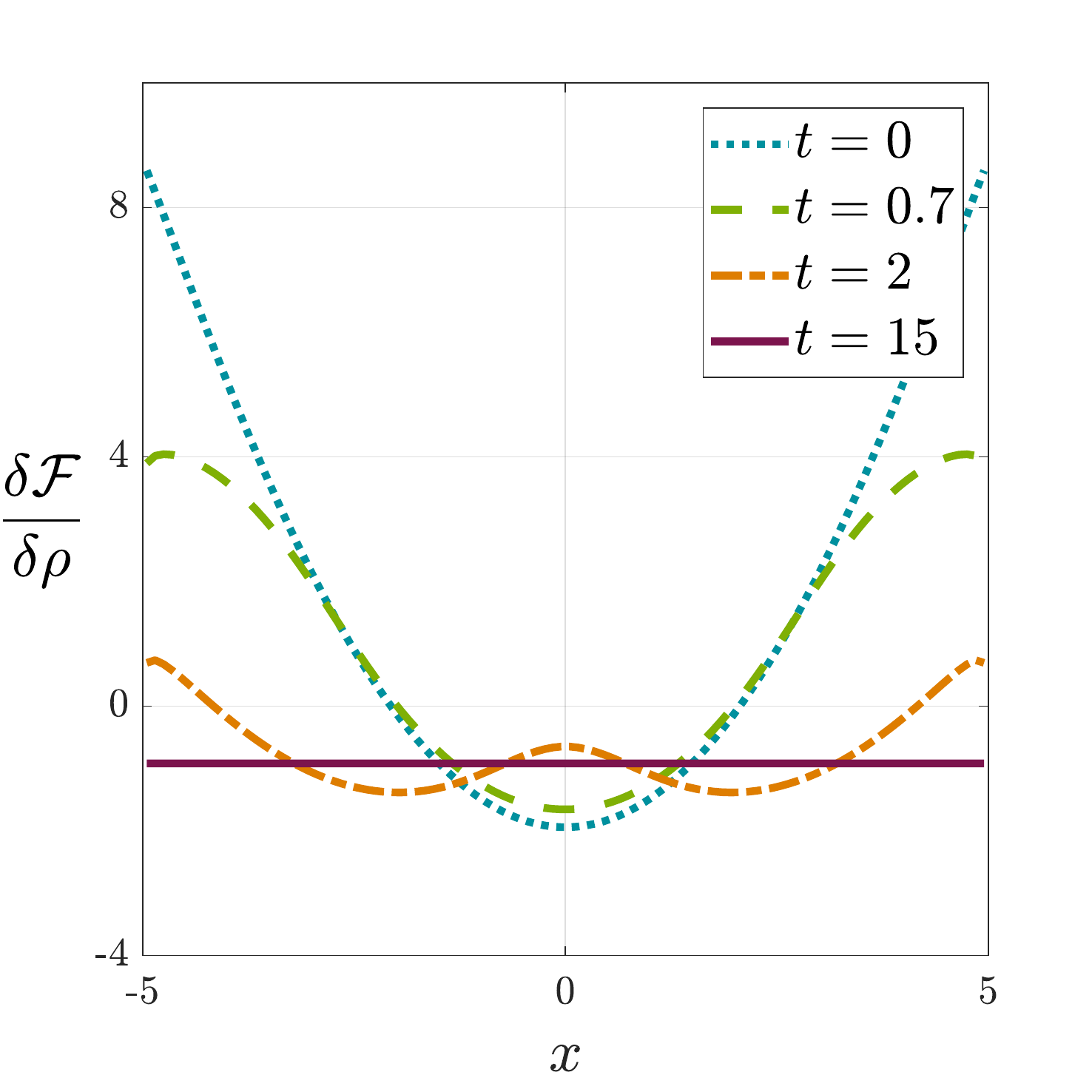}
}
\subfloat[Evolution of the total energy and free energy]{\protect\protect\includegraphics[scale=0.4]{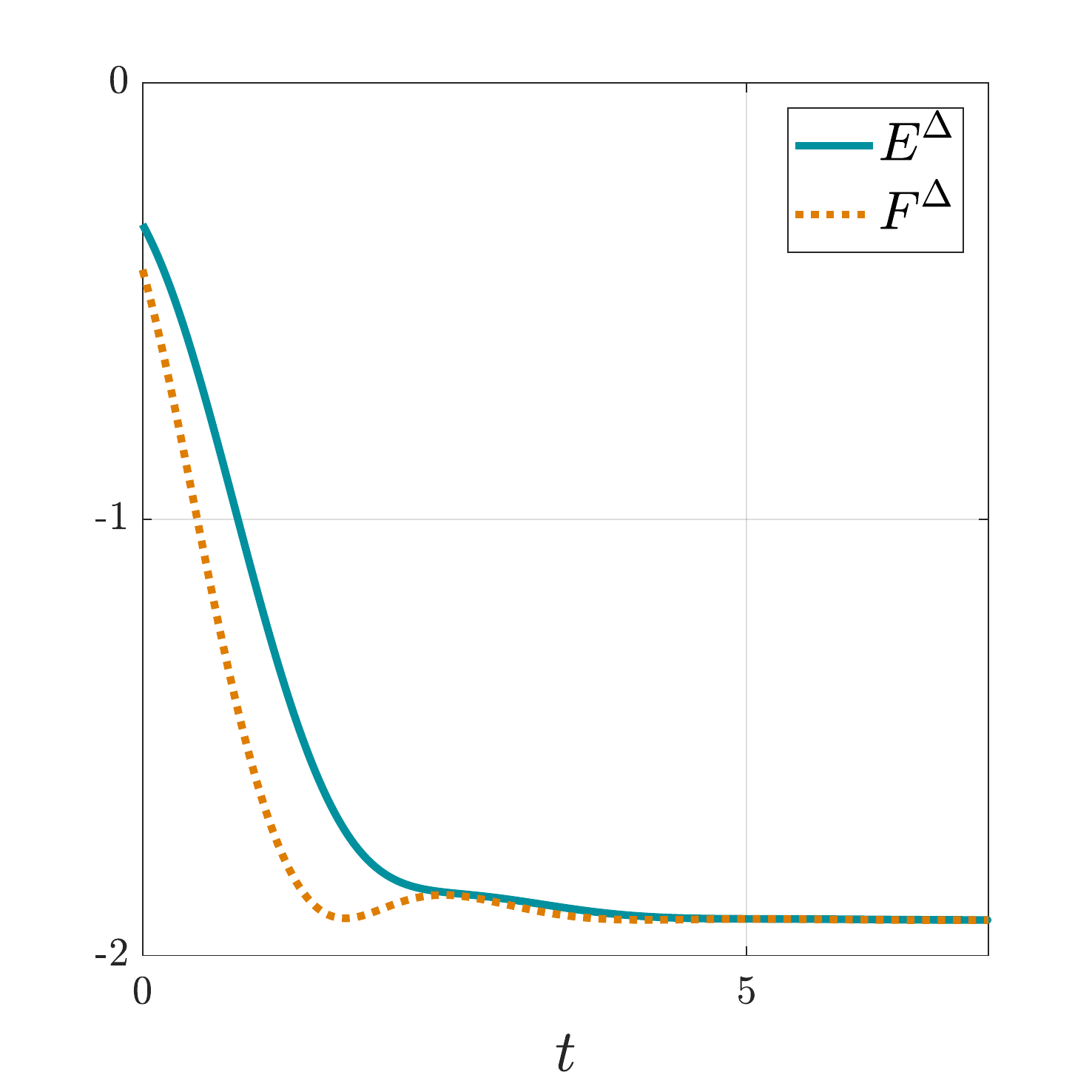}
}
\end{center}
\protect\protect\caption{\label{fig:idealpot} Temporal evolution of Example \ref{ex:idealpot}.}
\end{figure}

\begin{table}[!htbp]
\centering
\caption{Accuracy test for Example \ref{ex:idealpot} with the first and second-order schemes\blue{, at $t=0.3$}}
\label{table:idealpot}
\begin{tabular}{c c c c c }
\hline
\multirow{2}{*}{\begin{tabular}[c]{@{}c@{}}Number of \\ cells\end{tabular}} & \multicolumn{2}{c}{ First-order} & \multicolumn{2}{c}{Second-order}   \\ \cline{2-5}
                                                                                & $L^1$ error    & order    & $L^1$ error           & order          \\ \hline
50                                                                              &               6.8797E-03 &     -     &      7.6166E-04                &     -             \\ \hline
100                                                                             &               3.4068E-03 &     1.01     &        2.0206E-04               &           1.91        \\ \hline
200                                                                             &               1.6826E-03 &    1.02      &        5.0308E-05               &         2.01         \\ \hline
400                                                                             &               8.3104E-04 &    1.02      &        1.2879E-05               &         1.97        \\ \hline
\end{tabular}

\end{table}

\end{examplecase}


\begin{examplecase}[Ideal-gas pressure, attractive potential and Cucker-Smale damping terms]\label{ex:idealpotCS}
In this example the pressure satisfies $P(\rho)=\rho$ and there is an
external potential of the form $V(x)=\frac{x^2}{2}$. The difference with
example \ref{ex:idealpot} is that the Cucker-Smale damping terms are
included, and the linear damping term $-\rho u$ excluded.

The relation holding in the steady state is expressed in \eqref{eq:constantidealpot} and the steady state satisfies \eqref{eq:steadyidealpot}. The initial conditions are also \eqref{eq:icidealpot}. The order of accuracy test from this example is shown in table \ref{table:idealpotCS}, and the evolution of the density, momentum, variation of the free energy with respect to the density, total energy and free energy are depicted in figure \ref{fig:idealpotCS}. The lack of linear damping leads to higher oscillations in the momentum plots in comparison to figure \ref{fig:idealpot}. There is also an exchange of kinetic and free energy during the temporal evolution, which could be noticed from the oscillations of the discrete free energy in figure \ref{fig:idealpotCS} (D).

\begin{figure}[ht!]
\begin{center}
\subfloat[Evolution of the density]{\protect\protect\includegraphics[scale=0.4]{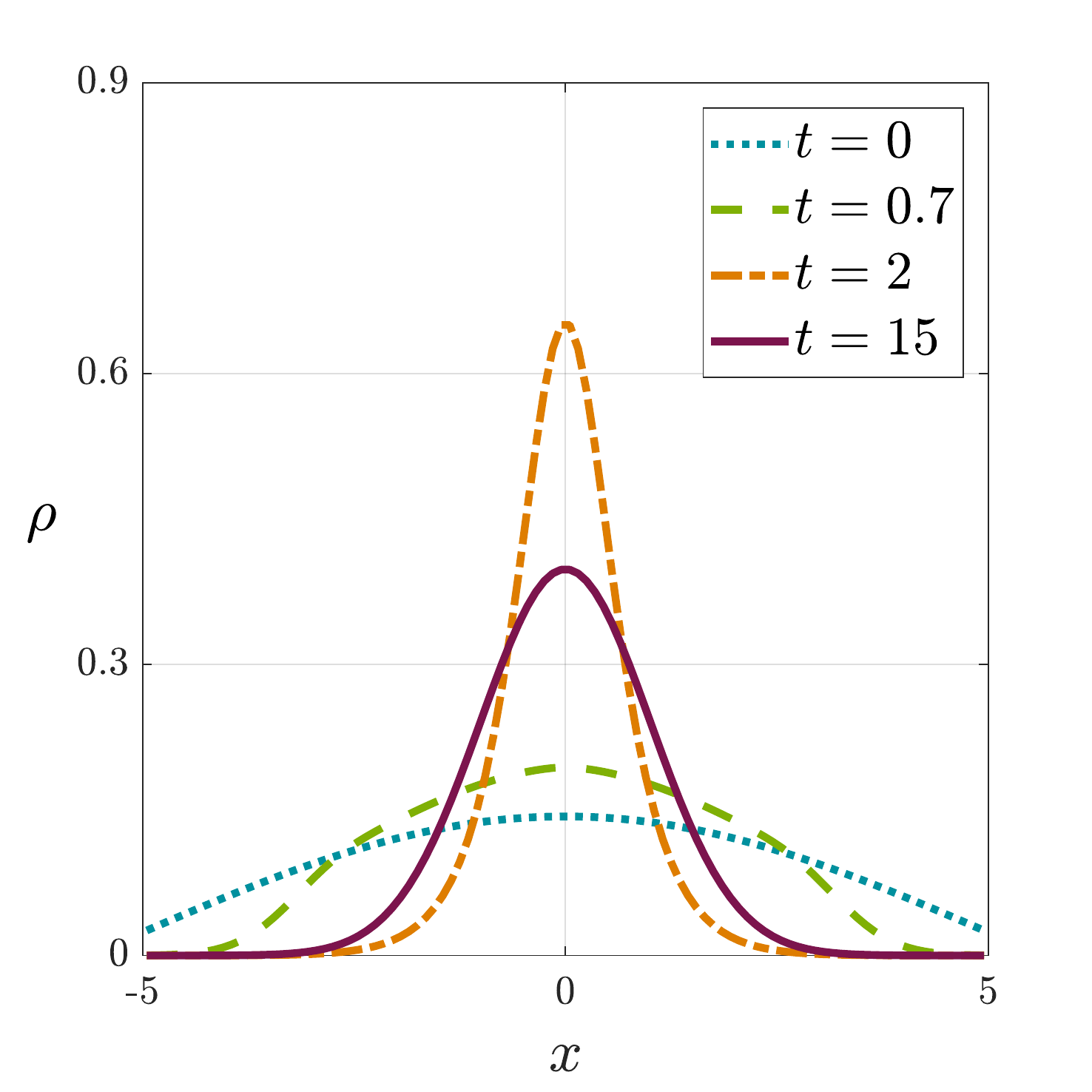}
}
\subfloat[Evolution of the momentum]{\protect\protect\includegraphics[scale=0.4]{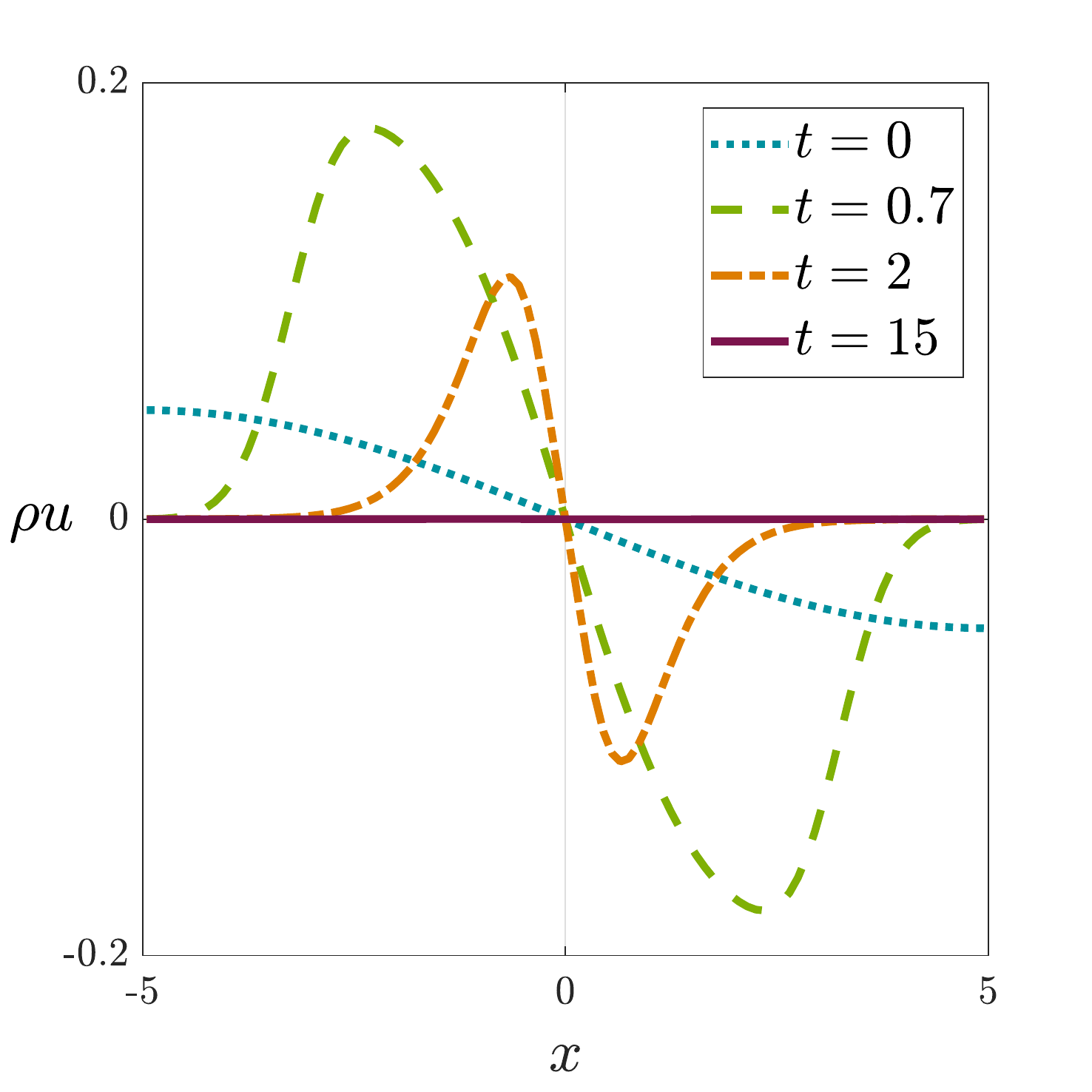}
}\\
\subfloat[Evolution of the variation of the free energy]{\protect\protect\includegraphics[scale=0.4]{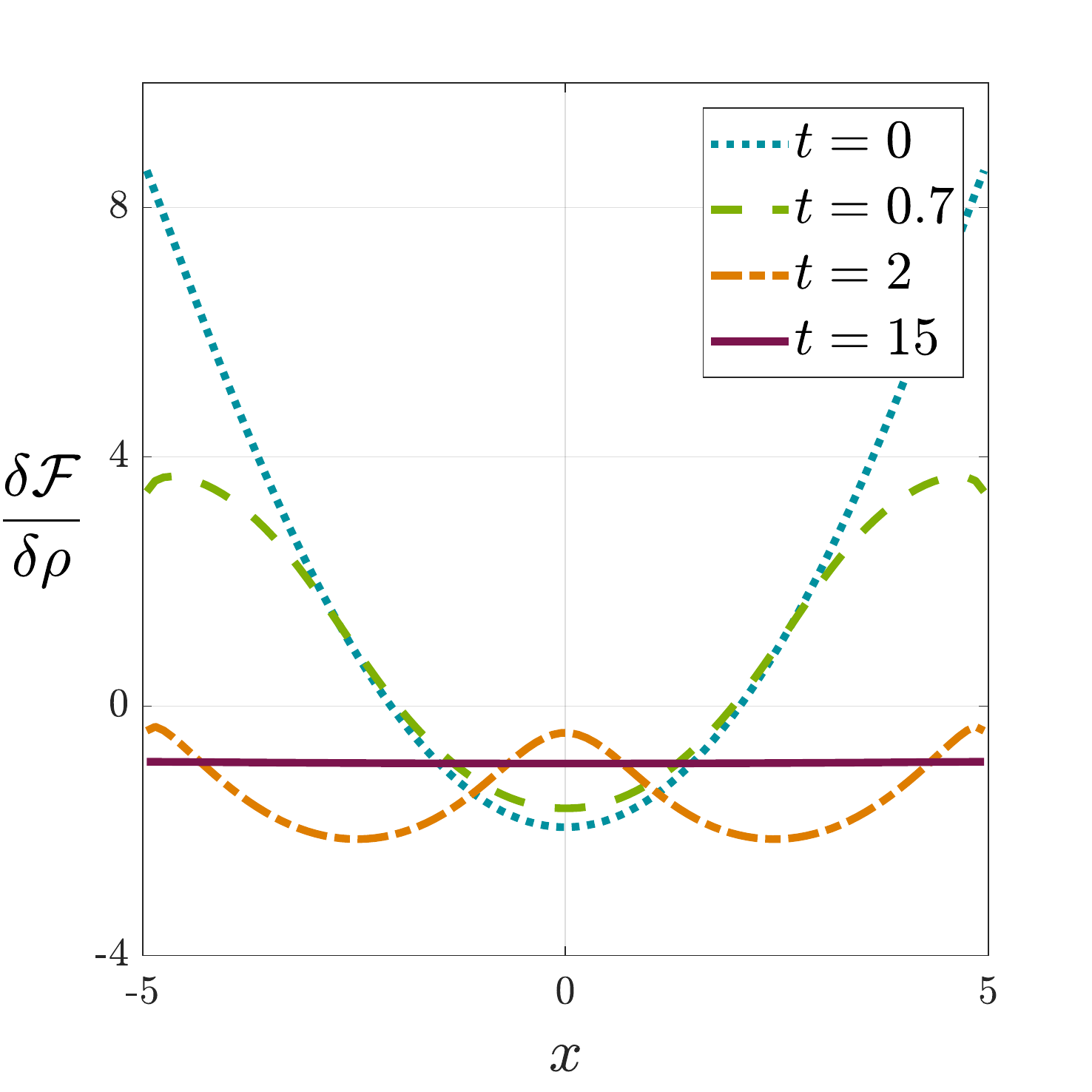}
}
\subfloat[Evolution of the total energy and free energy]{\protect\protect\includegraphics[scale=0.4]{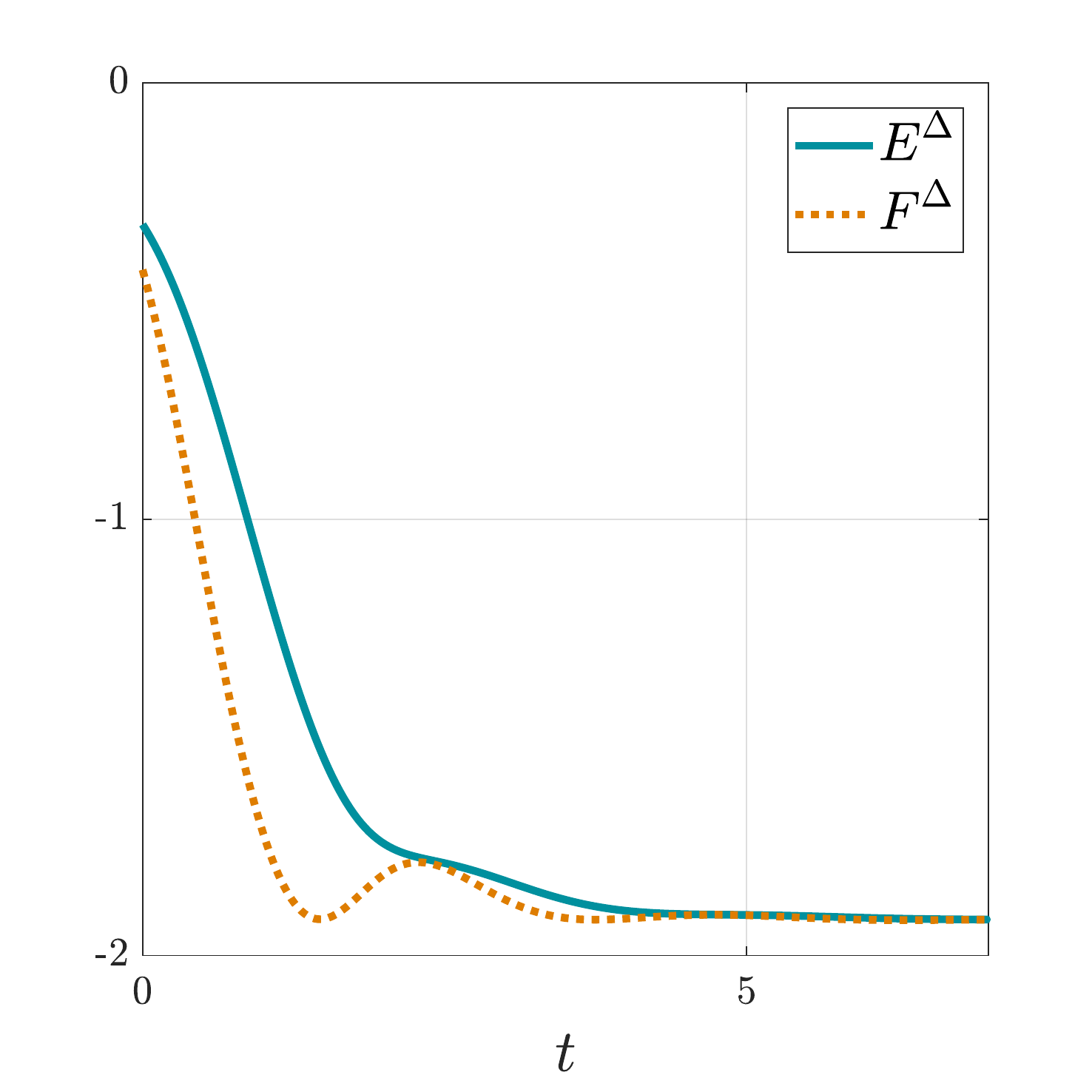}
}
\end{center}
\protect\protect\caption{\label{fig:idealpotCS} Temporal evolution of Example \ref{ex:idealpotCS}.}
\end{figure}

\begin{table}[!htbp]
\centering
\caption{Accuracy test for Example \ref{ex:idealpotCS} with the first and second-order schemes\blue{, at $t=0.3$}}
\label{table:idealpotCS}
\begin{tabular}{c c c c c }
\hline
\multirow{2}{*}{\begin{tabular}[c]{@{}c@{}}Number of \\ cells\end{tabular}} & \multicolumn{2}{c}{ First-order} & \multicolumn{2}{c}{Second-order}   \\ \cline{2-5}
                                                                                & $L^1$ error    & order    & $L^1$ error           & order          \\ \hline
50                                                                              &   6.3195E-03             &     -     &     7.3045E-04                 &     -             \\ \hline
100                                                                             &      3.2658E-03          &    0.95      &       1.9462E-04                &                  1.91 \\ \hline
200                                                                             &    1.6373E-03            &     1.00     &            4.8629E-05           &                 2.00 \\ \hline
400                                                                             &      8.7771E-04          &    1.01      &         1.2468E-05              &                1.97 \\ \hline
\end{tabular}

\end{table}

\end{examplecase}


\begin{examplecase}[Ideal-gas pressure and attractive kernel]\label{ex:idealker} In this case study the pressure
satisfies $P(\rho)=\rho$ and there is an interaction potential with a kernel
of the form $W(x)=\frac{x^2}{2}$. The steady state for a general total mass
$M_0$ is again equal to the steady states from examples \ref{ex:idealpot} and
\ref{ex:idealpotCS} with unit mass. The linear damping coefficient $\gamma$
has been reduced, $\gamma=0.01$, in order to compare the evolution with
respect to the previous examples.

The initial conditions for the order of accuracy test are the ones from example \ref{ex:idealpot}  in  \eqref{eq:icidealpot}. The order of accuracy test from this example is shown in table \ref{table:idealker}, and the evolution of the density, momentum, variation of the free energy with respect to the density, total energy and free energy are depicted in figure \ref{fig:idealker}. Due to the low value of $\gamma$ in the linear damping, there is a repeated exchange of free energy and kinetic energy during the temporal evolution, which can be noticed from the oscillations of the free energy plot in figure \ref{fig:idealker} (D). In the previous examples the linear damping term dissipates the momentum in a faster timescale and these exchanges only last for a few oscillations. One can also notice that the time to reach the steady state is higher than in the previous examples.
\begin{figure}[ht!]
\begin{center}
\subfloat[Evolution of the density]{\protect\protect\includegraphics[scale=0.4]{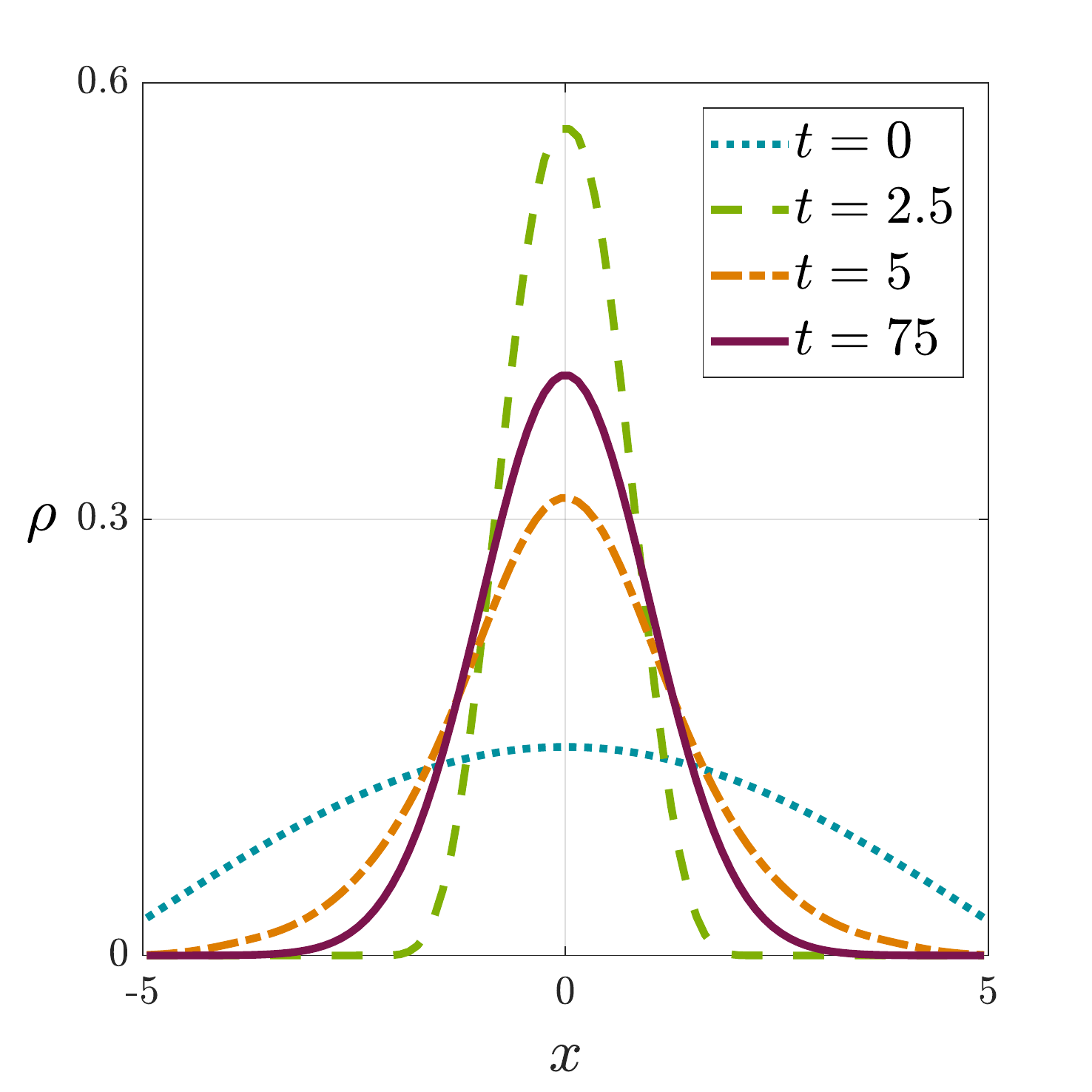}
}
\subfloat[Evolution of the momentum]{\protect\protect\includegraphics[scale=0.4]{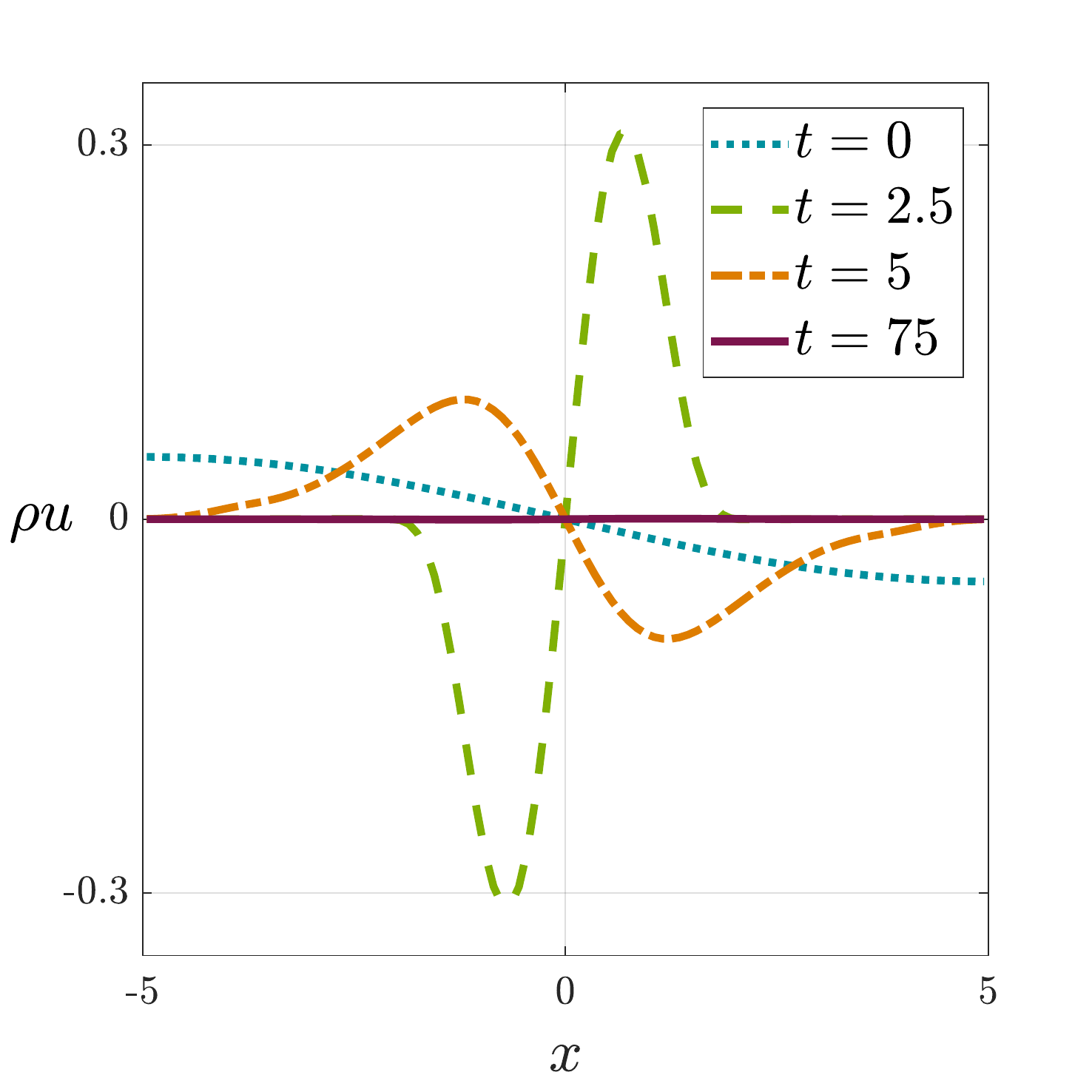}
}\\
\subfloat[Evolution of the variation of the free energy]{\protect\protect\includegraphics[scale=0.4]{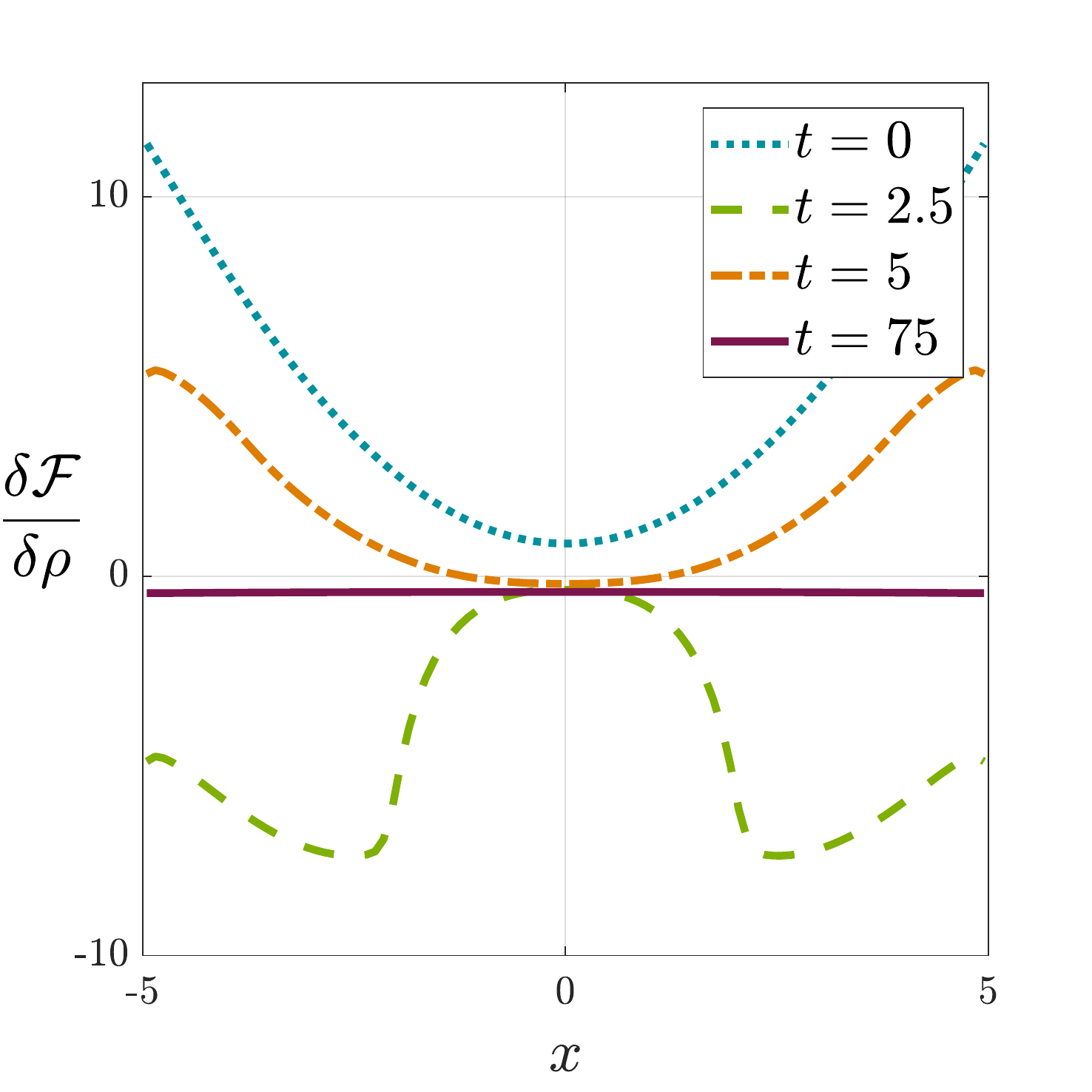}
}
\subfloat[Evolution of the total energy and free energy]{\protect\protect\includegraphics[scale=0.4]{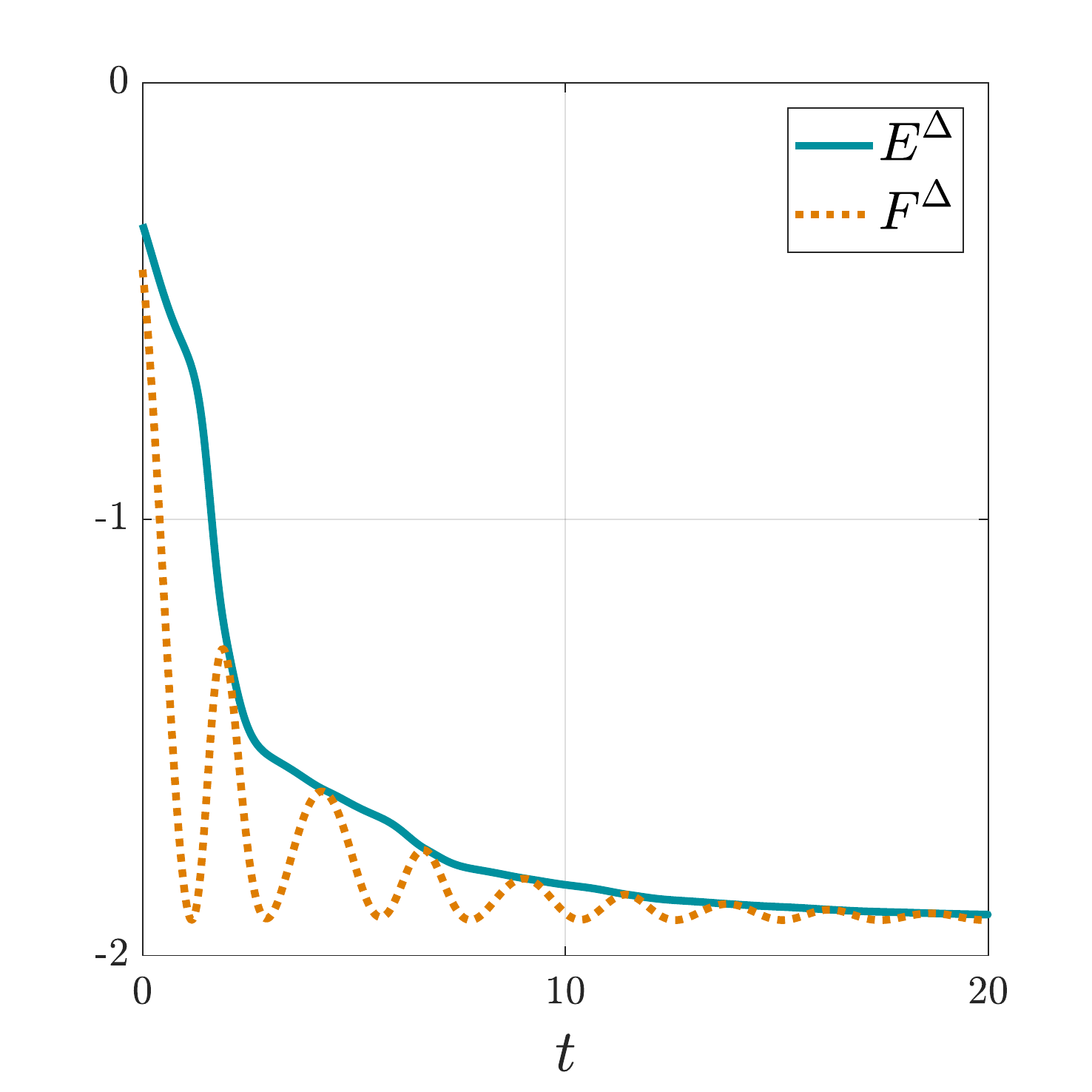}
}
\end{center}
\protect\protect\caption{\label{fig:idealker} Temporal evolution of Example \ref{ex:idealker}.}
\end{figure}

\begin{table}[!htbp]
\centering
\caption{Accuracy test for Example \ref{ex:idealker} with the first and second-order schemes\blue{, at $t=0.3$}}
\label{table:idealker}
\begin{tabular}{c c c c c }
\hline
\multirow{2}{*}{\begin{tabular}[c]{@{}c@{}}Number of \\ cells\end{tabular}} & \multicolumn{2}{c}{ First-order} & \multicolumn{2}{c}{Second-order}   \\ \cline{2-5}
                                                                                & $L^1$ error    & order    & $L^1$ error           & order          \\ \hline
50                                                                              &               6.6938E-03 &     -     &     7.6135E-04                 &     -             \\ \hline
100                                                                             &               3.4702E-03 &   0.95      &         2.0207E-04              &        1.91            \\ \hline
200                                                                             &               1.7410E-03 &   1.00       &         5.0306E-05              &       2.01           \\ \hline
400                                                                             &               8.6890E-04 &   1.00       &         1.2879E-05              &         1.97        \\ \hline
\end{tabular}

\end{table}

\end{examplecase}


\begin{examplecase}[Pressure proportional to square of density and attractive potential]\label{ex:squarepot} For this example the pressure satisfies $P(\rho)=\rho^2$ and there is an external potential of the form $V(x)=\frac{x^2}{2}$. Contrary to the previous examples \ref{ex:idealpot}, \ref{ex:idealpotCS} and \ref{ex:idealker}, the choice of $P(\rho)=\rho^2$ implies that regions of vacuum where $\rho=0$ appear in the evolution and steady solution of the system. As explained in the introduction of this section, the numerical flux employed for this case is a kinetic solver based on \cite{bouchut2004nonlinear}.

The steady state for this example with an initial mass of $M_0$  satisfies
\[
\rho_\infty (x) = \left\{ \begin{array}{ll}
\displaystyle -\frac{1}{4}\left(x+\sqrt[3]{3 M_0}\right)\left(x-\sqrt[3]{3  M_0}\right) \quad & \textrm{for} \quad x\in\left[-\sqrt[3]{3  M_0},\sqrt[3]{3  M_0}\right],\\[3mm]
 0 & \textrm{otherwise}.
  \end{array} \right.
\]

The initial conditions taken for the order of accuracy test are
\begin{equation*}
\rho(x,t=0)=M_0\frac{0.1+e^{-x^2}}{\int_\R \left(0.1+e^{-x^2} \right)dx},\quad \rho u(x,t=0)=-0.2 \sin\left(\frac{\pi x}{10}\right),\quad x\in [-5,5],
\end{equation*}
with $M_0$ being the mass of the system and equal to $1$.
The order of accuracy test from this example is shown in table \ref{table:squarepot}, and the evolution of the density, momentum, variation of the free energy with respect to the density, total energy and free energy are depicted in figure \ref{fig:squarepot}. The initial kinetic energy represents a large part of the initial total energy, and there is also an exchange between the kinetic energy and the free energy resulting in the oscillations for the plot of the discrete free energy.

As a remark, in this example the order of accuracy for the schemes with order higher than one is reduced to one \blue{both in the vacuum and interface regions}, as it is also pointed out in \cite{filbet2005approximation}. The orders showed in table  \ref{table:squarepot} \blue{are computed by considering only the cells in the support of the density that are away from the interface region, and the vacuum regions are not taken into consideration}.

\begin{figure}[ht!]
\begin{center}
\subfloat[Evolution of the density]{\protect\protect\includegraphics[scale=0.4]{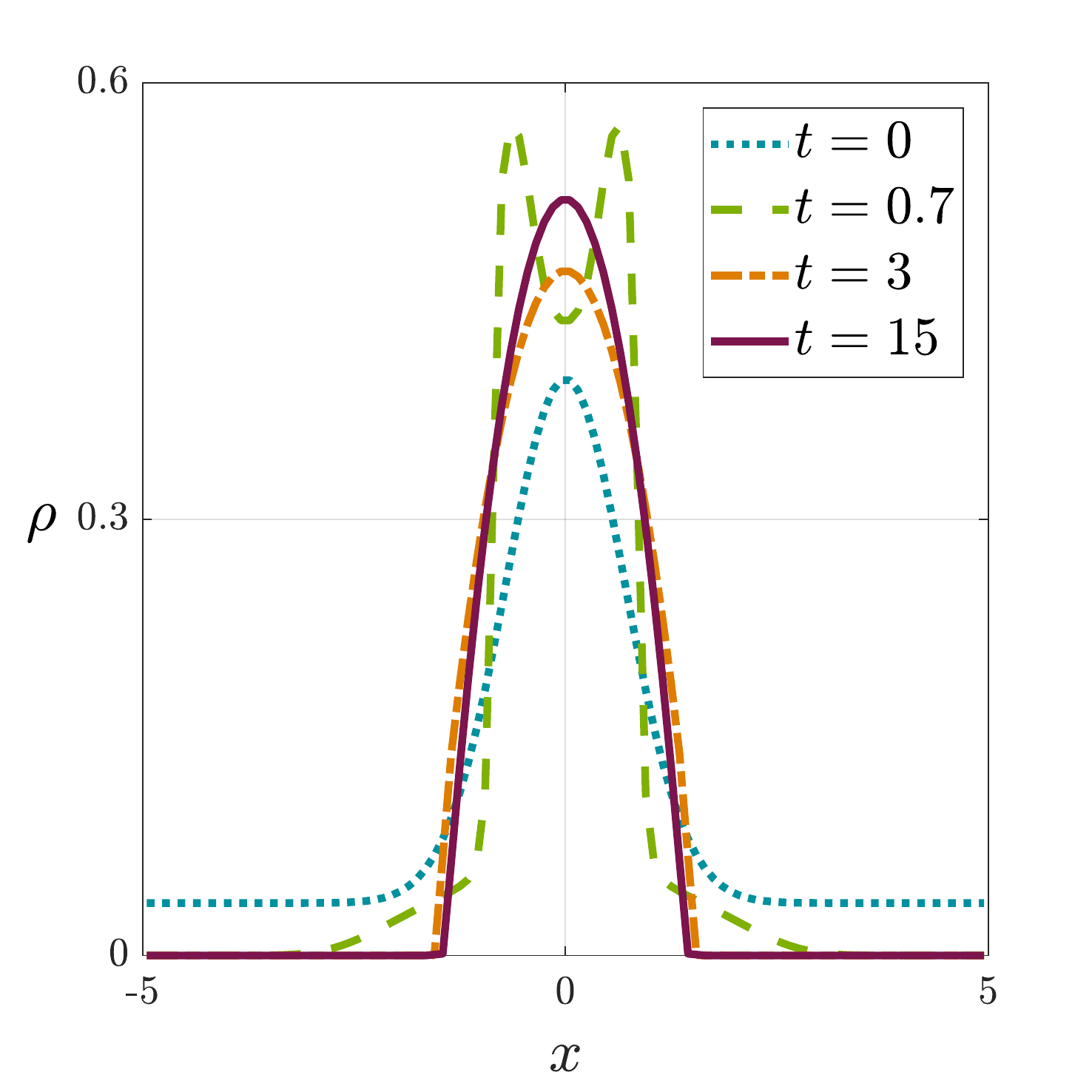}
}
\subfloat[Evolution of the momentum]{\protect\protect\includegraphics[scale=0.4]{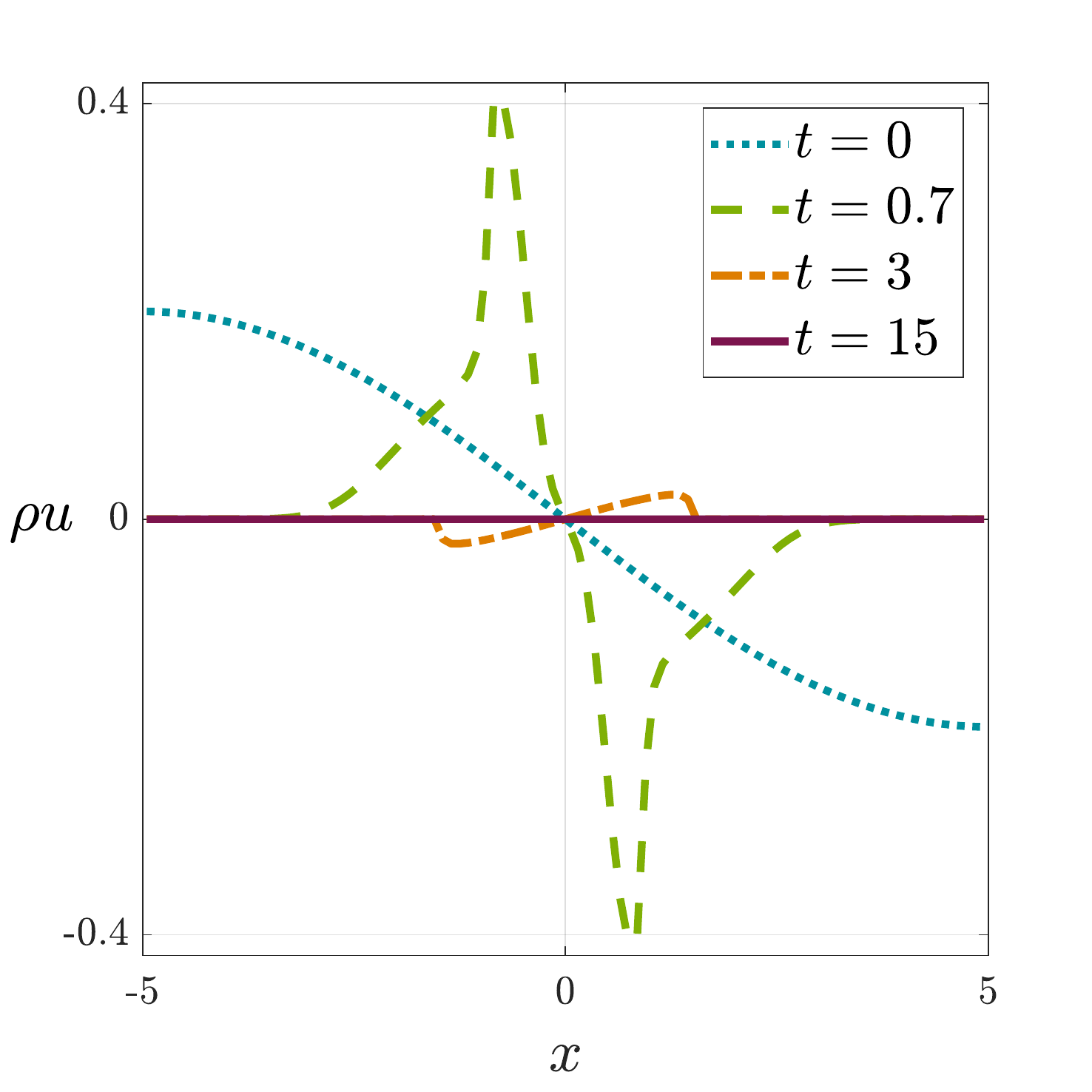}
}\\
\subfloat[Evolution of the variation of the free energy]{\protect\protect\includegraphics[scale=0.4]{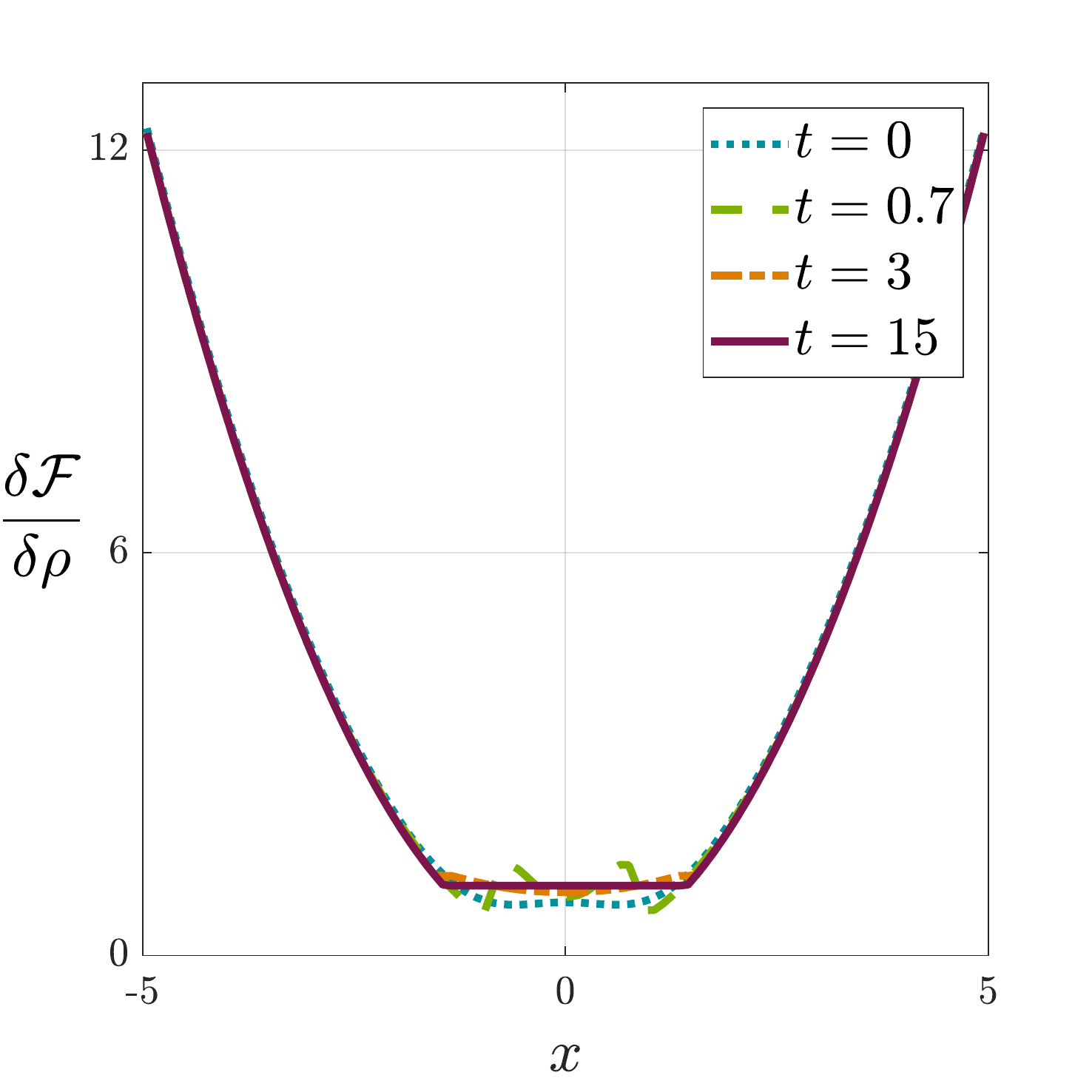}
}
\subfloat[Evolution of the total energy and free energy]{\protect\protect\includegraphics[scale=0.4]{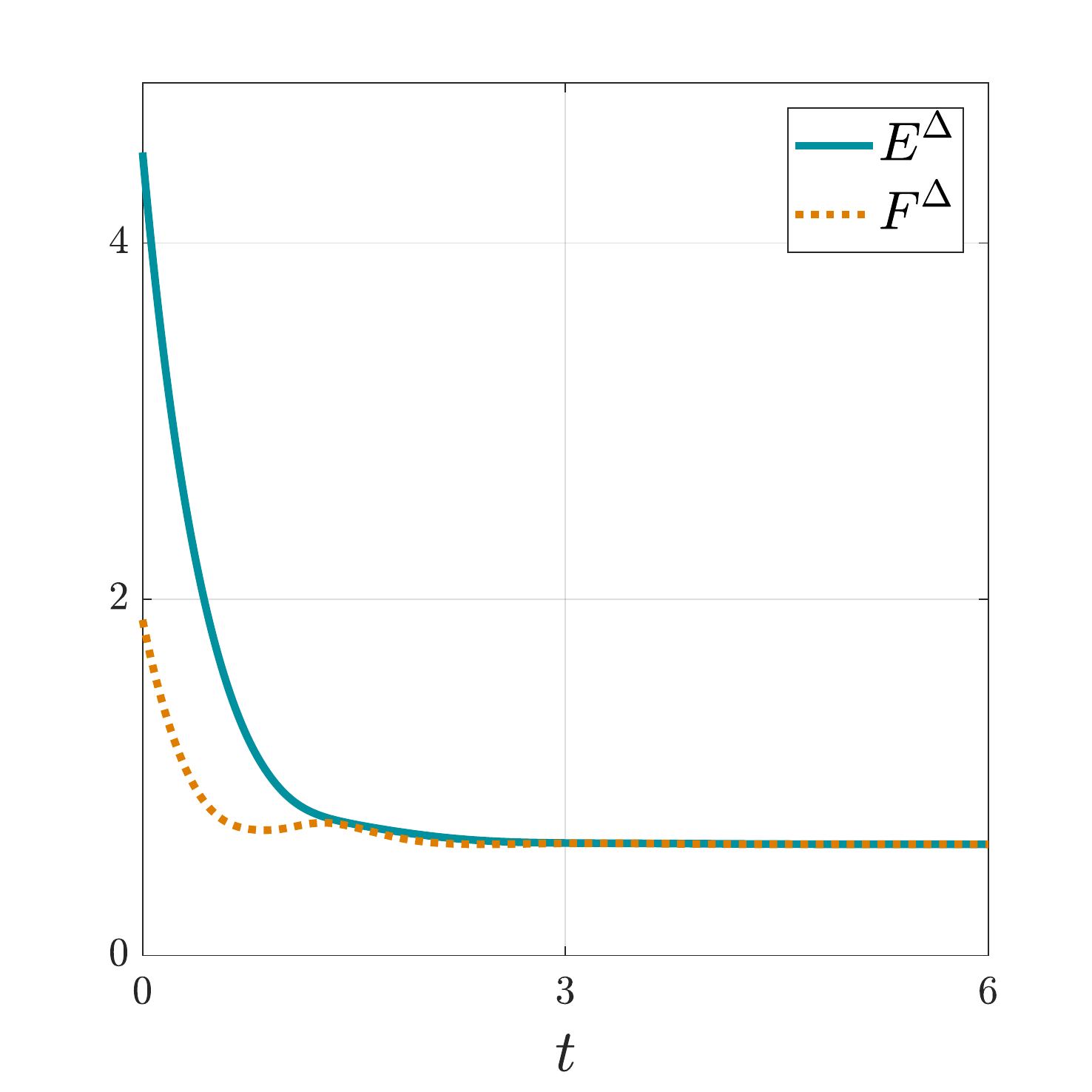}
}
\end{center}
\protect\protect\caption{\label{fig:squarepot} Temporal evolution of Example \ref{ex:squarepot}.}
\end{figure}

\begin{table}[!htbp]
\centering
\caption{Accuracy test for Example \ref{ex:squarepot} with the first and second-order schemes\blue{, at $t=0.3$}}
\label{table:squarepot}
\begin{tabular}{c c c c c }
\hline
\multirow{2}{*}{\begin{tabular}[c]{@{}c@{}}Number of \\ cells\end{tabular}} & \multicolumn{2}{c}{ First-order} & \multicolumn{2}{c}{Second-order}   \\ \cline{2-5}
                                                                                & $L^1$ error    & order    & $L^1$ error           & order          \\ \hline
50                                                                              &               6.8826E-03 &     -     &       1.0735E-03               &     -             \\ \hline
100                                                                             &               3.5106E-03 &    0.97      &         2.9188E-04              &       1.88            \\ \hline
200                                                                             &               1.7596E-03 &    1.00      &           7.6113E-05            &       1.94           \\ \hline
400                                                                             &               8.8184E-04 &   1.00       &           1.9103E-05            &      1.99           \\ \hline
\end{tabular}
\end{table}

\end{examplecase}


\begin{examplecase}[Moving steady state with ideal-gas pressure, attractive kernel and Cucker-Smale damping term]\label{ex:moving} The purpose of this example is to show that our scheme from section \ref{sec:numsch} preserves the order of accuracy for moving steady states of the form \eqref{eq:steadyvarener2}, where the velocity is not dissipated. As mentioned in the introduction, the generalization of well-balanced schemes to preserve moving steady states has proven to be quite complicated \cite{noelle2007high,xing2011advantage}, and it is not the aim of this work to construct such schemes.

For this example the pressure satisfies $P(\rho)=\rho$ and there is an interaction potential with a kernel of the form $W(x)=\frac{x^2}{2}$. The linear damping is eliminated and the Cucker-Smale damping term included. Under this configuration, there exists an explicit solution for system \eqref{eq:generalsys2} consisting in a travelling wave of the form

\begin{equation}\label{eq:travellingwave}
\rho(x,t)=M_0 \frac{e^{-(x-u t)^2/2}}{\int_\R e^{-x^2/2} dx},\quad u(x,t)=0.2,
\end{equation}
with $M_0$ equals to $1$ so that the total mass is unitary. As a result, the order of accuracy test can be accomplished by computing the error with respect to the exact reference solution, contrary to what was proposed in the previous examples. It should be remarked however that the velocity and the variation of the free energy with respect to the density profiles are not kept constant along the domain by our numerical scheme, since the well-balanced property for moving steady states is not satisfied.

The initial conditions for our simulation are \eqref{eq:travellingwave} at $t=0$, in a numerical domain with $ x\in [-8,9]$. The simulation is run until $t=3$. The table of errors for different number of cells is showed in table \ref{table:movingwater}, and a depiction of the evolution of the system is illustrated in figure \ref{fig:moving}. The velocity and the variation of the free energy plots are not included since they are not maintained constant with our scheme.
\begin{table}[!htbp]
\centering
\caption{Accuracy test for Example \ref{ex:moving} with the first and second-order schemes\blue{, at $t=3$}}
\label{table:movingwater}
\begin{tabular}{c c c c c }
\hline
\multirow{2}{*}{\begin{tabular}[c]{@{}c@{}}Number of \\ cells\end{tabular}} & \multicolumn{2}{c}{ First-order} & \multicolumn{2}{c}{Second-order}   \\ \cline{2-5}
                                                                                & $L^1$ error    & order    & $L^1$ error           & order          \\ \hline
50                                                                              &       9.84245E-03         &     -     &        2.78988E-03              &     -             \\ \hline
100                                                                             &       4.92029E-03         &   1.00      &         9.09342E-04              &       1.62            \\ \hline
200                                                                             &       2.44627E-03          &    1.01      &       2.55340E-04               &       1.83           \\ \hline
400                                                                             &       1.21228E-03         &   1.01       &            7.47905E-05           &      1.77           \\ \hline
\end{tabular}
\end{table}

\begin{figure}[ht!]
\begin{center}
\subfloat[Evolution of the density]{\protect\protect\includegraphics[scale=0.4]{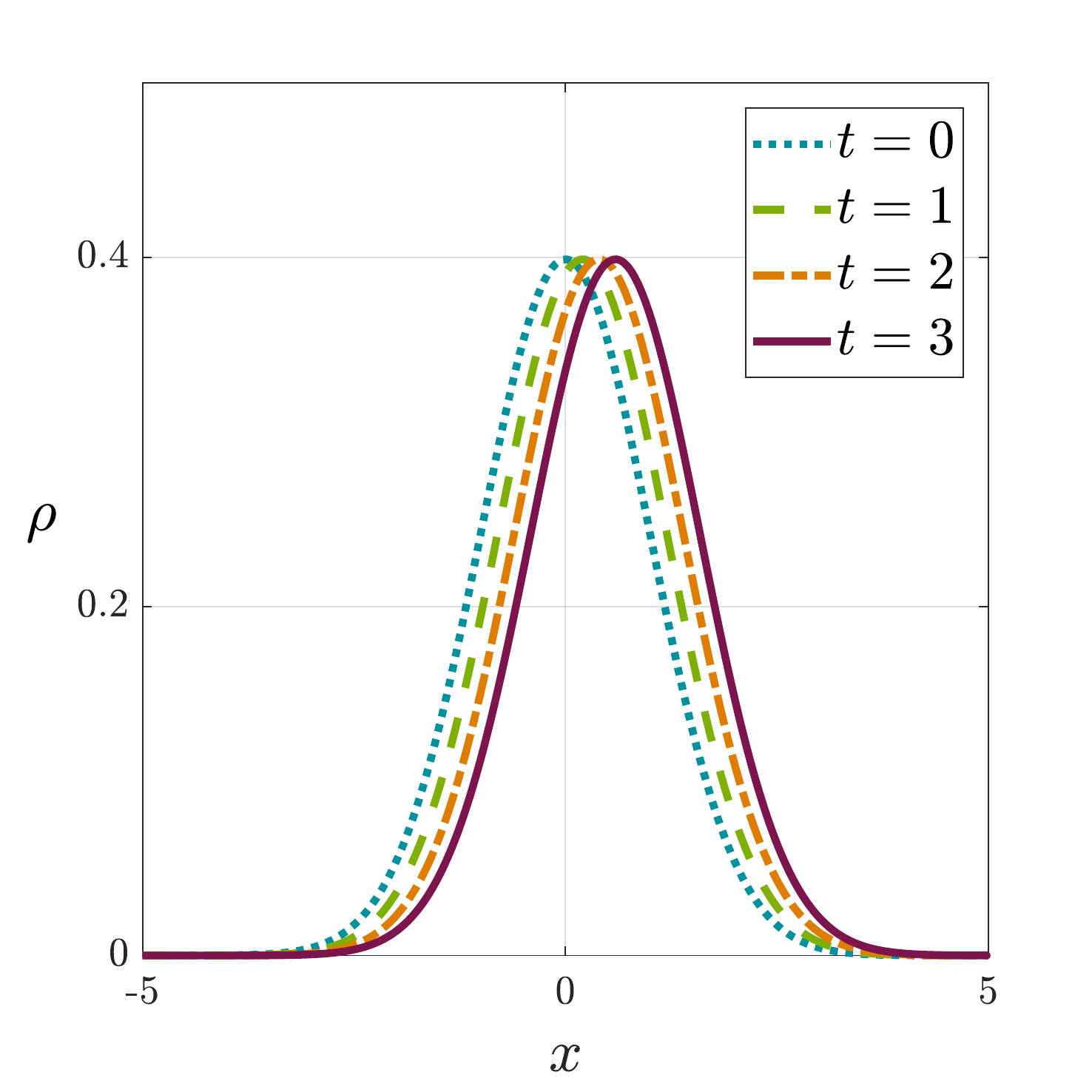}
}
\subfloat[Evolution of the total energy and free energy]{\protect\protect\includegraphics[scale=0.4]{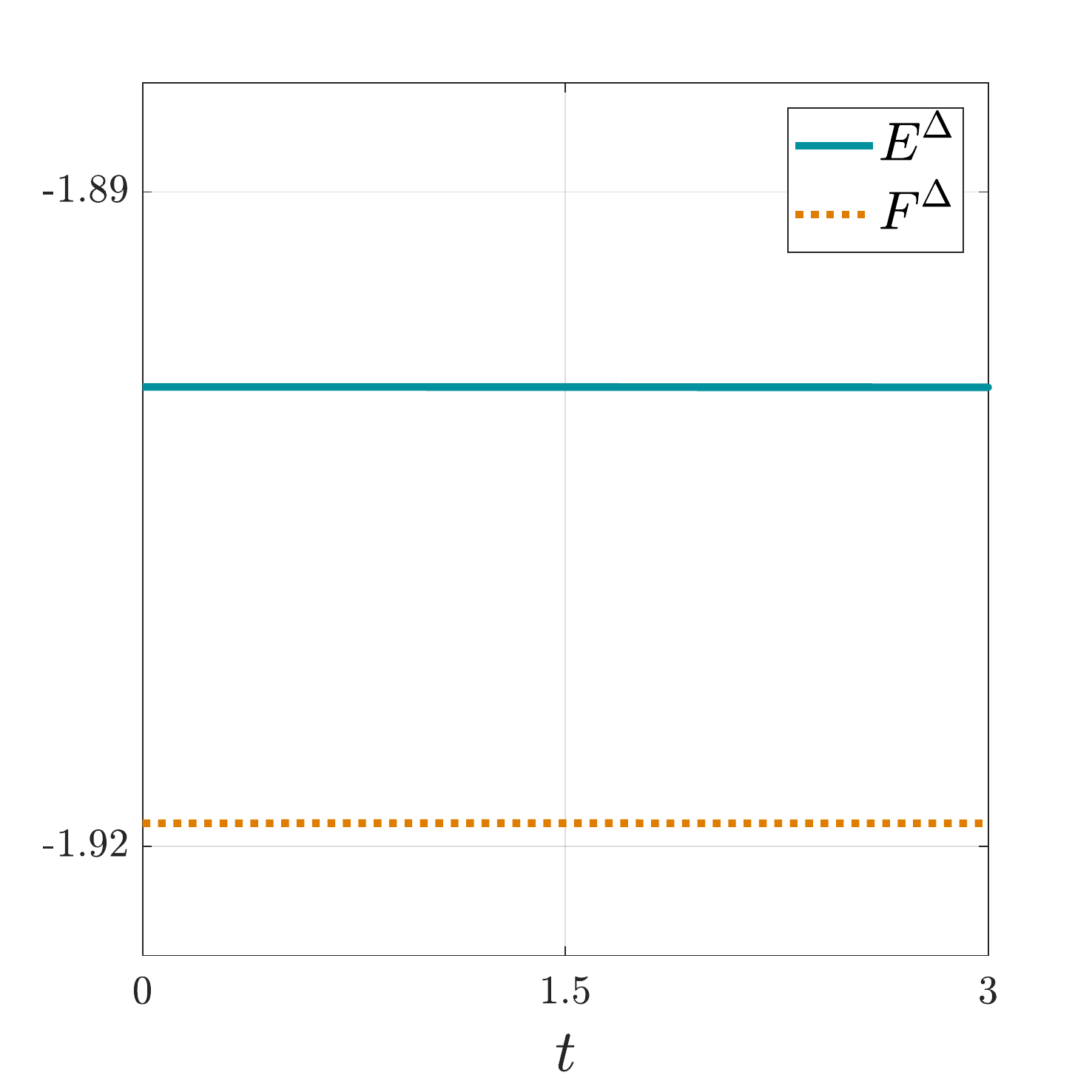}
}
\end{center}
\protect\protect\caption{\label{fig:moving} Temporal evolution of Example \ref{ex:moving}.}
\end{figure}
\end{examplecase}

%
%
\subsection{Numerical experiments}\label{subsec:numexp} This subsection applies the well-balanced scheme in section \ref{sec:numsch} to a variety of free energies from systems which have acquired an important consideration in the literature. Some of these systems have been mainly studied in their overdamped form, resulting when $\gamma \to \infty$, and as a result our well-balanced scheme can be useful in determining the role that inertia plays in those systems.

\begin{examplecase}[Pressure proportional to square of density and double-well potential]
\label{ex:squarepotdoublewell} In this example the pressure is taken as in example \ref{ex:squarepot}, with $P(\rho)=\rho^2$, thus leading to vacuum regions. The external potential are chosen to have a double-well shape of the form $V(x)=a\,x^4-b\,x^2$, with $a,\,b>0$. This system exhibits a variety of steady states depending on the symmetry of the initial condition, the initial mass and the shape of the external potential $V(x)$. The general expression for the steady states is
\begin{equation*}
\rho_{\infty}=\left(C(x)-V(x)\right)_{+}=\left(C(x)-a\,x^4+b\,x^2\right)_{+},
\end{equation*}
where $C(x)$ is a piecewise constant function, zero outside the support of the density. Notice that $C(x)$ can attain a different value in each connected component of the support of the density.

Three different initial data are simulated in order to compare the resulting long time asymptotics, i.e., we show that different steady states are achieved corresponding to different initial data. The initial conditions are
\begin{equation*}
\rho(x,t=0)=M_0\frac{0.1+e^{-(x-x_0)^2}}{\int_\R \left(0.1+e^{-(x-x_0)^2}\right)dx},\quad \rho u(x,t=0)=-0.2 \sin\left(\frac{\pi x}{10}\right),\quad x\in [-10,10],
\end{equation*}
with $M_0$ equal to $1$ so that the total mass is unitary. When $x_0=0$, the initial density is symmetric, and when $x_0\neq0$ the initial density is asymmetric.

\begin{enumerate}[label=\alph*.]
\item First case: The external potential satisfies
    $V(x)=\frac{x^4}{4}-\frac{3x^2}{2}$ and the initial density is
    symmetric with $x_0=0$. For this configuration the steady solution
    presents two disconnected bumps of density with the same mass in each
    of them, as it is shown in figure \ref{fig:doublewell} (A) and (B). The
    variation of the free energy with respect to the density presents the
    same constant value in the two disconnected supports of the density.
    The evolution is symmetric throughout.

\begin{figure}[]
\begin{center}
\subfloat[Density in first case]{\protect\protect\includegraphics[scale=0.4]{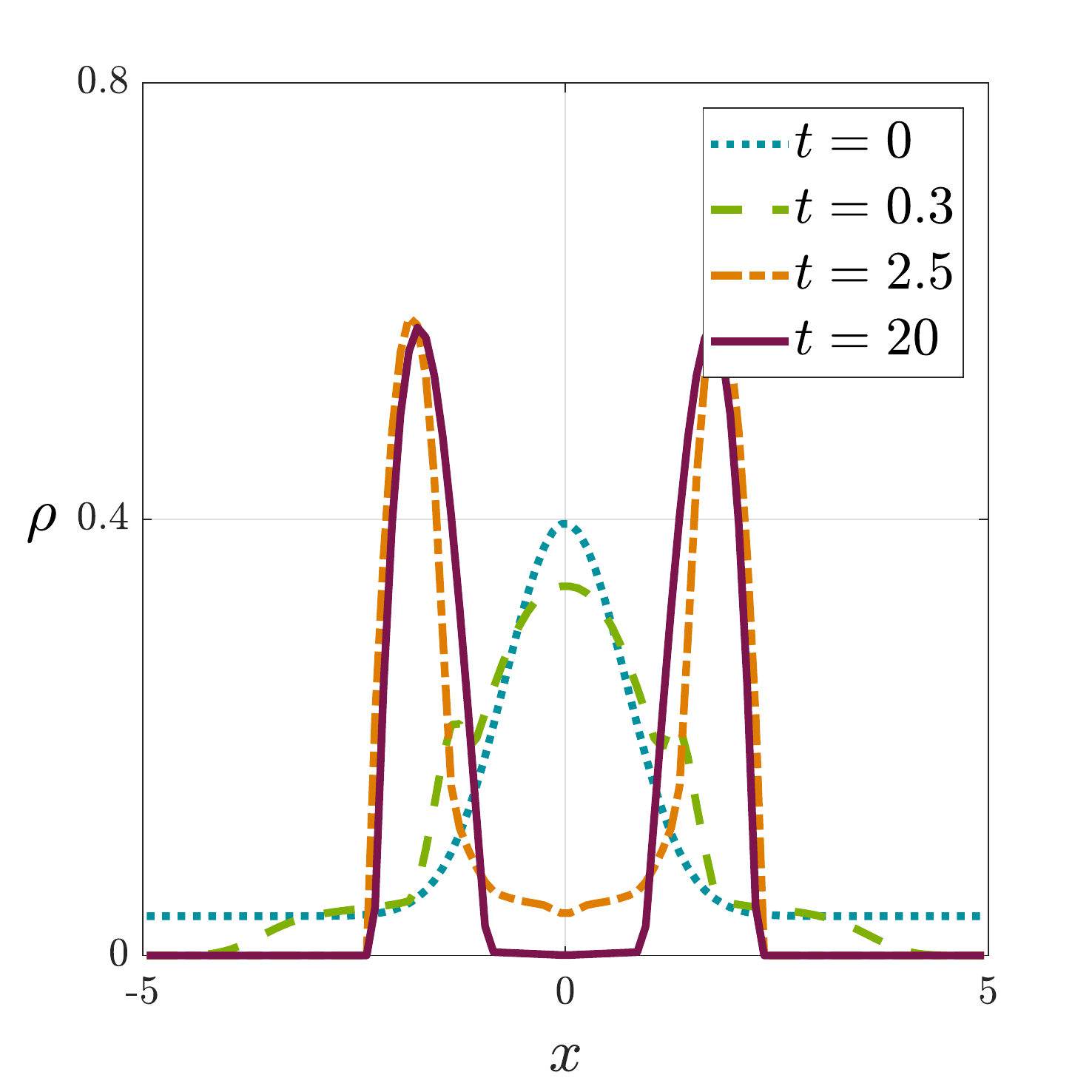}
}
\subfloat[Variation of the free energy in first case]{\protect\protect\includegraphics[scale=0.4]{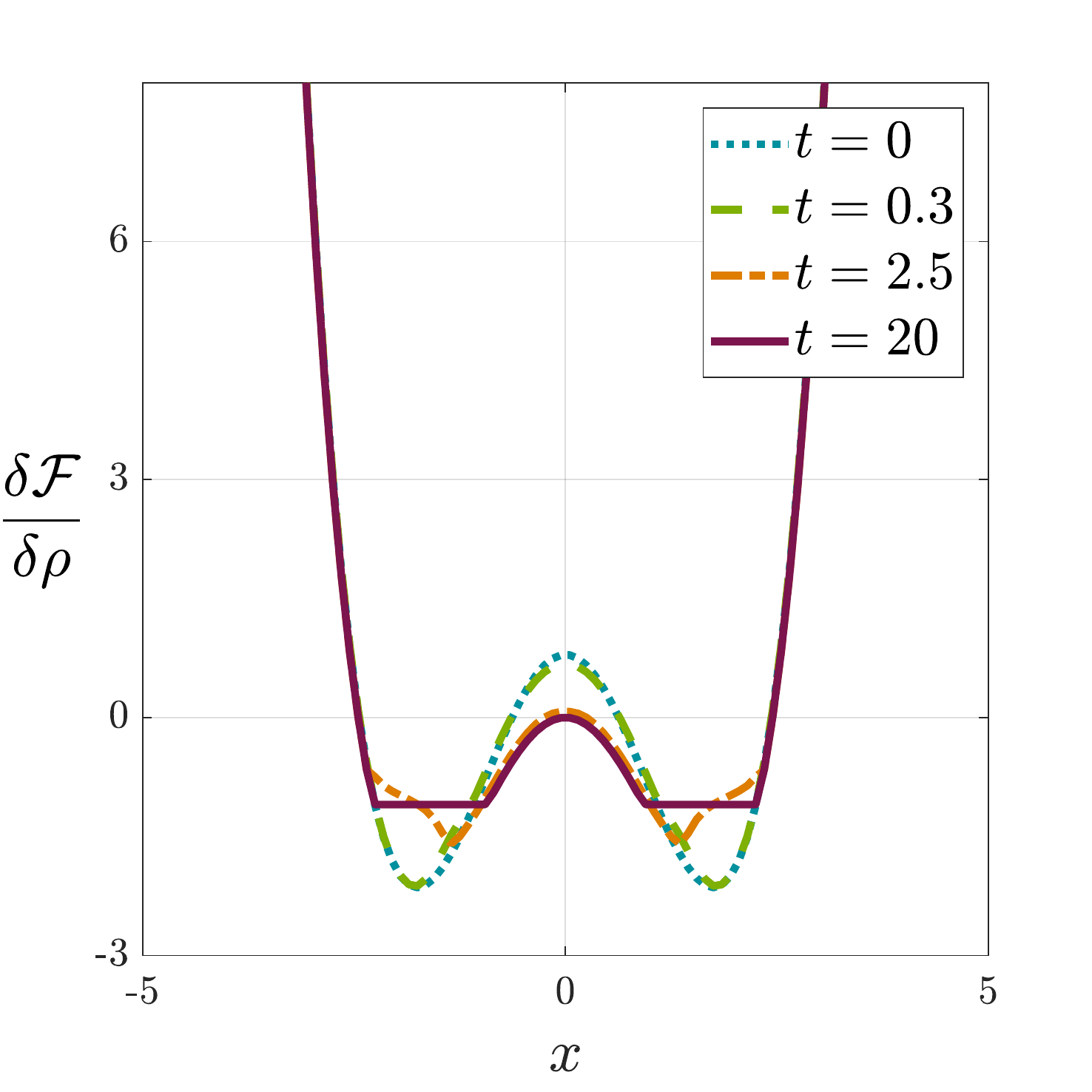}
}\\
\subfloat[Density in second case]{\protect\protect\includegraphics[scale=0.4]{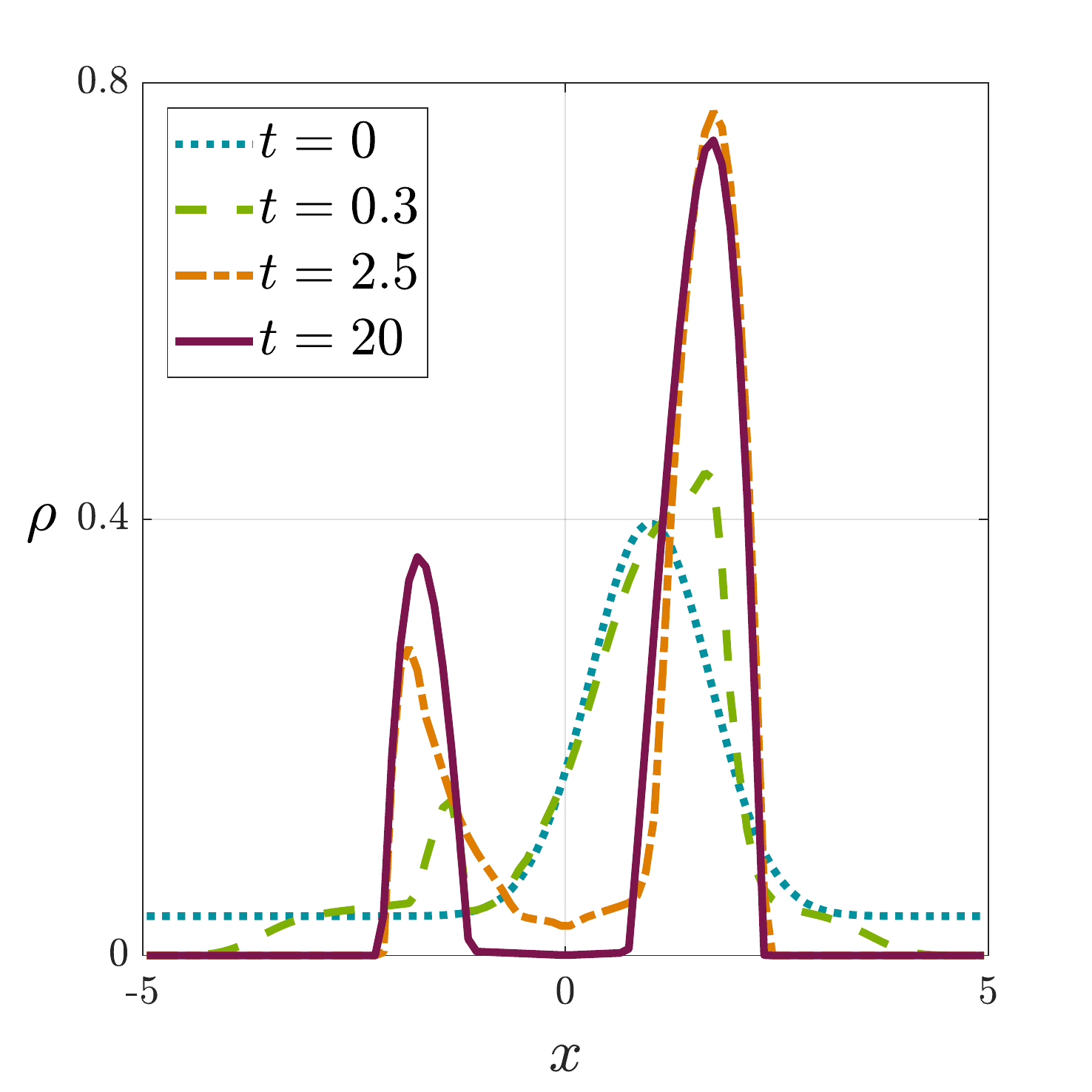}
}
\subfloat[Evolution of the variation of the free energy in second case]{\protect\protect\includegraphics[scale=0.4]{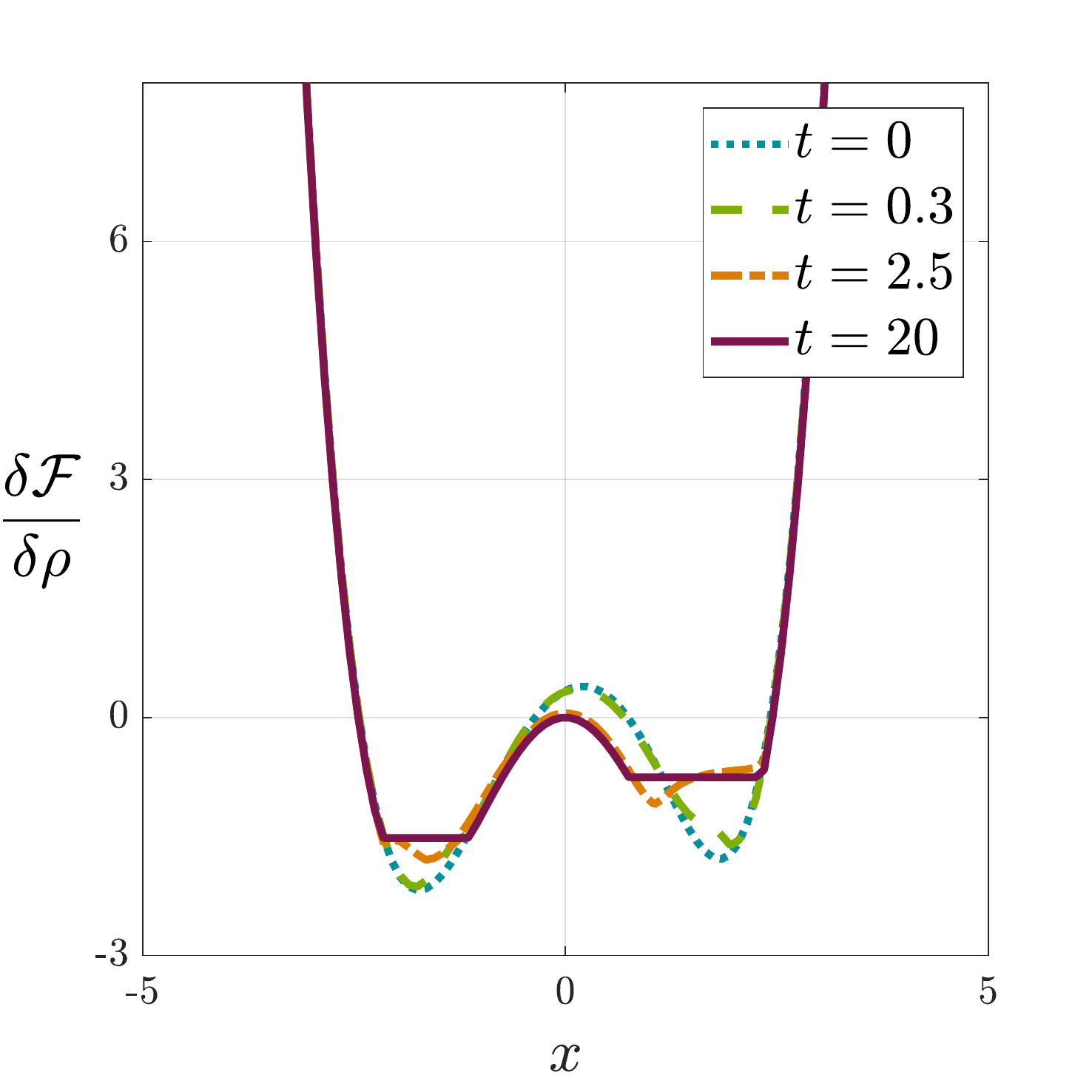}
}\\
\subfloat[Density in third case]{\protect\protect\includegraphics[scale=0.4]{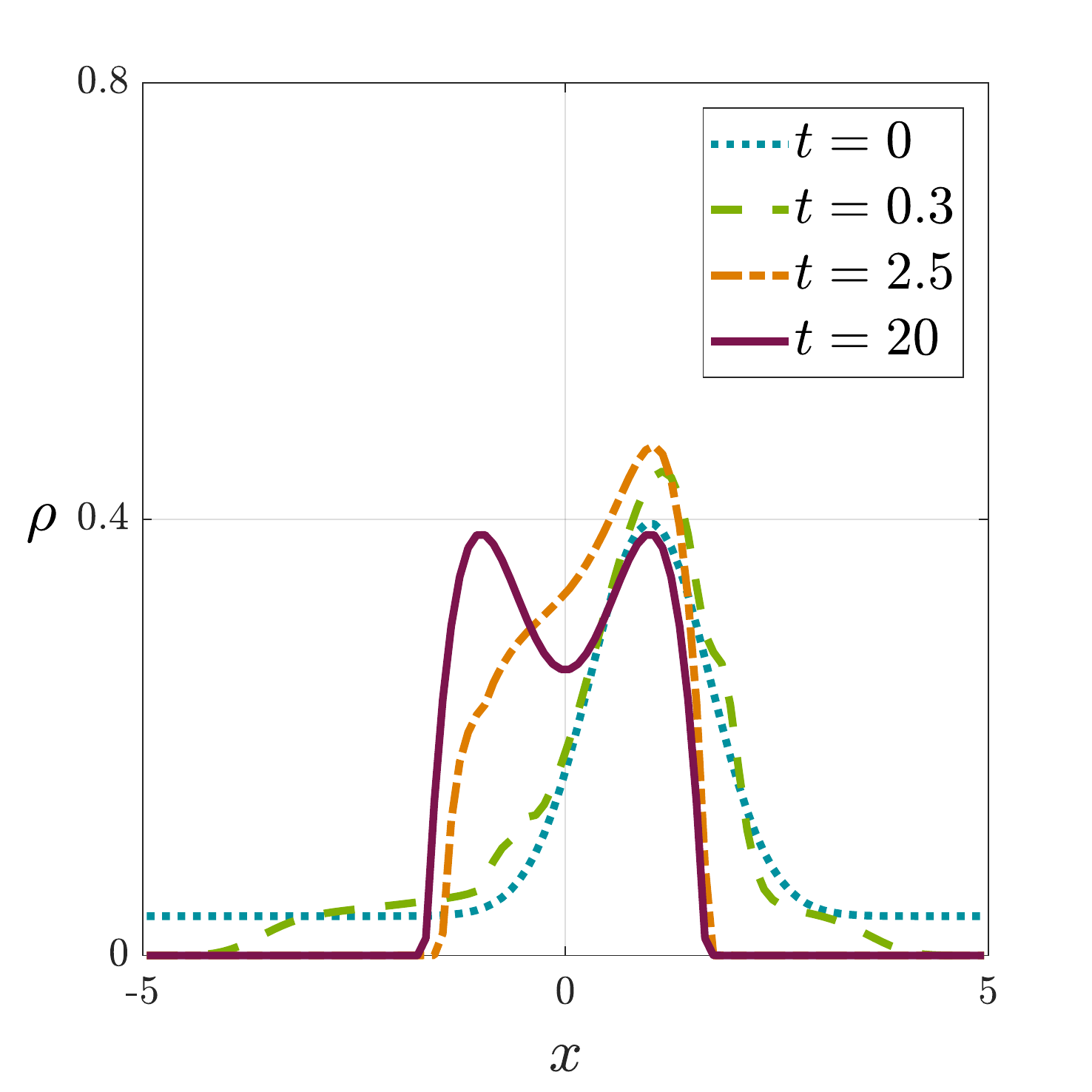}
}
\subfloat[Variation of the free energy in third case]{\protect\protect\includegraphics[scale=0.4]{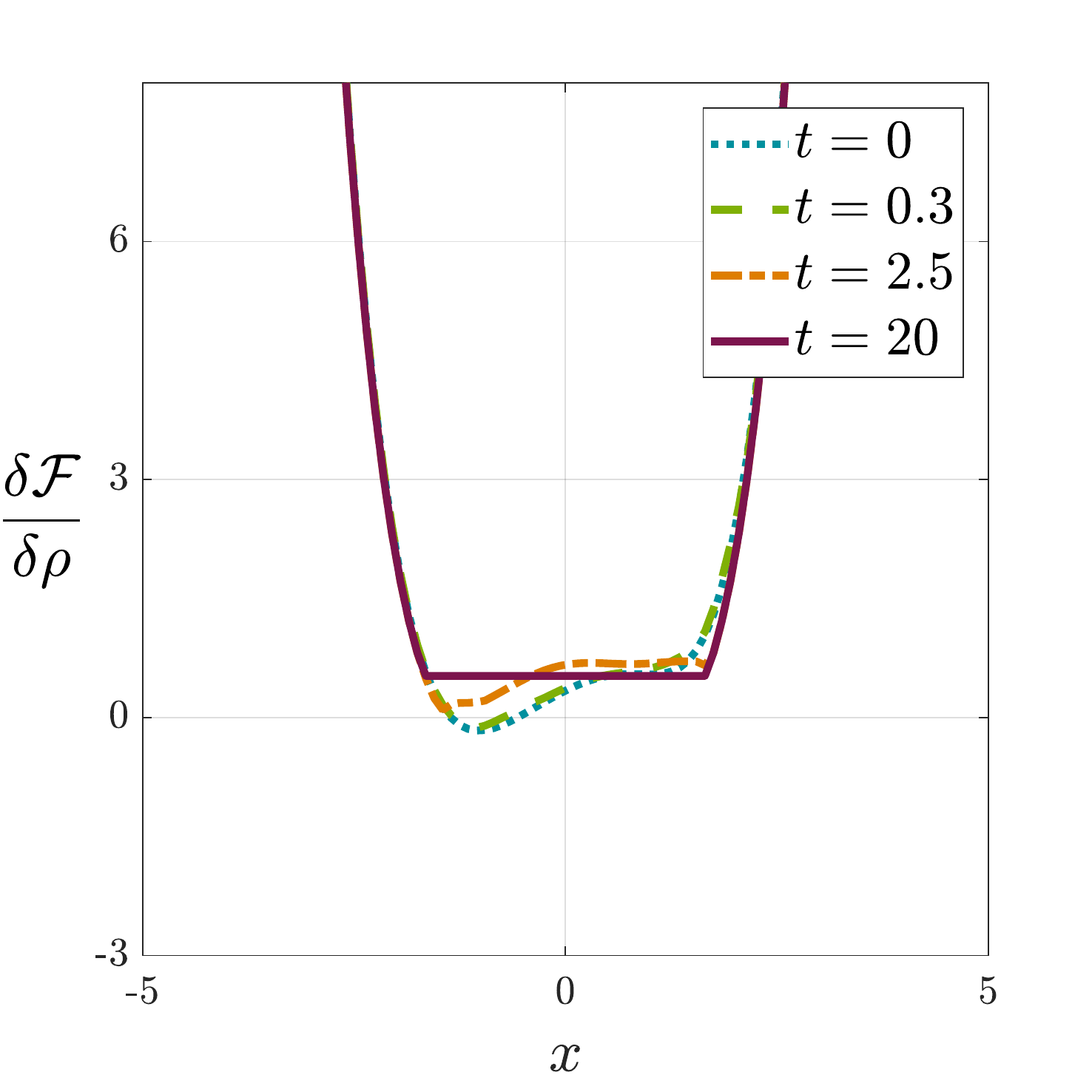}
}
\end{center}
\protect\protect\caption{\label{fig:doublewell} Temporal evolution of the first, second and third cases from example \ref{ex:squarepotdoublewell}.}
\end{figure}

\item Second case:  The external potential satisfies
    $V(x)=\frac{x^4}{4}-\frac{3x^2}{2}$ and the initial density is
    asymmetric with $x_0=1$. The final steady density is characterised
    again by the two disconnected supports but for this configuration the
    mass in each of them varies, as shown in figure \ref{fig:doublewell}
    (C) and (D). Similarly, the variation of the free energy with respect
    to the density presents different constant values in the two
    disconnected supports of the density.

\item Third case:  for this last configuration the external potential is varied and satisfies $V(x)=\frac{x^4}{4}-\frac{x^2}{2}$, and the initial density is asymmetric with $x_0=1$. For this case, even though the initial density is asymmetric, the final steady density is symmetric and compactly supported due to the shape of the potential, as it is shown in figure \ref{fig:doublewell} (E) and (F). The variation of the free energy with respect to the density presents constant value in all the support of the density.

\end{enumerate}

This behavior shows that this problem has a complicated bifurcation diagram and corresponding stability properties depending on the parameters, for instance the coefficient on the potential well controling the depth and support of the wells used above.

\end{examplecase}


\begin{examplecase}[Ideal pressure with noise parameter and its phase transition]\label{ex:bifurcation}
The model proposed for this example has a pressure satisfying $P(\rho)=\sigma \rho$, where $\sigma$ is a noise parameter, and external and interaction potentials chosen to be $V(x)=\frac{x^4}{4}-\frac{x^2}{2}$ and $W(x)=\frac{x^2}{2}$, respectively. The corresponding model in the overdamped limit has been previously studied in the context of collective behaviour \cite{barbaro2016phase}, mean field limits \cite{gomes2018mean}, and systemic risk \cite{garnier2013large}, see also \cite{tugaut2014phase} for the proof in one dimension.

We find that this hydrodynamic system exhibits a supercritical pitchfork
bifurcation in the center of mass  $\hat{x}$ of the steady state when varying
the noise parameter $\sigma$ as its overdamped limit counterpart discussed
above. For values of $\sigma$ higher than a certain threshold, all teady
states are symmetric and have the center of mass  $\hat{x}$ at $x=0$.
However, when $\sigma$ decreases below that threshold, the pitchfork
bifurcation takes place. On the one hand, if the center of mass of the
initial density is at $x=0$, the final center of mass in the steady state
remains at $x=0$. On the other hand, if the center of mass of the initial
density is at $x\neq0$, the center of mass of the steady state approaches
asymptotically to $x=1$ or $x=-1$ as $\sigma \to 0$, depending on the sign of
the initial center of mass. Finally, when $\sigma=0$, the steady state turns
into a Dirac delta at $x=0$, $x=1$ or $x=-1$, depending on the initial
density.
\begin{figure}[ht!]
\begin{center}
\subfloat[Bifurcation diagram]{\protect\protect\includegraphics[scale=0.4]{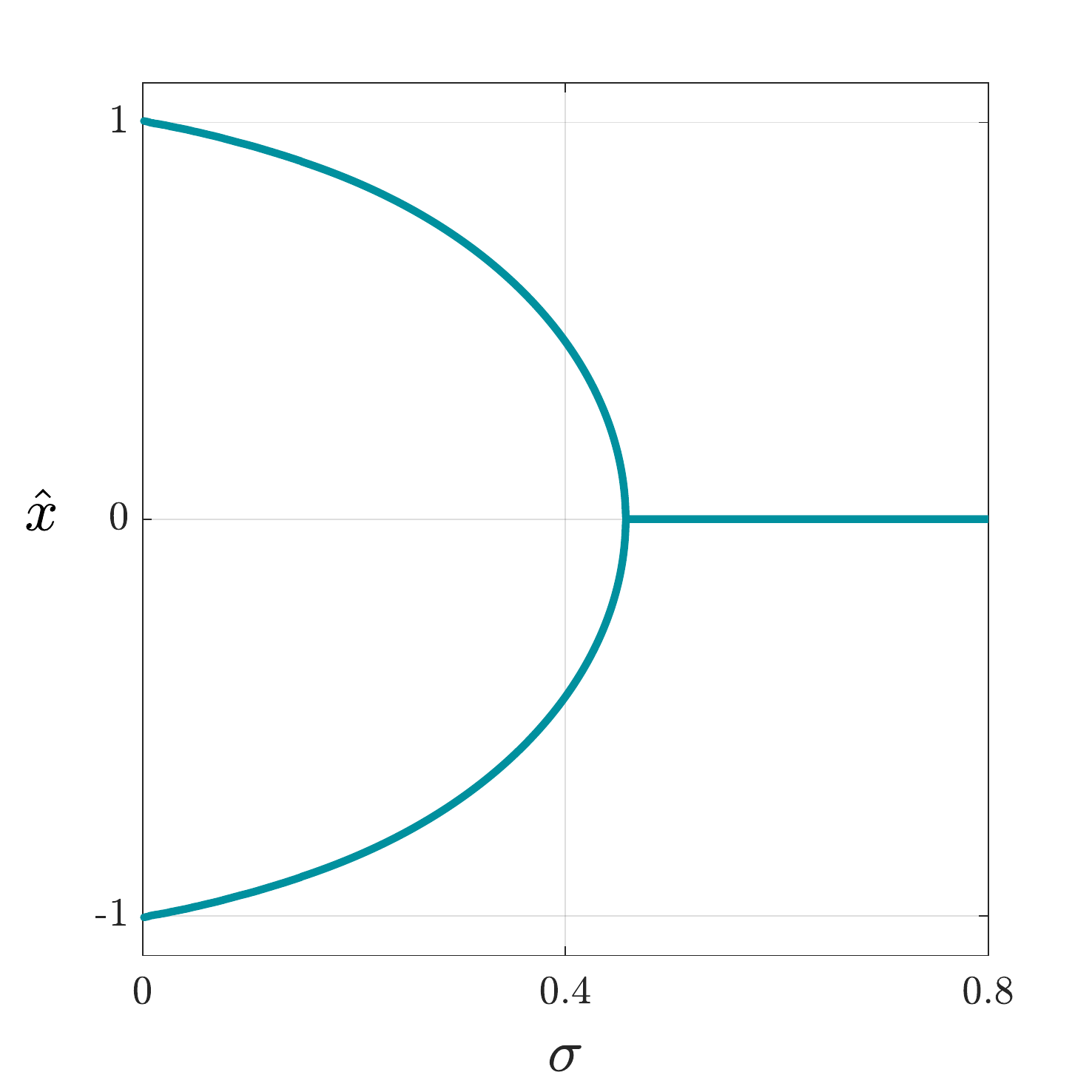}
}
\subfloat[Steady state profiles for different $\sigma$]{\protect\protect\includegraphics[scale=0.4]{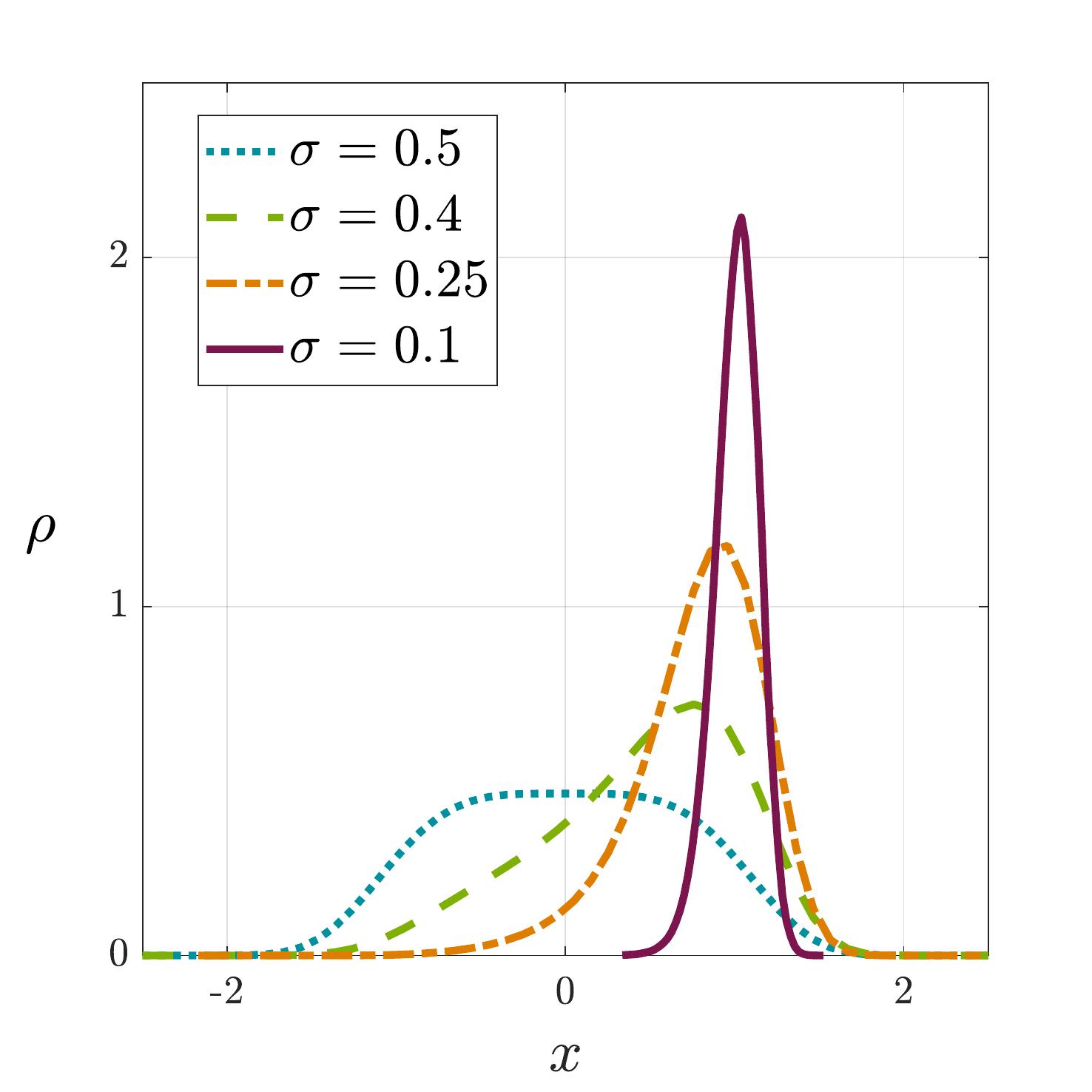}
\llap{\shortstack{%
        \includegraphics[scale=.14,trim={0 0 1.4cm 0cm},clip]{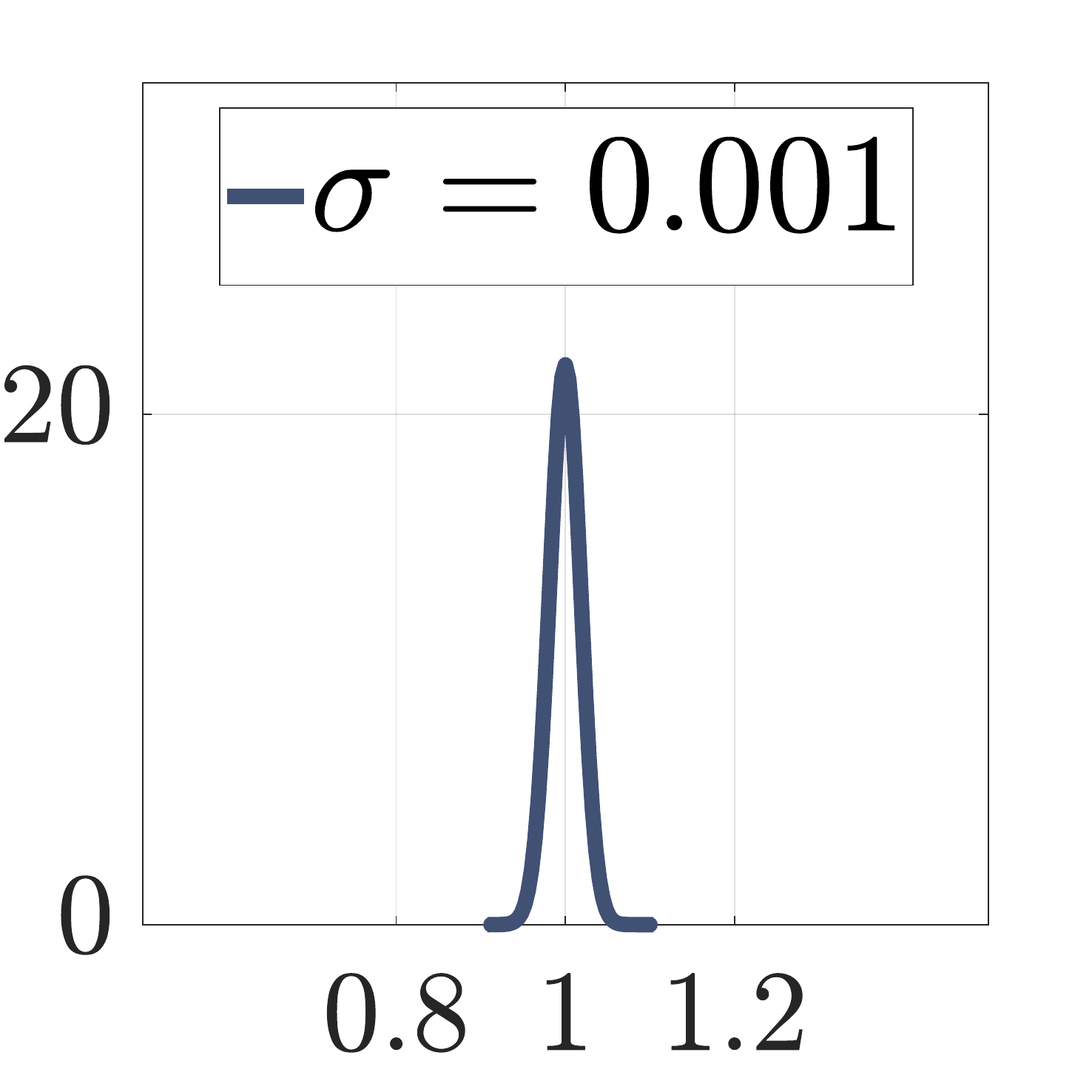}\\
        \rule{0ex}{0.65in}%
      }
  \rule{1.288in}{0ex}}
}
\end{center}
\protect\protect\caption{\label{fig:bifurcation} Bifurcation diagram (A) and steady states for different values of the noise parameter $\sigma$ (B) from Example \ref{ex:bifurcation}}
\end{figure}
The pitchfork bifurcation is supercritical since the branch of the bifurcation corresponding to $\hat{x}=0$ is unstable. This means that any deviation from an initial center of mass at $x=0$ leads to a steady center of mass located in one of the two branches of the parabola in the bifurcation state.

The numerical scheme outlined in section \ref{sec:numsch} captures this
bifurcation diagram for the evolution of the hydrodynamic system. The results
are shown in figure \ref{fig:bifurcation}. In it, (A) depicts the bifurcation
diagram of the final centre of mass when the noise parameter $\sigma$ is
varied, and for an initial center of mass at $x\neq0$. For a symmetric
initial density and antisymmetric velocity, the centre of mass numerically
remains at $x=0$ for an adequate stopping criterion, since property (vi) in
subsection \ref{subsec:firstorderprop} holds. However, any slight error in
the numerical computation unavoidably leads to a steady state deviating
towards any of the two stable branches, due to the strong unstable nature of
the branch with $x=0$. In (B) of figure \ref{fig:bifurcation} there is an
illustration of the steady states resulting from an initial center of mass
located at $x>0$, for different choices of the noise parameter $\sigma$. For
$\sigma=0.001$, which is the smallest value of $\sigma$ simulated, the
density profile approaches the theoretical Dirac delta expected at $x=1$ when
$\sigma \to 0$. When $\sigma=0$ the hyperbolicity of the system in
\eqref{eq:generalsys2} is lost since the pressure term vanishes, and as a
result the numerical approach in section \ref{sec:numsch} cannot be applied.

The numerical strategy followed to recover the bifurcation diagram is based
on the so-called differential continuation. It simply means that, as $\sigma
\to 0$, the subsequent simulations with new and lower values of $\sigma$ have
as initial conditions the previous steady state from the last simulation.
This allows to complete the bifurcation diagram, since otherwise the
simulations with really small $\sigma$ take long time to converge for general
initial conditions. In addition, to maintain sufficient resolution for the
steady states close to the Dirac delta, the mesh is adapted for each
simulation. This is accomplished by firstly interpolating the previous steady
state with a piecewise cubic hermite polynomial, which preserves the shape
and avoids oscillations, and secondly by creating a new and narrower mesh
where the interpolating polynomial is employed to construct the new initial
condition for the differential continuation.

\end{examplecase}


\begin{examplecase}[Hydrodynamic generalization of the Keller-Segel system - Generalized Euler-Poisson systems]\label{ex:kellersegel} The original Keller-Segel model has been widely employed in chemotaxis, which is usually defined as the directed movement of cells and organisms in response to chemical gradients \cite{keller1970initiation}. These systems also find their applications in astrophysics and gravitation \cite{sire2002thermodynamics,deng2002solutions}. It is a system of two coupled drift-diffusion differential equations for the density $\rho$ and the chemoattractant concentration $S$,
\begin{equation*}
\begin{cases}
{\displaystyle\partial_{t}\rho= \nabla \cdot \left(\nabla P(\rho)-\chi\rho\nabla S\right),
}\\[2mm]
{\displaystyle \partial_{t}S=D_s\Delta S-\theta S+\beta \rho.}
\end{cases}
\end{equation*}
In this system $P(\rho)$ is the pressure, and the biological/physical meaning of the constants $\chi$, $D_s$, $\alpha$ and $\beta$ can be reviewed in the literature \cite{horstmann20031970,bellomo2015toward,hoffmann2017keller}. For this example they are simplified as usual so that $\chi=D_s=\beta=1$ and $\theta=0$. A further assumption usually taken in the literature is that $\pa_{t} \rho$ is very big in comparison to $\pa_{t} S$ \cite{hoffmann2017keller}, leading to a simplification of the equation for the chemoattractant concentration $S$, which becomes the Poisson equation $-\Delta S=\rho$. Hydrodynamic extensions of the model, which include inertial effects, have also been proven to be essential for certain applications\cite{chavanis2006jeans,chavanis2007kinetic,gamba2003percolation}, leading to a hyperbolic system of equations with linear damping which in one dimension reads as
\begin{equation*}
\begin{cases}
\partial_{t}\rho+\pa_{x}\left(\rho u\right)=0,
\\[2mm]
{\displaystyle \pa_{t}(\rho u)\!+\pa_{x}(\rho u^2)\!=-\pa_{x} P(\rho)+\pa_{x}S -\gamma\rho u
,}\\[2mm]
{-\pa_{xx} S=\rho
.}
\end{cases}
\end{equation*}
By using the fundamental solution of the Laplacian in one dimension, this
equation becomes $2S =|x| \star \rho$. This term, after neglecting the
constant, can be plugged in the momentum equation so that the last equation
for $S$ can be removed. As a result, the hydrodynamic Keller-Segel model is
reduced to the system of equations \eqref{eq:generalsys} considered in this
work, with $W(x)=|x|/2$, $V(x)=0$ and $\psi\equiv 0$. As a final
generalization  \cite{carrillo2015finite}, the original interaction potential
$W(x)=|x|/2$ can be extended to be a homogeneous kernel
$W(x)=|x|^\alpha/\alpha$, where $\alpha>-1$. By convention, $W(x)=\ln|x|$ for
$\alpha=0$. Further generalizations are Morse-like potentials as in
\cite{carrillo2014explicit,carrillo2015finite} where
$W(x)=1-\exp(-|x|^\alpha/\alpha)$ with $\alpha> 0$.

The solution of this system can present a rich variety of behaviours due to
the competition between the attraction from the local kernel $W(x)$ and the
repulsion caused by the diffusion of the pressure $P(\rho)$, as reviewed in
\cite{calvez2017equilibria,calvez2017geometry}. By appropriately tuning the
parameters $\alpha$ in the kernel $W(x)$ and $m$ in the pressure $P(\rho)$,
one can find compactly supported steady states, self-similar behavior, or
finite-time blow up. Three different regimes have been studied in the
overdamped generalized Keller-Segel model \cite{carrillo2015finite}:
diffusion dominated regime ($m>1-\alpha$), balanced regime ($m=1-\alpha$)
where a critical mass separates self-similar and blow-up behaviour, and
aggregation-dominated regime ($m<1-\alpha$). These three regimes have not
been so far analytically studied for the hydrodynamic system except for few
particular cases \cite{carrillo2016critical,carrillo2016pressureless}, and
the presence of inertia indicates that the initial momentum profile  plays a
role together with the mass of the system to separate diffusive from blow-up
behaviour.

The well-balanced scheme provided in section \ref{sec:numsch} is a useful
tool to effectively reach the varied steady states resulting from different
values of $\alpha$ and $m$. The objective of this example is to provide some
numerical experiments to show the richness of possible behaviors. This scheme
can be eventually employed to numerically validate the theoretical studies
concerning the existence of the different regimes for the hydrodynamic system
for instance, or how the choice of the initial momentum or the total mass can
lead to diffusive or blow-up behaviour. This will be explored further
elsewhere.

\begin{figure}[ht!]
\begin{center}
\subfloat[Evolution of the density]{\protect\protect\includegraphics[scale=0.4]{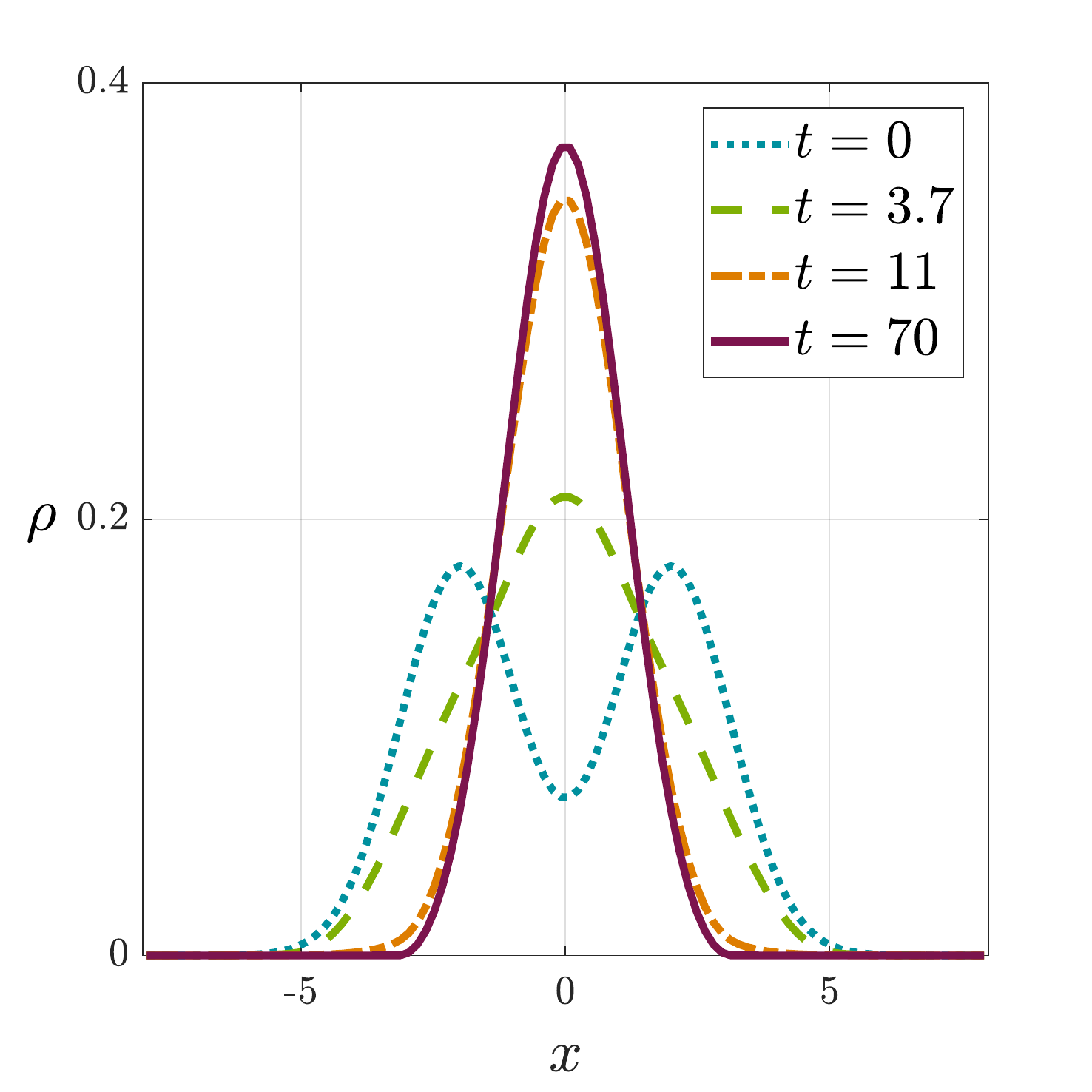}
}
\subfloat[Evolution of the momentum]{\protect\protect\includegraphics[scale=0.4]{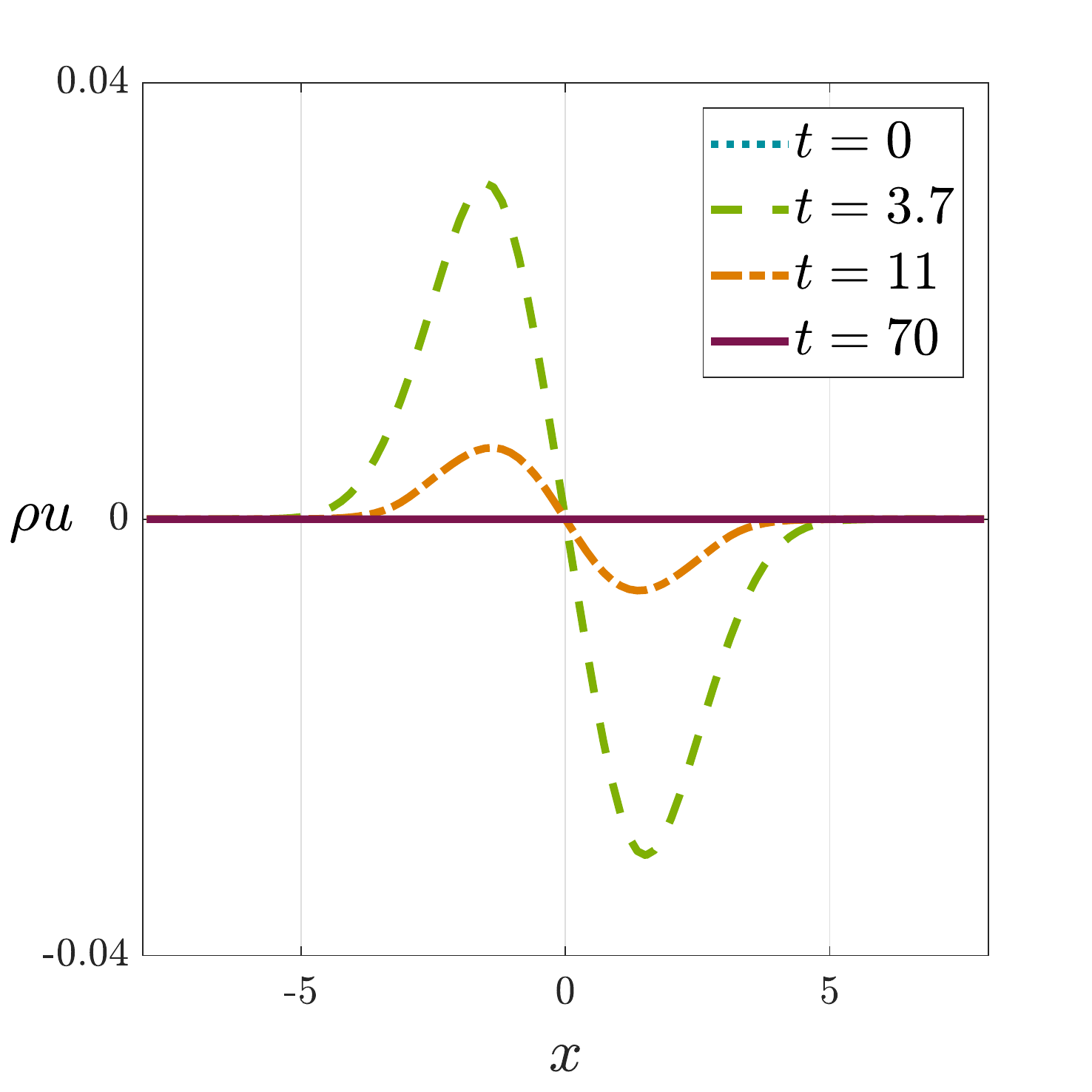}
}\\
\subfloat[Evolution of the variation of the free energy]{\protect\protect\includegraphics[scale=0.4]{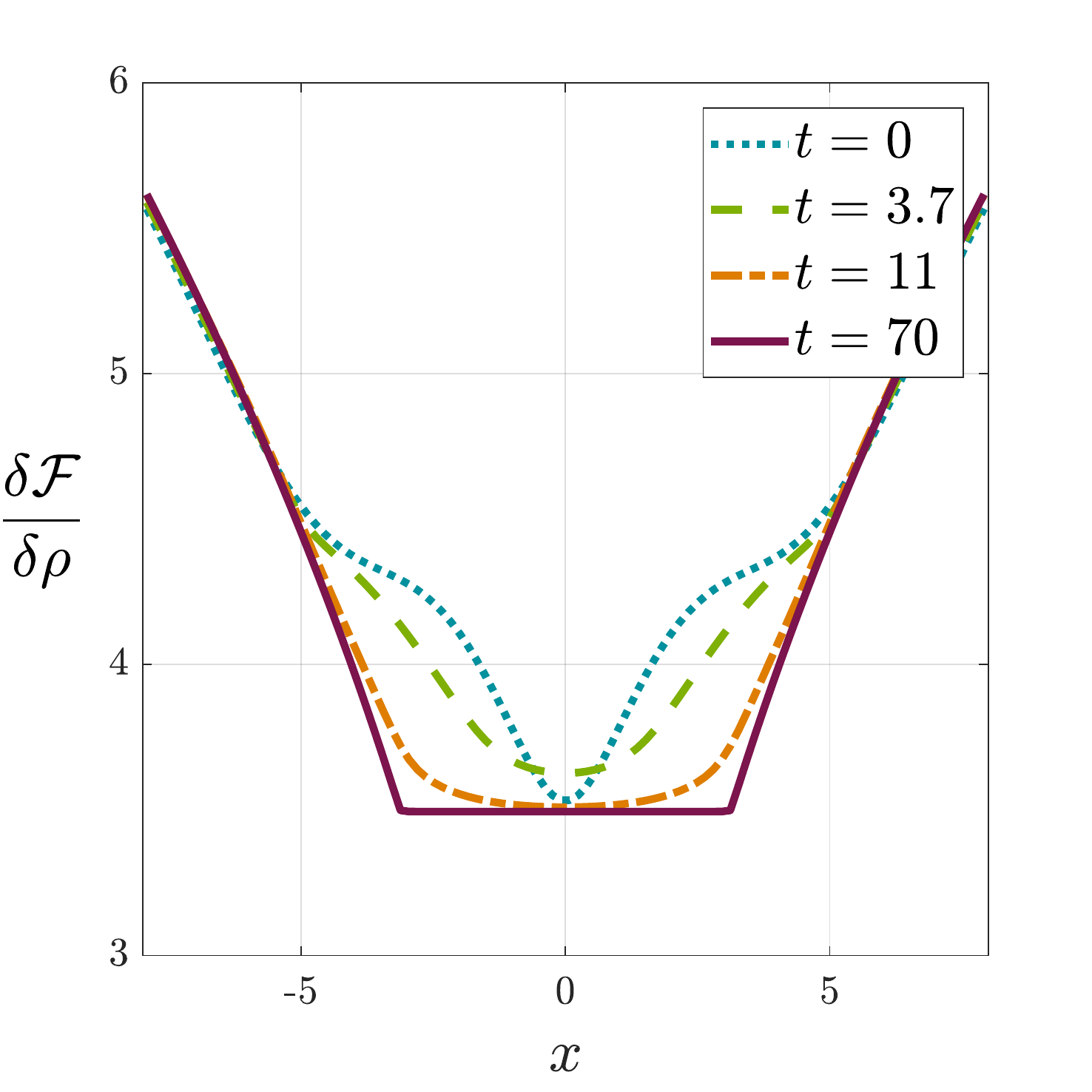}
}
\subfloat[Evolution of the total energy and free energy]{\protect\protect\includegraphics[scale=0.4]{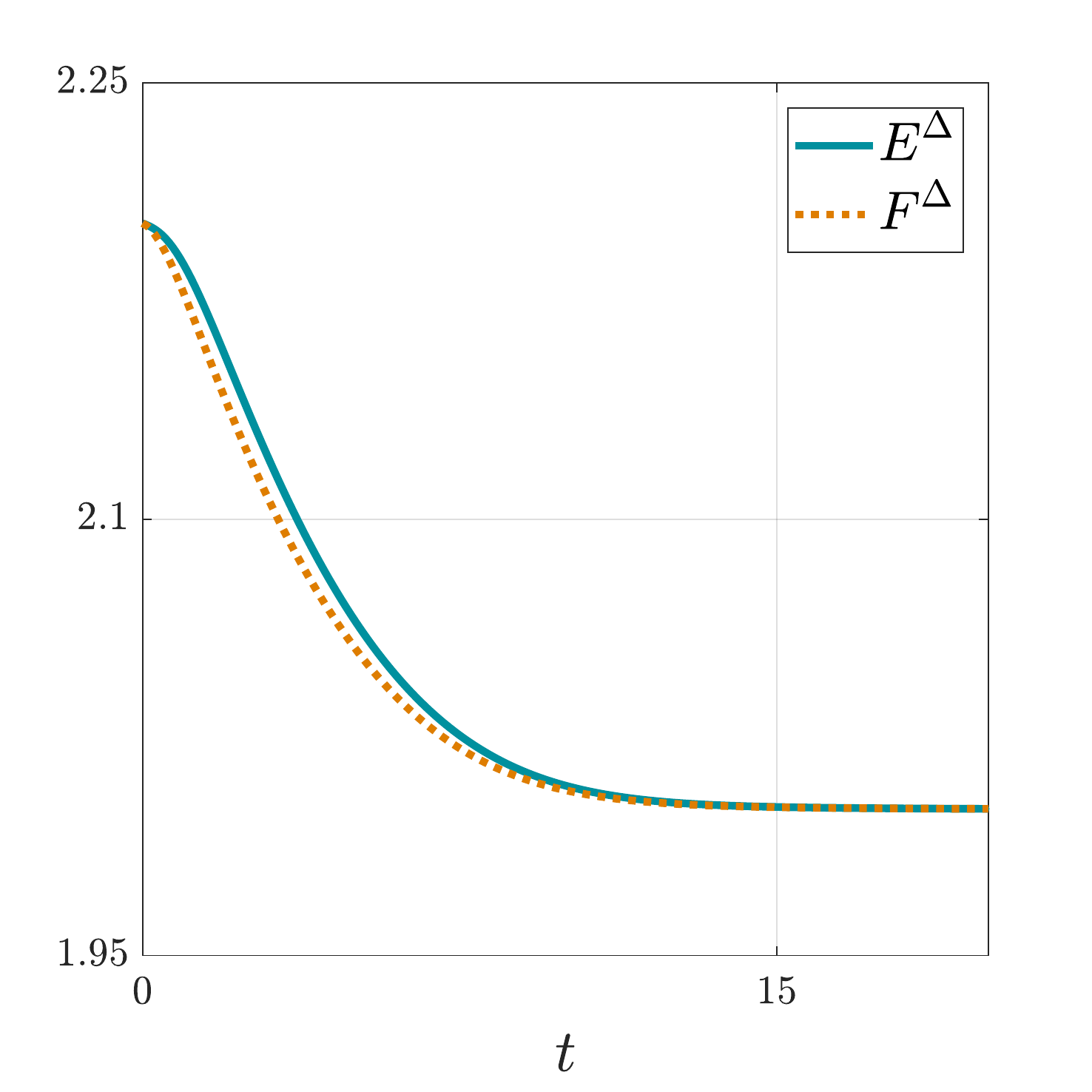}
}
\end{center}
\protect\protect\caption{\label{fig:kellersegel1} Temporal evolution of Example \ref{ex:kellersegel} with compactly-supported steady state.}
\end{figure}

We have conducted two simulations with different choices of the paramenters
$\alpha$ and $m$. In both $m>1$, so that a proper numerical flux able to deal
with vacuum regions has to be implemented. As emphasised in the introduction
of this section, the kinetic scheme developed in \cite{perthame2002kinetic}
is employed. Both of the simulations share the same initial conditions,
\begin{equation*}
\rho(x,t=0)=M_0\frac{e^{-\frac{4(x+2)^2}{10}}+e^{-\frac{4(x-2)^2}{10}}}{\int_\R \left(e^{-\frac{4(x-2)^2}{10}}+e^{-\frac{4(x+2)^2}{10}}\right)dx},\quad \rho u(x,t=0)=0,\quad x\in [-8,8],
\end{equation*}
where the total mass $M_0$ of the system is $1$.

In the first simulation the choice of parameters is $\alpha=0.5$ and $m=1.5$.
According to the regime classification for the overdamped system, this would
correspond to the diffusion-dominated regime. In the overdamped limit,
solutions exist globally in time, and the steady state is compactly
supported. The results are depicted in figure \ref{fig:kellersegel1} and
adequately agree with this regime. In the steady state the variation of the
free energy with respect to density has a constant value only in the support
of the density, as expected. The total energy decreases in time and there is no exchange between the free energy and the kinetic energy
since the free energy in figure \ref{fig:kellersegel1} (D) does not
oscillate.

\begin{figure}[ht!]
\begin{center}
\subfloat[Evolution of the density]{\protect\protect\includegraphics[scale=0.4]{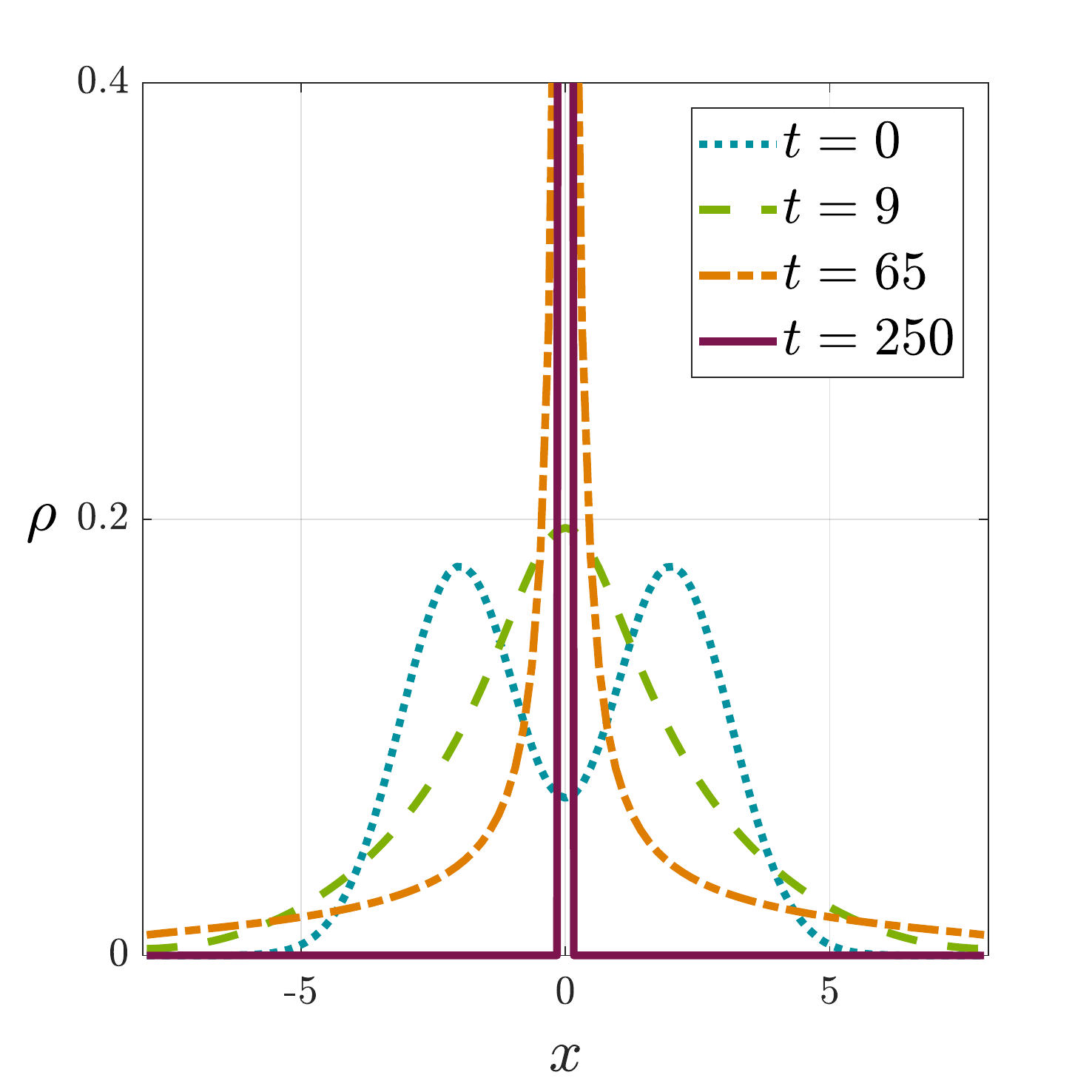}
\llap{\shortstack{%
        \includegraphics[scale=.14,trim={0 0 1.4cm 0cm},clip]{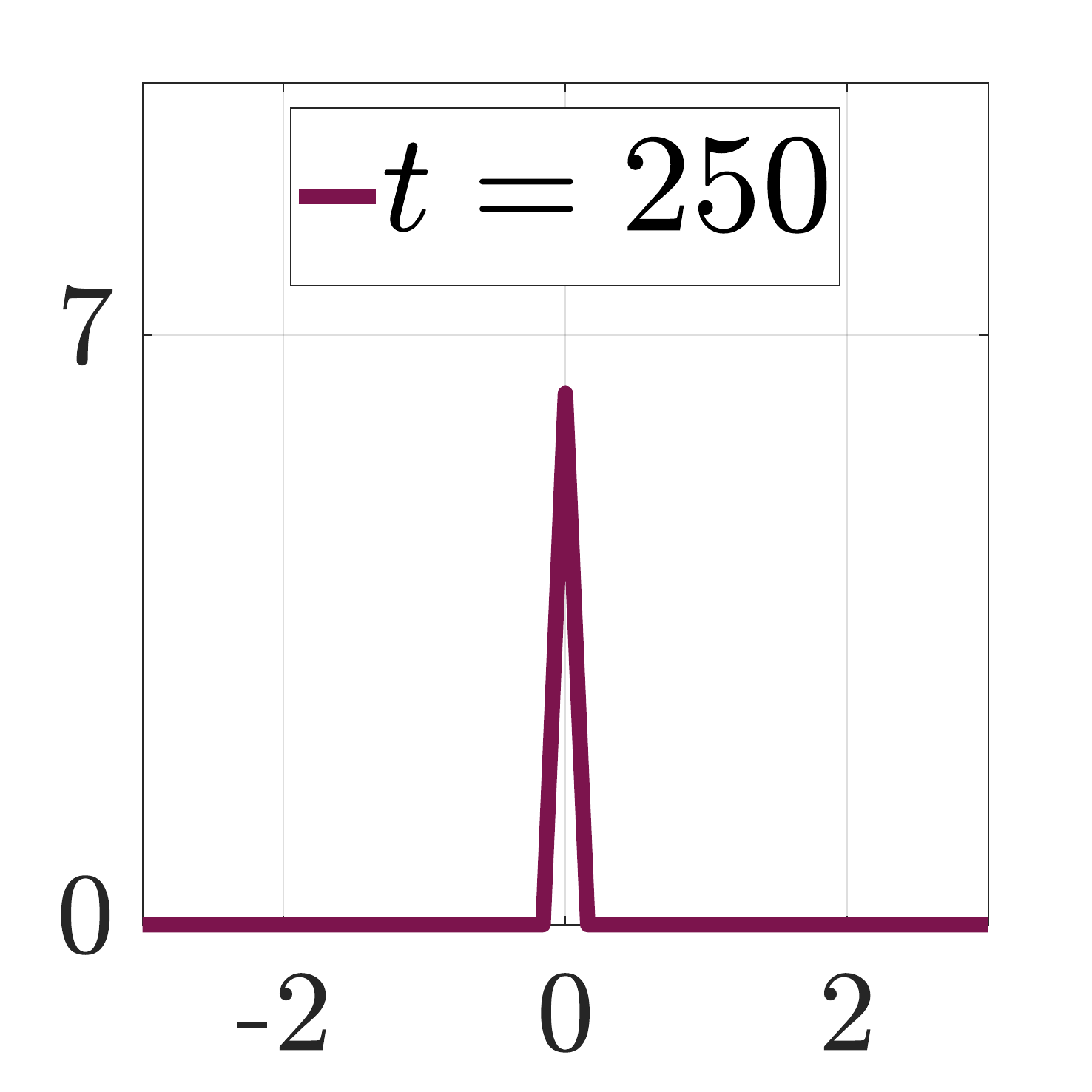}\\
        \rule{0ex}{1.28in}%
      }
  \rule{1.288in}{0ex}}

}
\subfloat[Evolution of the momentum]{\protect\protect\includegraphics[scale=0.4]{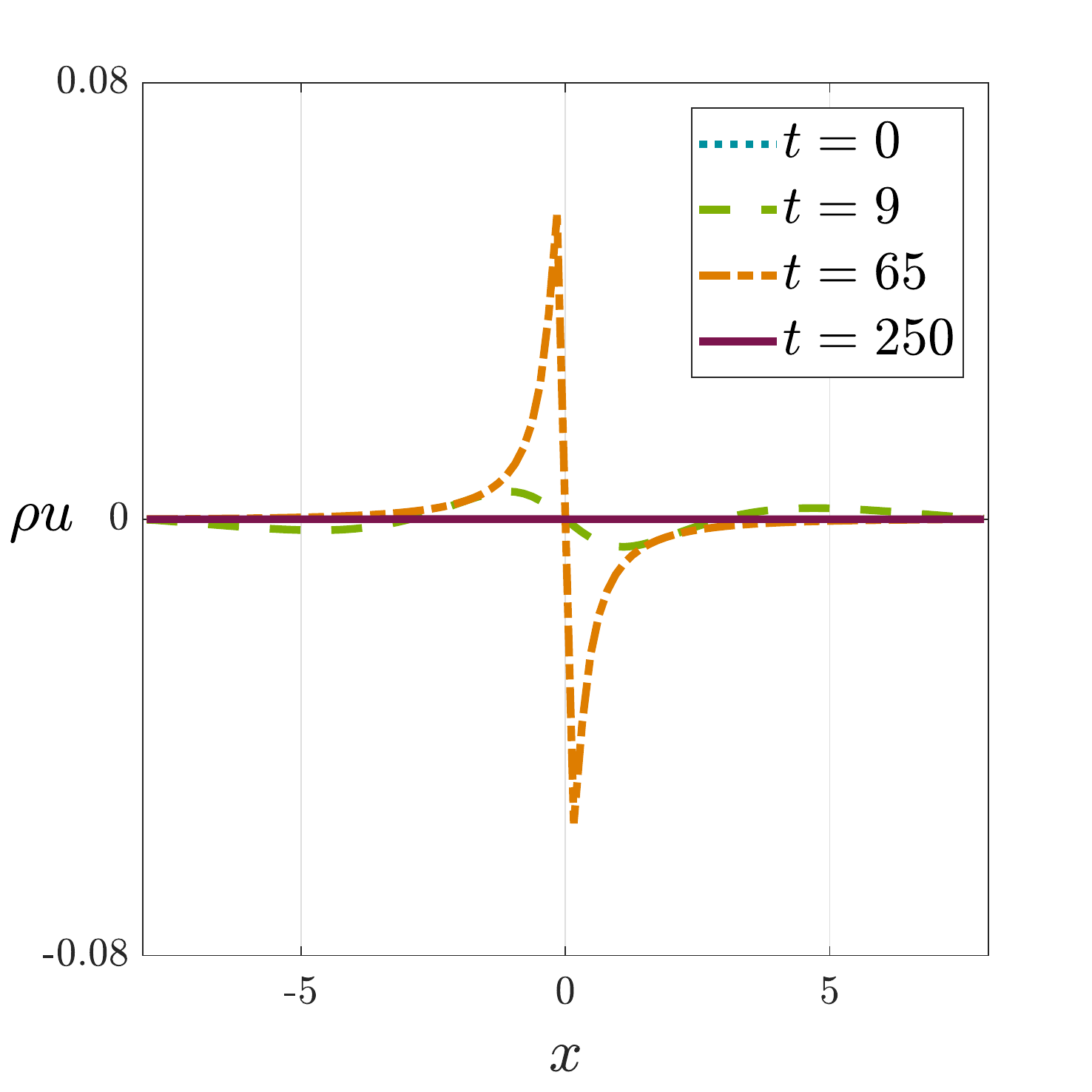}
}\\
\subfloat[Evolution of the variation of the free energy]{\protect\protect\includegraphics[scale=0.4]{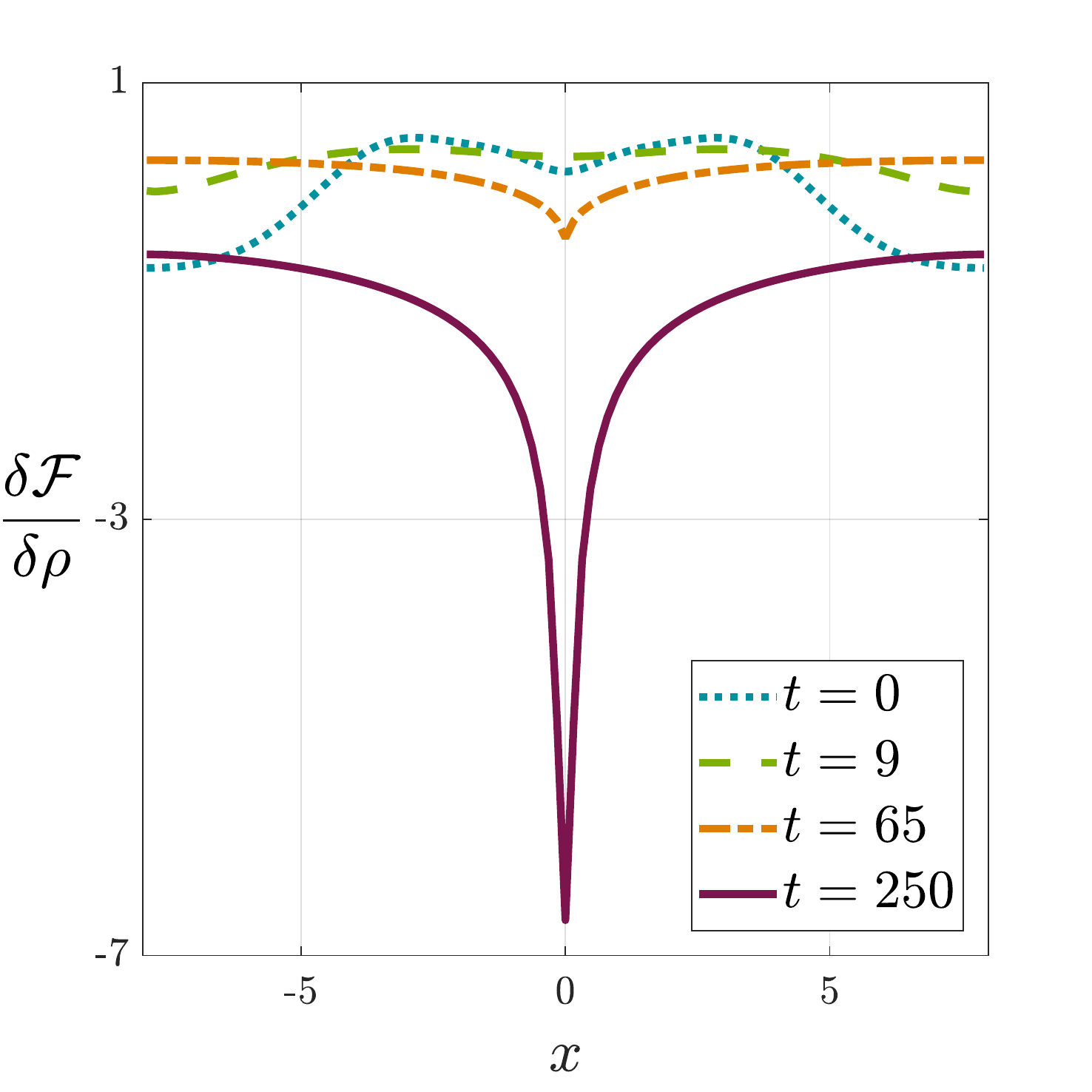}
}
\subfloat[Evolution of the total energy and free energy]{\protect\protect\includegraphics[scale=0.4]{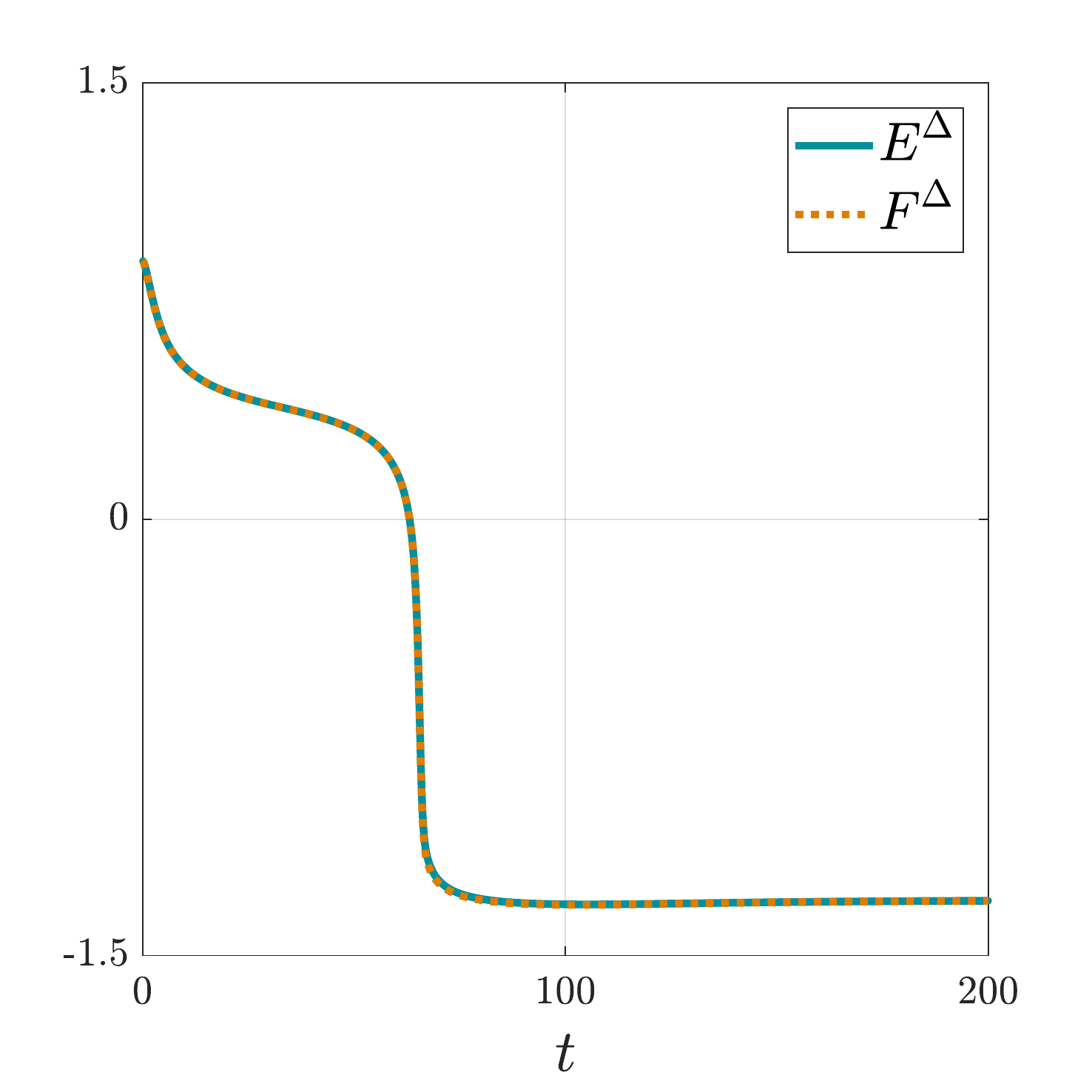}
}
\end{center}
\protect\protect\caption{\label{fig:kellersegel2} Temporal evolution of Example \ref{ex:kellersegel} with finite-time blow up.}
\end{figure}
\end{examplecase}

The second simulation has a choice of parameters of $\alpha=-0.5$ and
$m=1.3$. In the case of the overdamped system this would correspond to the
aggregation-dominated regime, where blow-up and diffusive behaviour coexist
and depend on the initial density profile. The results from this simulation
of the hydrodynamic system are illustrated in figure \ref{fig:kellersegel2}.
For this particular initial condition there is analytically finite-time blow up. Our scheme,
due to the conservation of mass of the finite volume scheme, concentrates all the mass in one single cell in finite time, that is, the scheme achieves in finite time the better approximation to a Dirac Delta at a point with the chosen mesh. Once this happens, this artificial numerical steady state depending on the mesh is kept for all times. From figure \ref{fig:kellersegel2} (C) it is evident that the variation of the
free energy with respect to density does not reach a constant value, and in
figure \ref{fig:kellersegel2} (D) the free energy presents a sharp decay when
the concentration in one cell is produced (around $t\approx65$). The value of
the slope in the free energy plot theoretically tends to $-\infty$ due to the
blow up, but in the simulation the decay is halted due to conservation of
mass and the artificial steady state. This agrees with the fact that the
expected Dirac delta profile in the density at the blow up time is obviously not
reached numerically. It was also checked that this phenomena repeats for all meshes
leading to more concentrated artificial steady states with more negative free energy values
for more refined meshes.
For other more spreaded initial conditions our scheme produces diffusive
behaviour as expected from theoretical considerations.

\begin{figure}[ht!]
\begin{center}
\subfloat[Evolution of the density]{\protect\protect\includegraphics[scale=0.4]{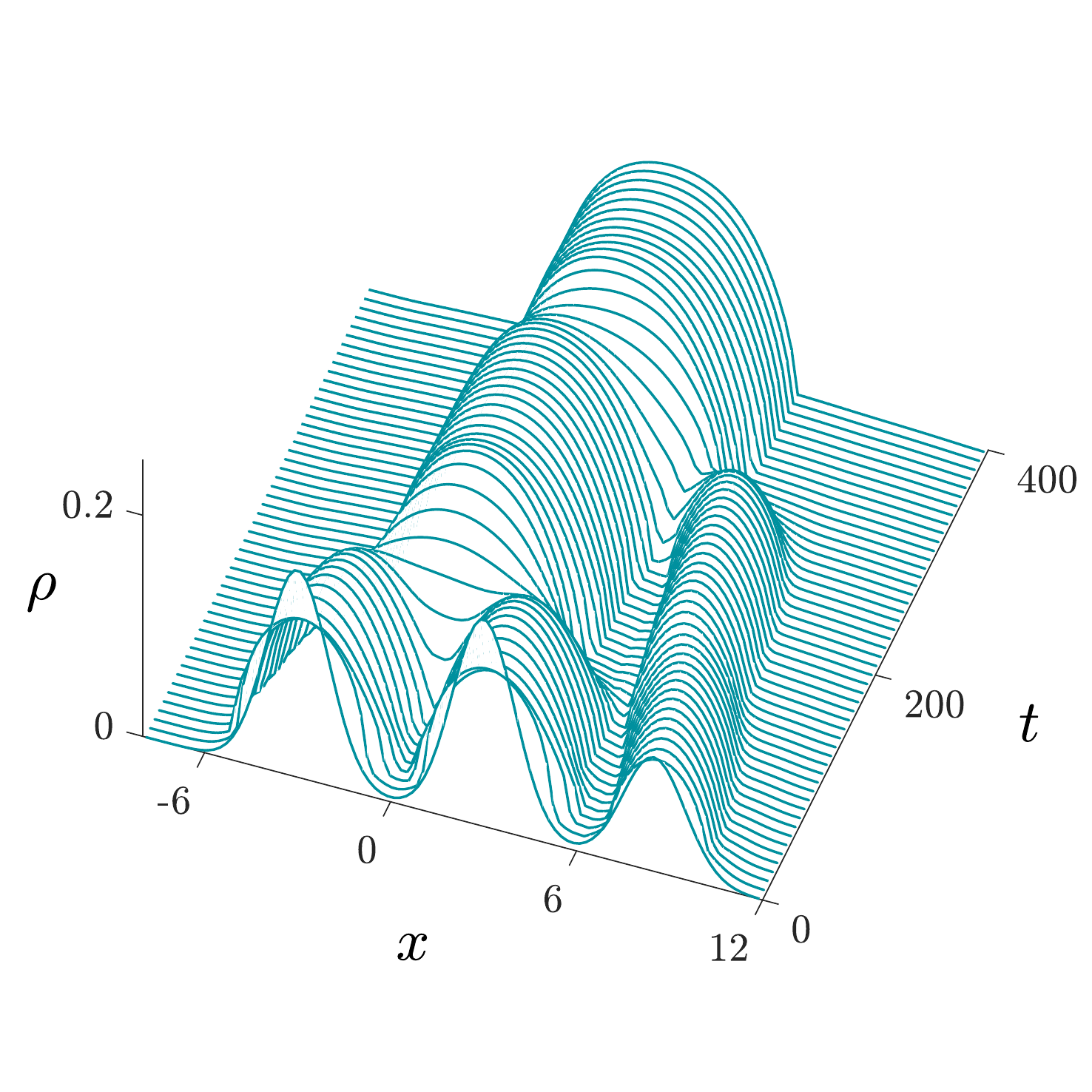}
}
\subfloat[Evolution of the momentum]{\protect\protect\includegraphics[scale=0.4]{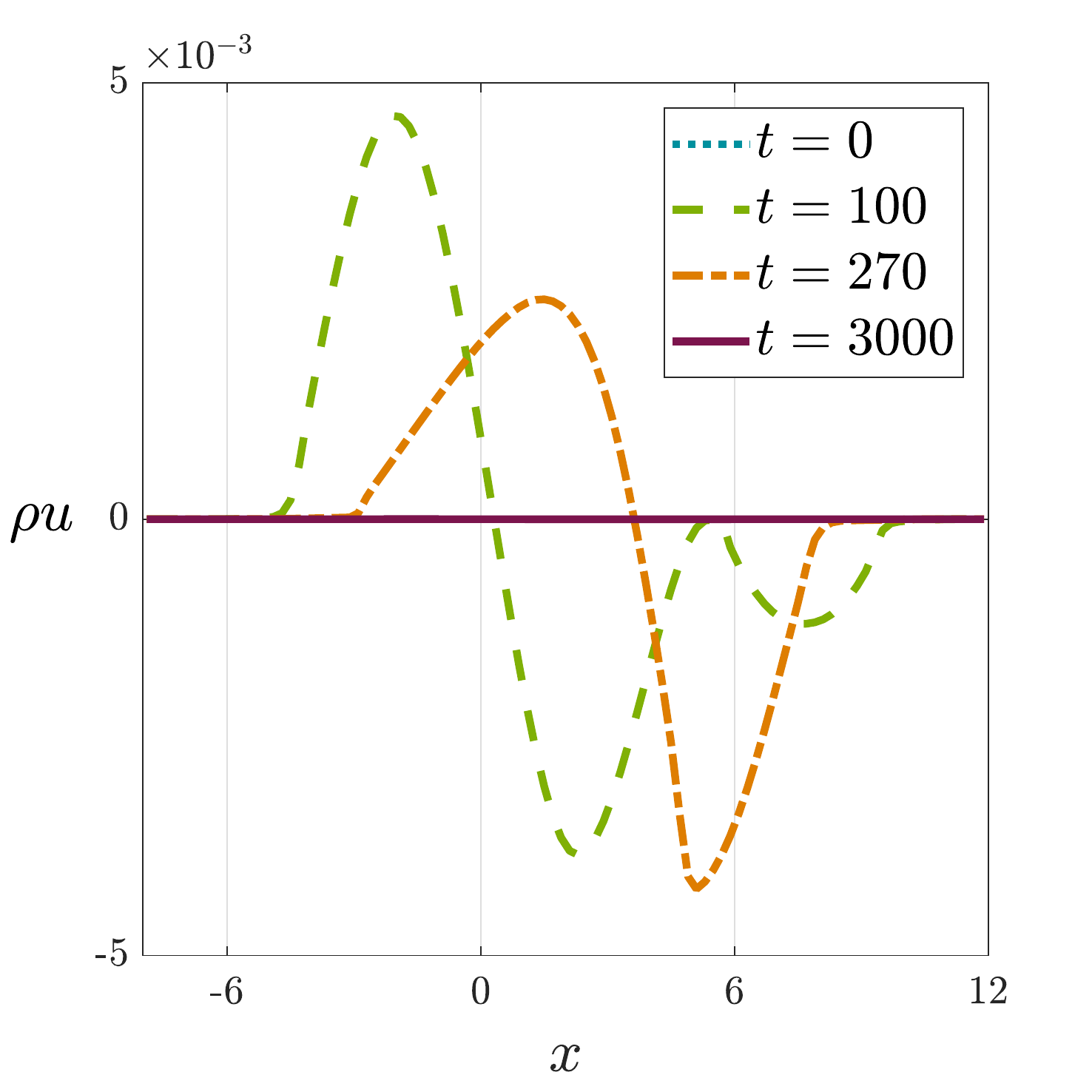}
}\\
\subfloat[Evolution of the variation of the free energy]{\protect\protect\includegraphics[scale=0.4]{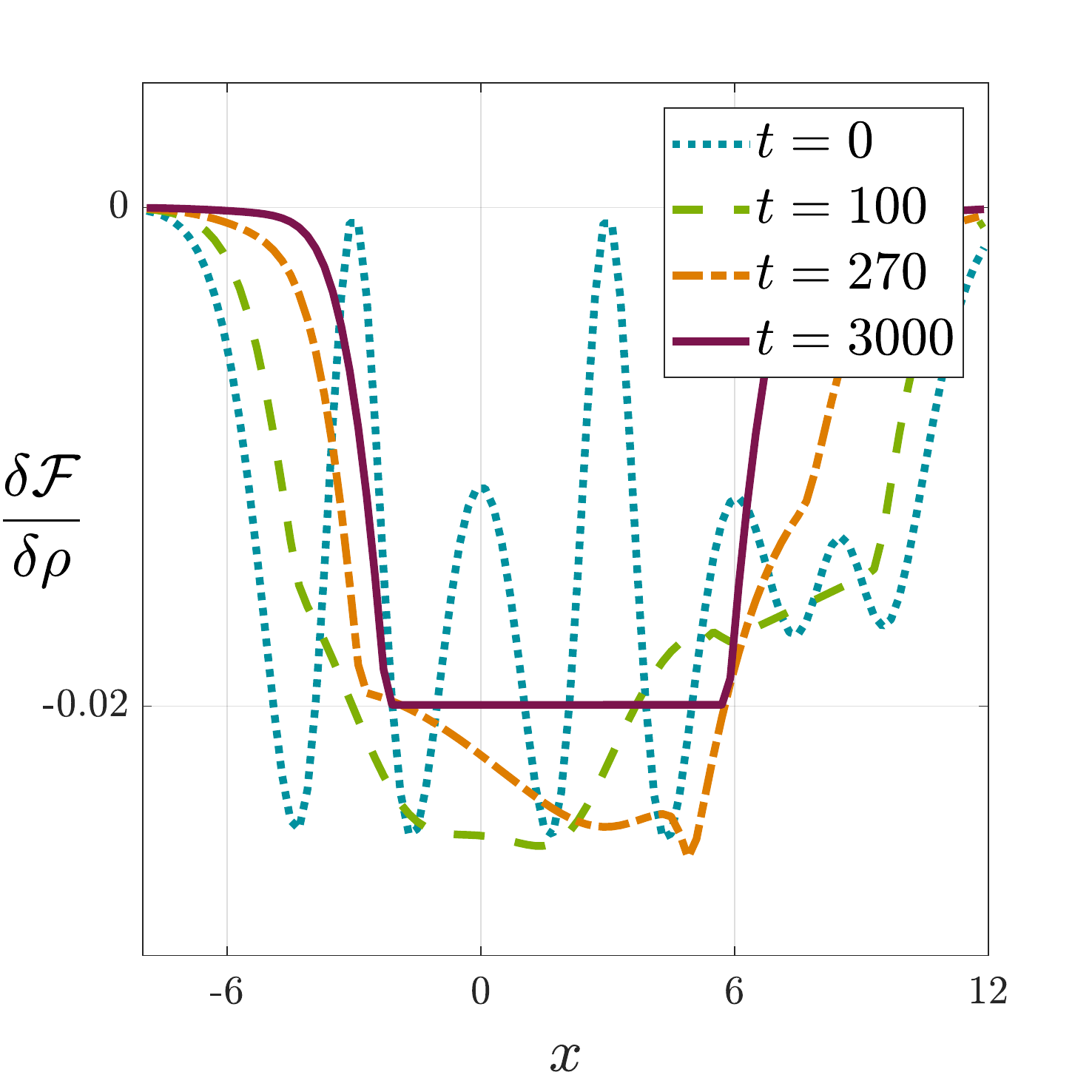}
}
\subfloat[Evolution of the total energy and free energy]{\protect\protect\includegraphics[scale=0.4]{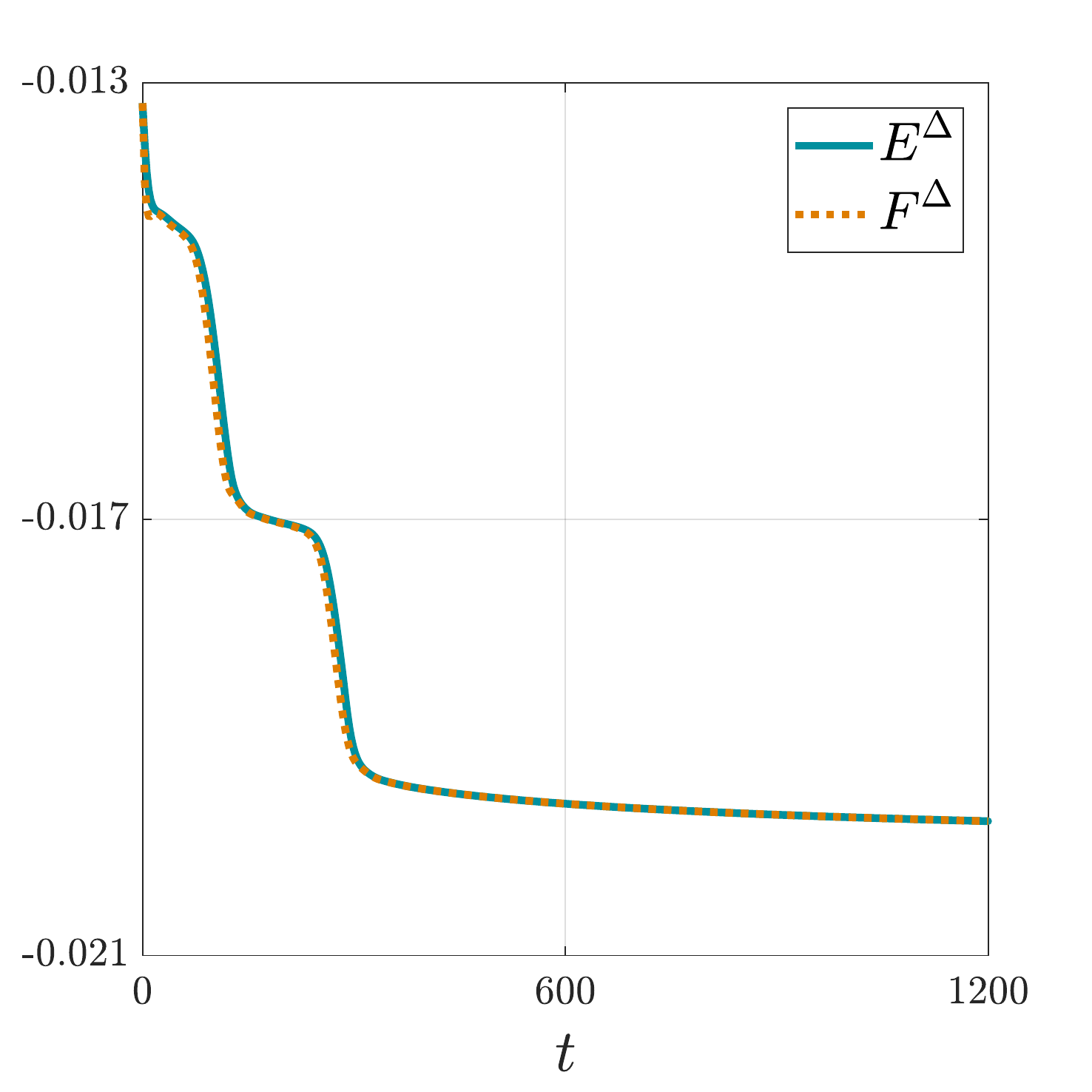}
}
\end{center}
\protect\protect\caption{\label{fig:kellersegel3} Temporal evolution of Example \ref{ex:kellersegel} with Morse-type potential and three initial density bumps.}
\end{figure}

A further simulation is carried out to explore the convergence in time towards equilibration with a Morse-type potential of the form $W(x)=-e^{-|x|^2/2}/\sqrt{2\pi}$. With this potential the attraction between two bumps of density separated at a considerable distance is quite small. However, when enough time has passed and the bumps get closer, they merge in an exponentially fast pace due to the convexity of the Gaussian potential, and a new equilibrium is reached with just one bump. The interesting fact about this system is therefore the existence of two timescales: the time to get the bumps of density close enough, which could be arbitrarily slow, and the time to merge the bumps, which is exponentially fast in time.

We have set up a simulation whose initial state presents three bumps of density, with the initial conditions satisfying
\begin{equation*}
\rho(x,t=0)=M_0\frac{e^{-\frac{(x+3)^2}{2}}+e^{-\frac{(x-3)^2}{2}}+0.55e^{-\frac{(x-8.5)^2}{2}}}{\int_\R \left(e^{-\frac{(x+3)^2}{2}}+e^{-\frac{(x-3)^2}{2}}+0.55e^{-\frac{(x-8.5)^2}{2}}\right)dx},\quad \rho u(x,t=0)=0,\quad x\in [-8,12],
\end{equation*}
and the total mass of the system equal to $M_0=1.2$. The parameter $m$ in the pressure satisfies $m=3$, and the effect of the linear damping is reduced by assigning $\gamma=0.05$.

The results are depicted in \ref{fig:kellersegel3}. In (A) one can observe
how the two central bumps of density merge after some time, and how the third
bump, with less mass, starts getting closer in time until it also blends.
This is also reflected in the evolution of the free energy in figure
\ref{fig:kellersegel3} (D), where there are two sharp and exponential decays
corresponding to the merges of the bumps.


\begin{examplecase}[DDFT for 1D hard rods]\label{ex:hardrods} Classical (D)DFT is a theoretical framework provided by
nonequilibrium statistical mechanics but has increasingly become a
widely-employed method for the computational scrutiny of the microscopic
structure of both uniform and non-uniform fluids
\cite{duran2016dynamical,goddard2012unification,yatsyshin2012spectral,yatsyshin2013geometry,lutsko2010recent}.
The DDFT equations have the same form as in \eqref{eq:generalsys2} when the
hydrodynamic interactions are neglected. The starting point in (D)DFT is a
functional $\mathcal{F}[\rho]$ for the fluid's free energy which encodes all
microscopic information such as the ideal-gas part, short-range repulsive
effects induced by molecular packing, attractive interactions and external
fields. This functional can be exactly derived only for a limited number of
applications, for instance the one-dimensional hard rod system from Percus
\cite{percus1976equilibrium}. However, in general it has to be approximated
by making appropriate assumptions, as e.g. in the so-called
fundamental-measure theory of Rosenfeld \cite{rosenfeld1989free}. These
assumptions are usually validated by carrying out appropriate test
simulations (e.g. of the underlying stochastic dynamics) to compare e.g. the
DDFT system with the approximate free-energy functional to the microscopic
reference system~\cite{goddard2012general}.

The objective of this example is to show that the numerical scheme in section
\ref{sec:numsch} can also be applied to the physical free-energy functionals
employed in (D)DFT, which satisfy the more complex expression for the free
energy described in \eqref{eq:freeenergygeneral}, and with a variation
satisfying \eqref{eq:varfreeenergygeneral}. For this example the focus is on
the hard rods system in one dimension. Its free energy has a part depending
on the local density and which satisfies the classical form for an ideal gas,
with $P(\rho)=\rho$. It is therefore usually denoted as the ideal part of the
free energy,
\begin{equation*}
\mathcal{F}_{id}[\rho]=\int \Pi(\rho)dx=\int \rho(x) \left( \ln \rho -1 \right) dx.
\end{equation*}

There is also a part of general free energy in \eqref{eq:freeenergygeneral}
which contains the non-local dependence of the density, and has different
exact or approximative forms depending of the system under consideration. In
(D)DFT it is denoted as the excessive free energy, and for the hard rods
satisfies
\begin{align*}
\mathcal{F}_{ex}[\rho]&=\frac{1}{2}\int K\left(W(\bm{x})\star \rho(\bm{x})\right) \rho(\bm{x})d\bm{x}\\
&=-\frac{1}{2}\int  \rho(x+\sigma/2)\ln\left(1-\eta(x)\right)dx -\frac{1}{2}\int  \rho(x-\sigma/2)\ln\left(1-\eta(x)\right)dx,
\end{align*}
\begin{figure}[ht!]
\begin{center}
\subfloat[Evolution of the density]{\protect\protect\includegraphics[scale=0.4]{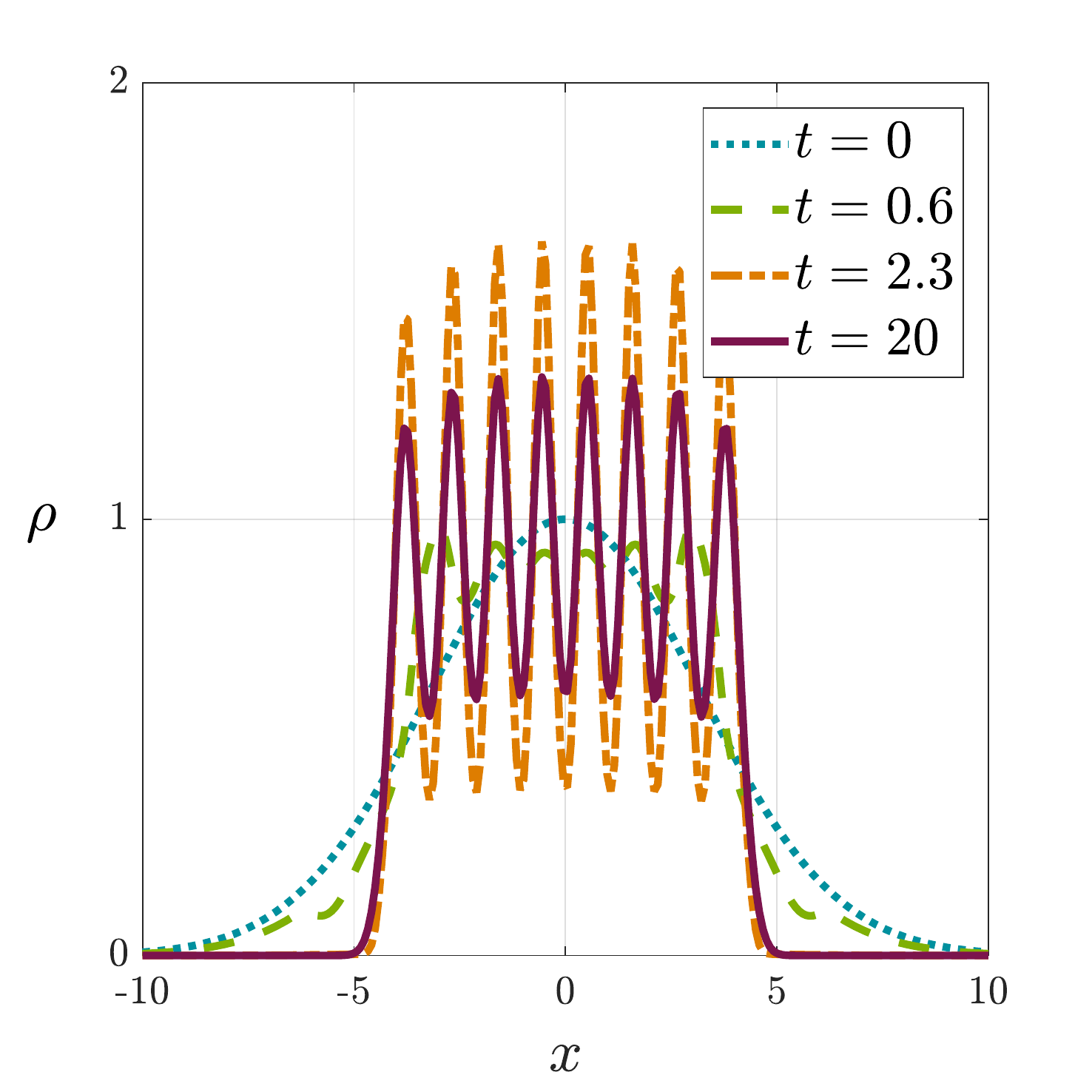}
}
\subfloat[Evolution of the momentum]{\protect\protect\includegraphics[scale=0.4]{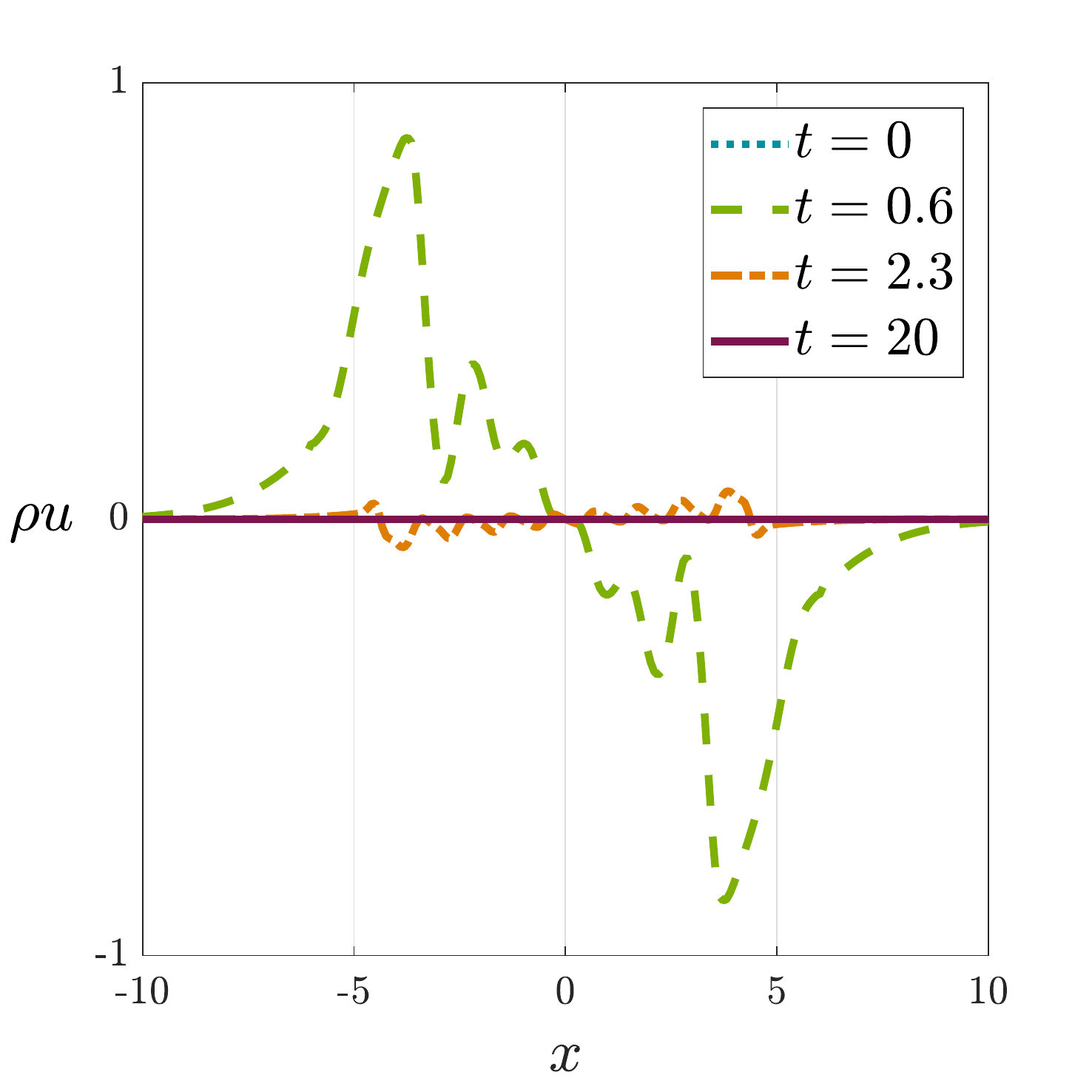}
}\\
\subfloat[Evolution of the variation of the free energy]{\protect\protect\includegraphics[scale=0.4]{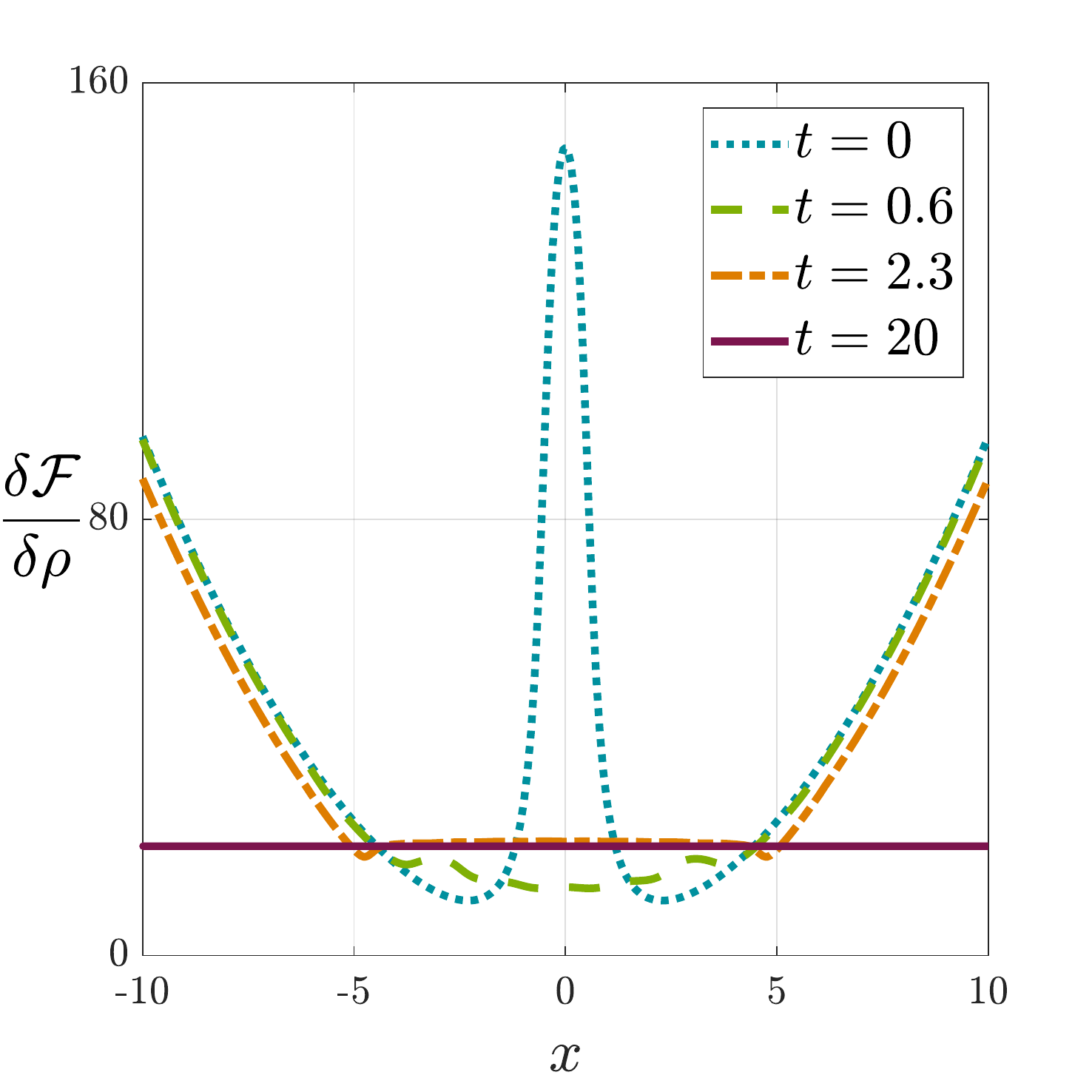}
}
\subfloat[Evolution of the total energy and free energy]{\protect\protect\includegraphics[scale=0.4]{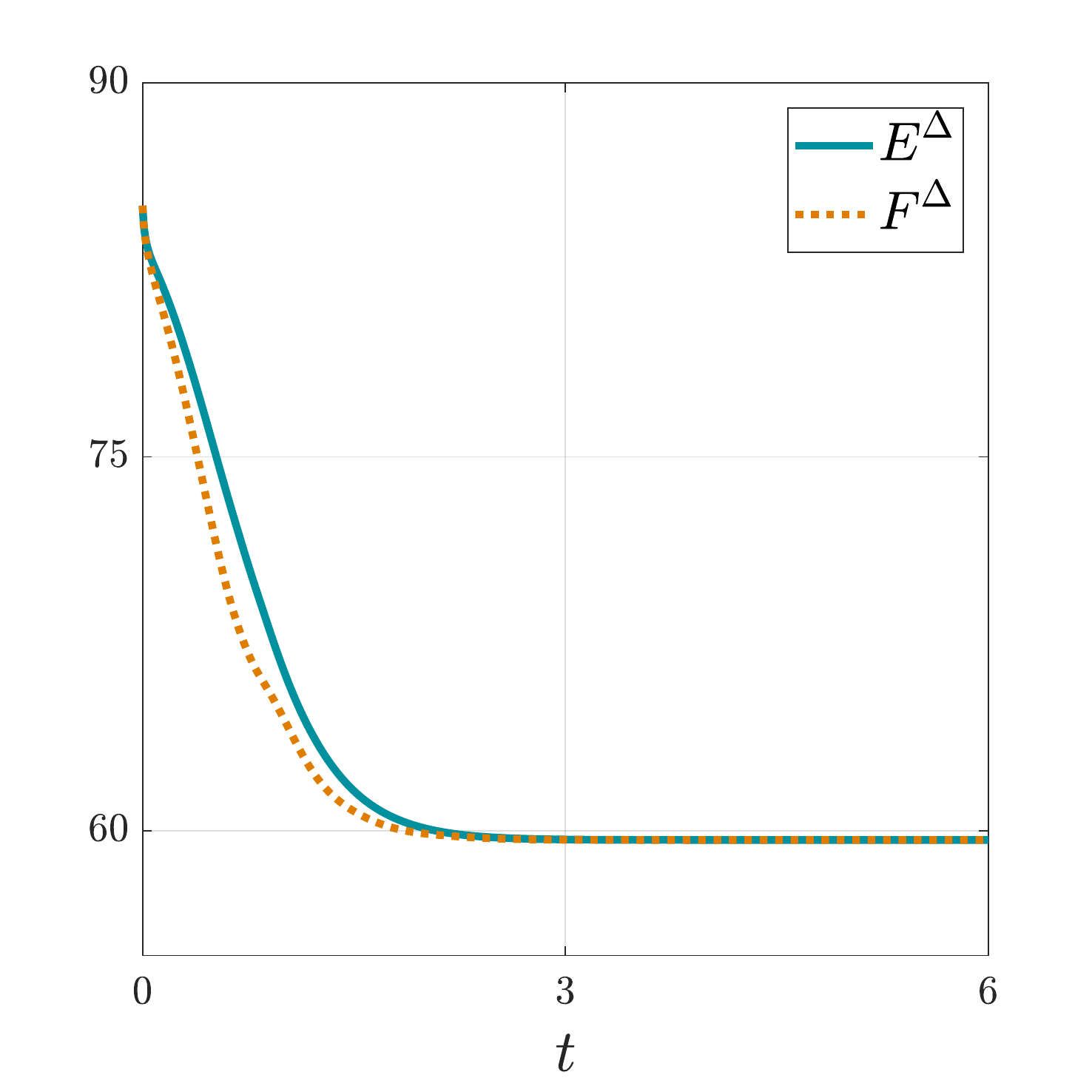}
}
\end{center}
\protect\protect\caption{\label{fig:hardrods1} Temporal evolution of Example \ref{ex:hardrods} with $8$ hard rods and a confining potential.}
\end{figure}
where $\sigma$ is the length of a hard rod and $\eta(x)$ the local packing fraction representing the probability that a point $x$ is covered by a hard rod,
\begin{equation*}\label{eq:packing}
\eta(x)=\int_{-\frac{\sigma}{2}}^{\frac{\sigma}{2}}\rho(x+y) dy .
\end{equation*}
The function $K(x)$ in this case satisfies $K(x)=\ln(1-x)$ and the kernel $W(x)$ takes the form of a characteristic function which limits the interval of the packing function \eqref{eq:packing}. To obtain the excessive free energy for the hard rods one has to also consider changes of variables in the integrals. The last part of the general free energy in \eqref{eq:freeenergygeneral} corresponds to the effect of the external potential $V(x)$. On the whole, the variation of the free energy in  \eqref{eq:freeenergygeneral} with respect to the density, for the case of hard rods, satisfies
\begin{align*}
\frac{\delta \mathcal{F}[\rho]}{\delta \rho}=&\frac{\delta \mathcal{F}_{id}[\rho]}{\delta \rho}+\frac{\delta \mathcal{F}_{ex}[\rho]}{\delta \rho}+V(x)\\
=&\ln(\rho)-\frac{1}{2}\ln\left(1-\int_{x-\sigma}^{x}\rho(y)dy\right)-\frac{1}{2}\ln\left(1-\int_{x}^{x+\sigma}\rho(y)dy\right)\\
& +\frac{1}{2}\int_{x-\sigma/2}^{x+\sigma/2}\left(\frac{\rho(x+\sigma/2)+\rho(x-\sigma/2)}{1-\eta(x)}\right)dx+V(x).
\end{align*}

This system can be straightforwardly simulated with the well-balanced scheme from section \ref{sec:numsch} by gathering the excessive part of the free energy and the external potentials under the term $H(x,\rho)$, so that
\begin{equation*}
H(x,\rho)=\frac{\delta \mathcal{F}_{ex}[\rho]}{\delta \rho}+V(x).
\end{equation*}
The first simulation seeks to capture the steady state reached by $8$ hard rods of unitary mass and length $\sigma=1$ under the presence of an external potential of the form $V(x)=x^2$. The initial conditions of the simulation are
\begin{equation*}
\rho(x,t=0)=e^{-\frac{x^2}{20.372}},\quad \rho u(x,t=0)=0,\quad x\in [-13,13],
\end{equation*}
where the density is chosen so that the total mass of the system is $8$. The
results are plotted in figure \ref{fig:hardrods1}. The steady state reached
for the density reveals layering due to the confining effects of the external
potential and the repulsion between the hard rods. These layering effects can
be amplified by increasing the coefficient in the external potential. It is
also observed how each of the $8$ peaks has a unitary width. This is due to
the fact that the length of the hard rods $\sigma$ was taken as $1$. The
variation of the free energy with respect to the density also reaches a
constant value in all the domain. For microscopic simulations of the
underlying stochastic dynamics for similar examples we refer the reader to
\cite{goddard2012unification}.

Starting from this last steady state, the second simulation performed for
this example shows how the hard rods diffuse when the confining potential is
removed. This simulation has as initial condition the previous steady state
from figure \ref{fig:hardrods1} and the external potential is set to
$V(x)=0$.  The results are depicted in figure \ref{fig:hardrods2}, and they
share the same features of the simulations in \cite{marconi1999dynamic}. The
final steady state of the density is   uniform profile resultant from the
diffusion of the hard rods, and in this situation the variation of the free
energy with respect to the density also reaches a constant value in the
steady state, as expected.

\begin{figure}[ht!]
\begin{center}
\subfloat[Evolution of the density]{\protect\protect\includegraphics[scale=0.4]{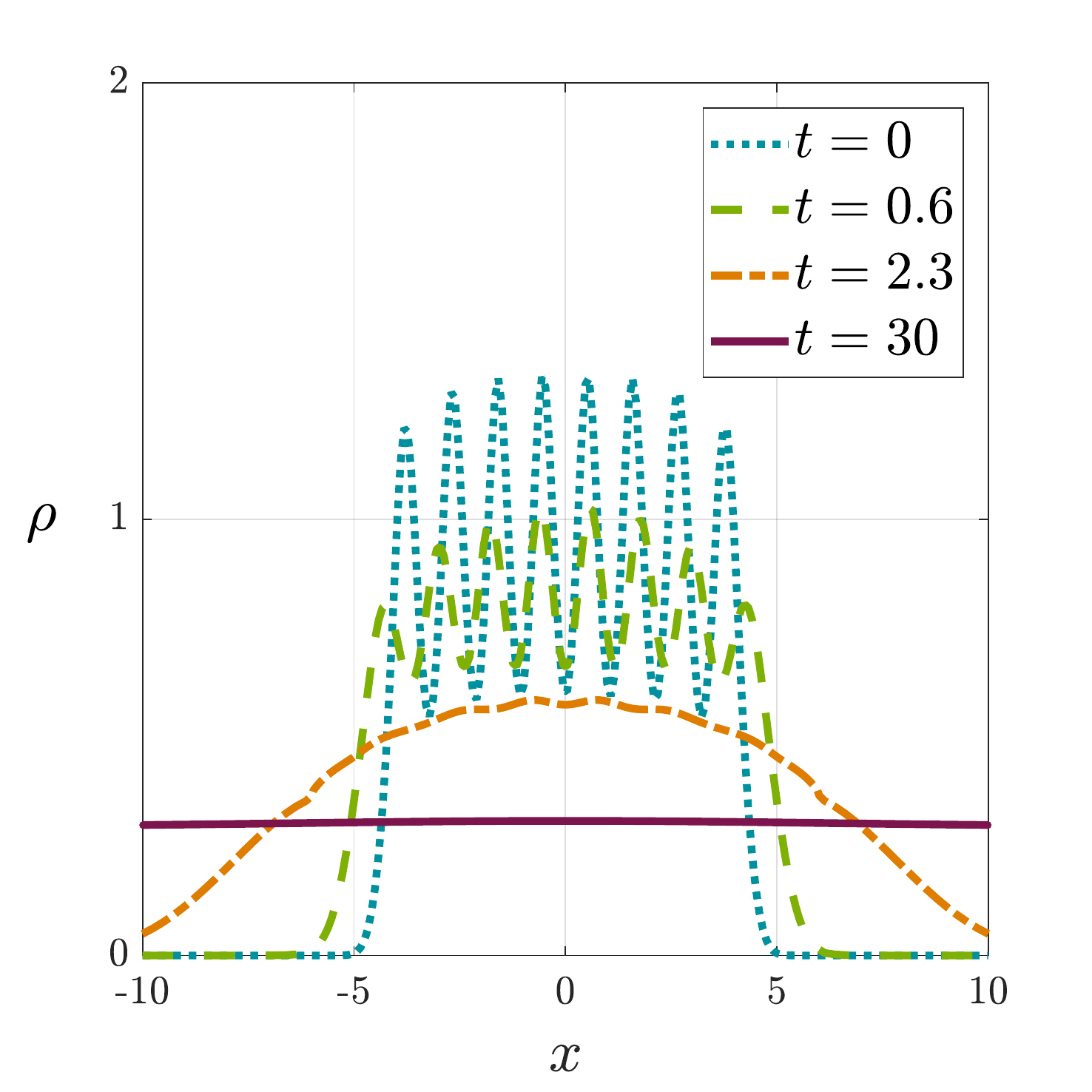}
}
\subfloat[Evolution of the momentum]{\protect\protect\includegraphics[scale=0.4]{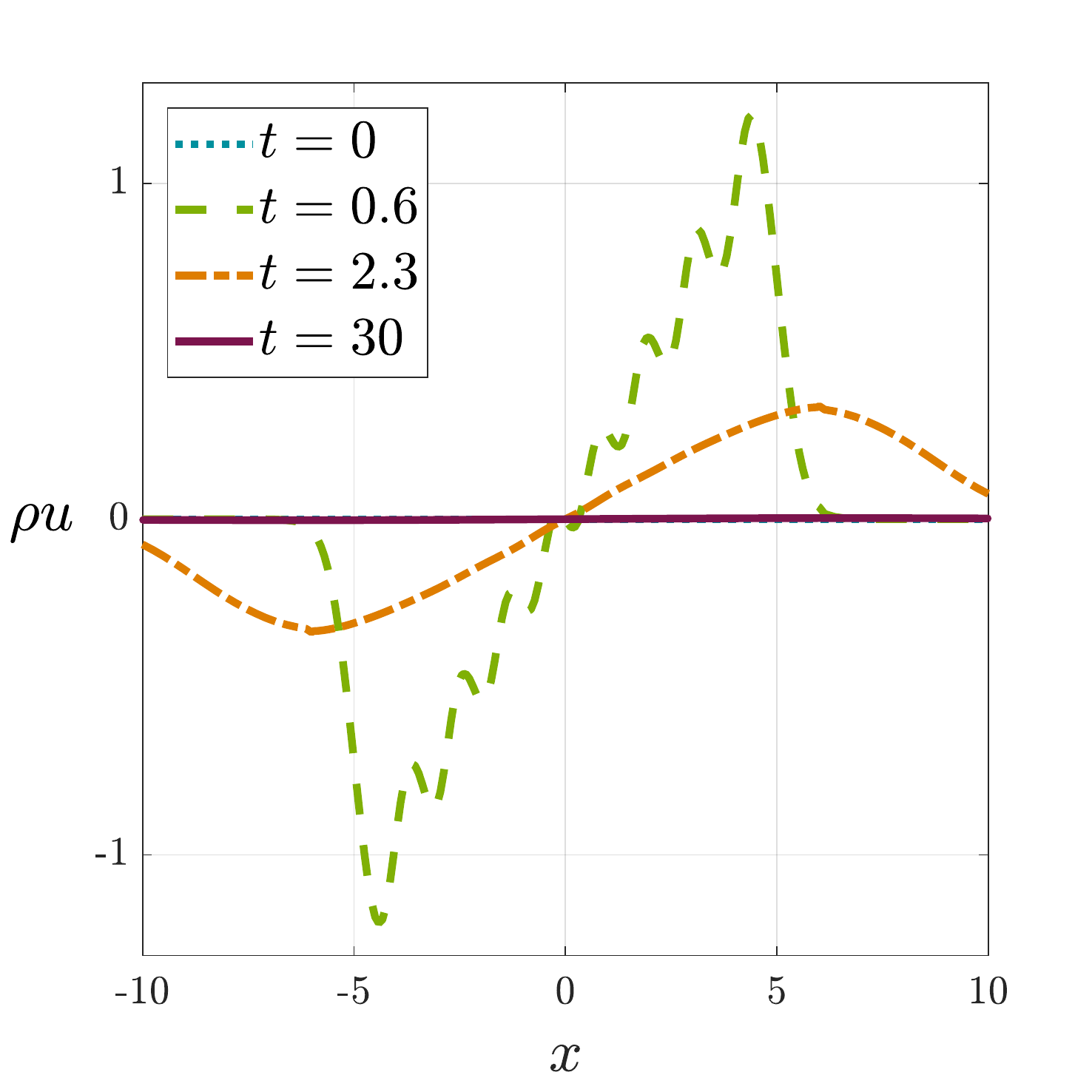}
}\\
\subfloat[Evolution of the variation of the free energy]{\protect\protect\includegraphics[scale=0.4]{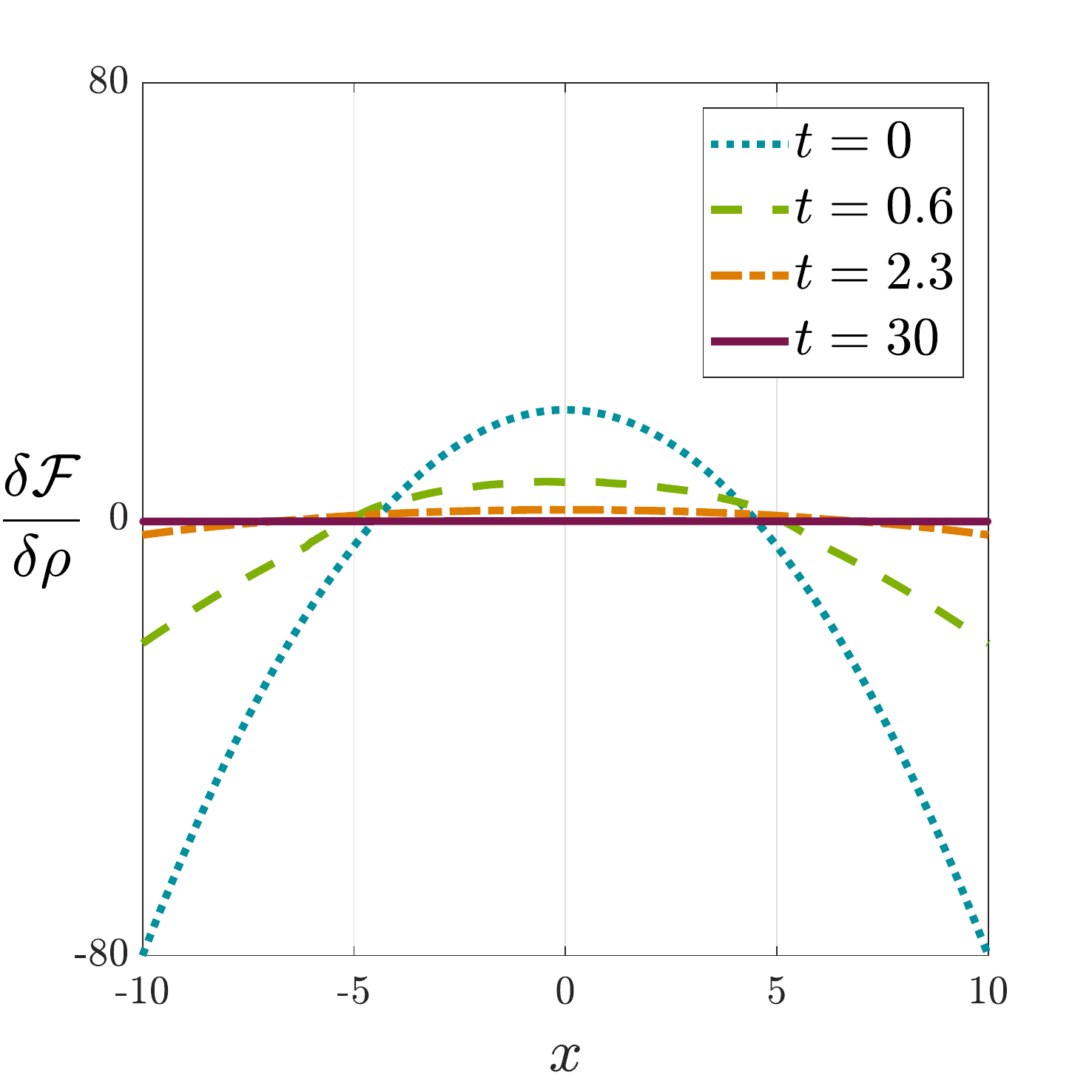}
}
\subfloat[Evolution of the total energy and free energy]{\protect\protect\includegraphics[scale=0.4]{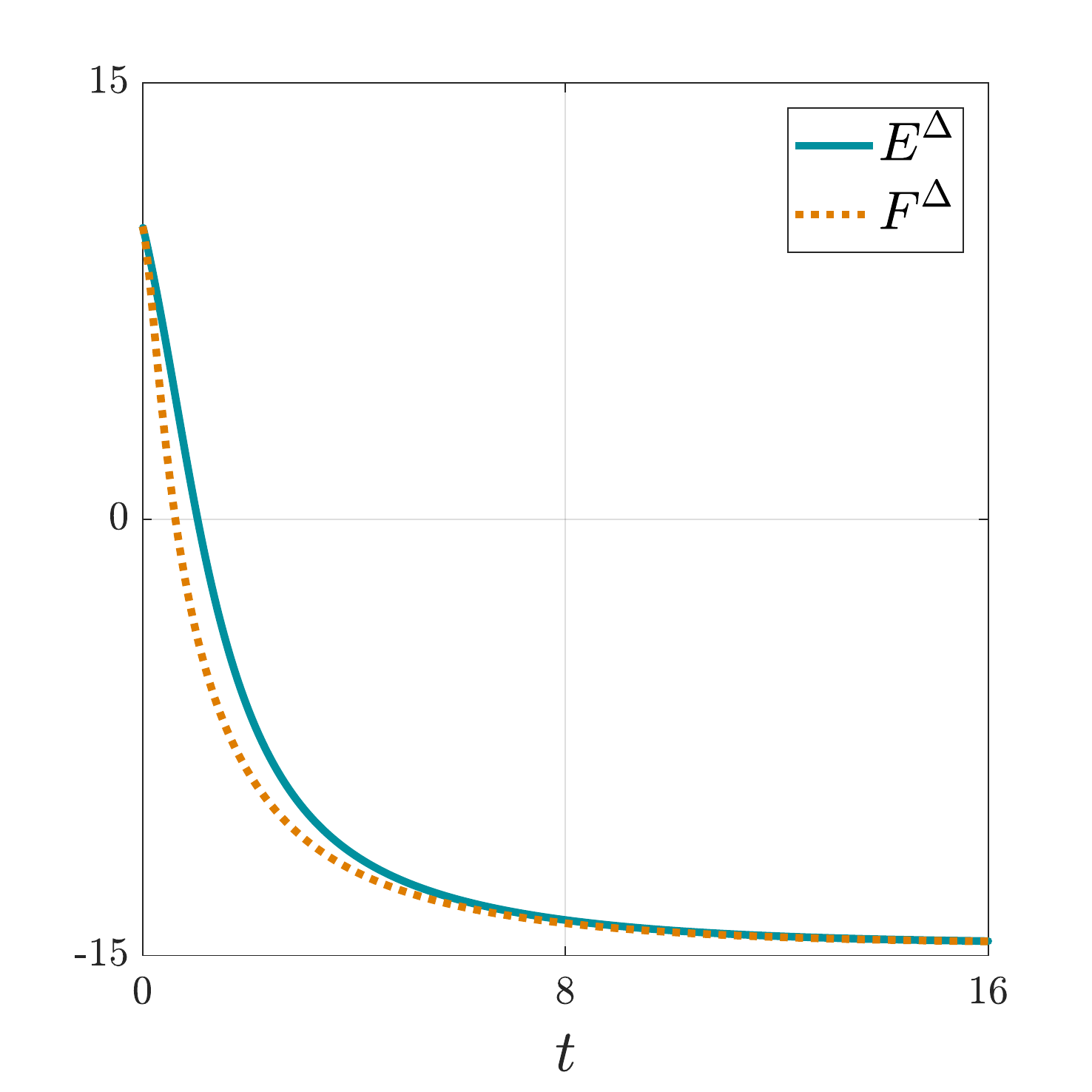}
}
\end{center}
\protect\protect\caption{\label{fig:hardrods2} Temporal evolution of Example \ref{ex:hardrods} with $8$ hard rods and no potential.}
\end{figure}

\end{examplecase}

%
%
%
%
\section{Conclusions}

We have introduced first- and second-order accurate finite volume schemes for
a large family of hydrodynamic equations with general free energy, positivity
preserving and free energy decaying properties. These hydrodynamic models
with damping naturally arise in dynamic density functional theories and the
accurate computation of their stable steady states is crucial to understand
their phase transitions and stability properties. The models possess a common
variational structure based on the physical free energy functional from
statistical mechanics. The numerical schemes proposed capture very well
steady states and their equilibration dynamics due to the crucial free energy
decaying property resulting into well-balanced schemes. The schemes were
validated in well-known test cases and the chosen numerical experiments
corroborate these conclusions for intricate phase transitions and complicated
free energies.

There are also several new avenues of possible future directions. Indeed, we
believe the computational framework and associated methodologies presented
here can be useful for the study of bifurcations and phase transitions for
systems where the free energy is known from experiments only, either physical
or in-silico ones, and then our framework can be adopted in a ``data-driven"
approach. Of particular extension would also be extension to
multi-dimensional problems. Two-dimensional problems in particular would be
of direct relevance to surface diffusion and therefore to technological
processes in materials science and catalysis. We shall examine these and
related problems in future studies.

%
%
%
%
\section*{Acknowledgements}
We are indebted to P. Yatsyshin and M. A. Dur\'an-Olivencia from the Chemical
Engineering Department of Imperial College (IC) for numerous stimulating
discussions on statistical mechanics of classical fluids and density
functional theory. J.~A.~Carrillo was partially supported by EPSRC via Grant
Number EP/P031587/1 and acknowledges support of the IBM Visiting
Professorship of Applied Mathematics at Brown University. S.~Kalliadasis was
partially supported by EPSRC via Grant Number EP/L020564/1. S.~P.~Perez
acknowledges financial support from the IC President's PhD Scholarship and
thanks Brown University for hospitality during a visit in April 2018. C.-W.
Shu was partially supported by NSF via Grant Number DMS-1719410.

\bibliographystyle{siam}
\bibliography{CKPSbib}

\appendix

\section{Numerical flux, temporal scheme, and CFL condition employed in the numerical simulations}\label{app:numerics}
This appendix aims to present the necessary details to compute the numerical flux, \red{boundary conditions}, the CFL condition, and the temporal discretization for the simulations in section \ref{sec:numtest}.

The pressure function in the simulations has the form of $P(\rho)=\rho^m$, with $m\geq1$. When $m=0$ the pressure satisfies the ideal-gas relation $P(\rho)=\rho$, and the density does not present vacuum regions during the temporal evolution. For this case the employed numerical flux is the versatile local Lax-Friedrich flux, which approximates the flux at the boundary $F_{i+\frac{1}{2}}$ in \eqref{eq:numflux} as
\begin{equation}
F_{i+\frac{1}{2}}=\mathscr{F}\left(U_{i+\frac{1}{2}}^-,U_{i+\frac{1}{2}}^+\right)=\frac{1}{2}\left(F\left(U_{i+\frac{1}{2}}^-\right)+F\left(U_{i+\frac{1}{2}}^+\right)-\lambda_{i+\frac{1}{2}} \left(U_{i+\frac{1}{2}}^+-U_{i+\frac{1}{2}}^-\right)\right),
\end{equation}
where $\lambda$ is taken as the maximum of the absolute value of the eigenvalues of the system,
\begin{equation}
\lambda_{i+\frac{1}{2}}=\max_{\left(U_{i+\frac{1}{2}}^-,U_{i+\frac{1}{2}}^+\right)}\left\{\left|u+\sqrt{P'(\rho)}\right|,\left|u-\sqrt{P'(\rho)}\right|\right\}.
\end{equation}
This maximum is taken locally for every node, resulting in different values of $\lambda$ along the lines of nodes. It is also possible to take the maximum globally, leading to the classical Lax-Friedrich scheme.
 
For the simulations where $P(\rho)=\rho^m$ and $m>1$ vacuum regions with $\rho=0$ are generated. This implies that the hiperbolicity of the system \eqref{eq:generalsys2} is lost in those regions, and the local Lax-Friedrich scheme fails. As a result, an appropiate numerical flux has to be implemented to handle the vacuum regions. In this case a kinetic solver based on \cite{perthame2001kinetic} is employed. This solver is constructed from kinetic formalisms applied in macroscopic models, and has already been employed in previous works for shallow-water applications \cite{audusse2005well}. The flux at the boundary $F_{i+\frac{1}{2}}$ in \eqref{eq:numflux} is computed from
\begin{equation}
F_{i+\frac{1}{2}}=\mathscr{F}\left(U_{i+\frac{1}{2}}^-,U_{i+\frac{1}{2}}^+\right)=A_-\left(U_{i+\frac{1}{2}}^-\right)+A_+\left(U_{i+\frac{1}{2}}^+\right),
\end{equation}
where
\begin{equation}
A_-\left(\rho,\rho u\right)=\int_{\xi\geq0}\xi \begin{pmatrix} 1  \\ \xi \end{pmatrix} M(\rho,u-\xi)\,d\xi, \quad  A_+\left(\rho,\rho u\right)=\int_{\xi\leq0}\xi \begin{pmatrix} 1  \\ \xi \end{pmatrix} M(\rho,u-\xi)\,d\xi.
\end{equation}
The function $M(\rho,\xi)$ is chosen accordingly to the kinetic representation of the macroscopic system, and for this case satisfies
\begin{equation}
M(\rho,\xi)=\rho^{\frac{2-m}{2}} \chi\left(\frac{\xi}{\rho^\frac{m-1}{2}}\right).
\end{equation}
The function $\chi(\omega)$ can be chosen in different ways. For this simulations we simply take it as a characteristic function, 
\begin{equation}\label{eq:chikinetic}
\chi(\omega)=\frac{1}{\sqrt{12}}\mathbbm{1}_{\left\{\left|\omega\right|\leq\sqrt{3}\right\}},
\end{equation} 
although \cite{perthame2001kinetic} presents other possible choices for $\chi(\omega)$. Further valid numerical fluxes able to treat vacuum, such as the Rusanov flux or the Suliciu relaxation solver, are reviewed in \cite{bouchut2004nonlinear}.

\red{The boundary conditions are taken to be no flux both for the density and the momentum equations. As a result, the evaluation of the numerical fluxes \eqref{eq:numflux} at the boundaries of the domain is taken as
\begin{equation}
F_{i-\frac{1}{2}}=0\enspace\text{if}\:i=1 \quad \text{and}\quad  F_{i+\frac{1}{2}}=0\enspace\text{if}\:i=n.
\end{equation}}

The time discretization is acomplished by means of the third order TDV Runge-Kutta method \cite{gottlieb1998total}. From \eqref{eq:compactsys} we can define $L(U)$ as $L(U)=S(x,U)-\pa_xF(U)$, so that $\pa_t U=L(U)$. Then, the third order TDV Runge-Kutta temporal scheme to advance from $U^n$ to $U^{n+1}$ with a time step $\Delta t$ reads
\begin{align*}
U^{(1)}&=U^n+\Delta t L\left(U^n\right),\\
U^{(2)}&=\frac{3}{4}U^n+\frac{1}{4}U^{(1)}+\frac{1}{4}\Delta t L\left(U^{(1)}\right),\\
U^{n+1}&=\frac{1}{3}U^n+\frac{2}{3}U^{(2)}+\frac{2}{3}\Delta t L\left(U^{(2)}\right).
\end{align*} 

The time step $\Delta t$ for the case of Lax-Friedrich flux is chosen from the CFL condition,
\begin{equation}
\Delta t = \text{CFL} \frac{\min_i \Delta x_i}{\max_{\forall\left(U_{i+\frac{1}{2}}^-,U_{i+\frac{1}{2}}^+\right)}\left\{\left|u+\sqrt{P'(\rho)}\right|,\left|u-\sqrt{P'(\rho)}\right|\right\}},
\end{equation}
and the $\Delta t$ for the kinetic flux, with a function $\chi(\omega)$ as in \eqref{eq:chikinetic}, is chosen as
\begin{equation}
\Delta t = \text{CFL} \frac{\min_i \Delta x_i}{\max_{\forall\left(U_{i+\frac{1}{2}}^-,U_{i+\frac{1}{2}}^+\right)}\left\{|u|+3^{\frac{m-1}{4}}\right\}}.
\end{equation}
The CFL number is taken as $0.7$ in all the simulations.

\end{document}